# Some exponential and trigonometric integrals involving $\log \Gamma(x)$

Donal F. Connon

dconnon@btopenworld.com

4 July 2022

**Abstract**

This paper follows on from [20] and considers some integrals where the integrand comprises the log gamma function or the digamma function multiplied by exponential or trigonometric functions. Some examples are shown below:

$$\int_0^1 x \log \Gamma(x) \cot(\pi x)\, dx = \frac{1}{2\pi}\left[\gamma_1 + \frac{1}{2}[\varsigma(2) + \gamma^2]\right]$$

$$\int_0^1 x \log \Gamma(x) \sin 2\pi x\, dx = \frac{\gamma}{4\pi}$$

$$\int_0^\infty \left[\frac{\sinh(xt)}{t(e^t - 1)} - \frac{x}{te^t}\right] dt = \frac{1}{2}[\log \Gamma(1-x) - \log \Gamma(1+x)]$$

We also determine the integral referenced by Glasser [33] in 1966

$$\int_0^1 \psi\left(\frac{x}{2}\right) x(1-x) \cos \pi x\, dx = -2\log A - \frac{1}{\pi^2}\left[2 + \frac{7}{2}\varsigma(3)\right] + \frac{1}{6}[\gamma + \log \pi]$$

$$-\frac{1}{2\pi^2}\sum_{n=2}^\infty \frac{1}{n^2} \log\left(1 - \frac{1}{4n^2}\right)$$

where $\gamma_1$ is the first Stieltjes constant and $A$ is the Glaisher–Kinkelin constant.

Some new Fourier series involving the logarithmic function are also presented.

CONTENTS

## 1. Introduction

Amdeberhan, Espinosa and Moll [3] showed in 2008 that for $p > 0$

(1.1) $$\int_0^1 e^{-px} \log \Gamma(x)\, dx = \frac{[\log(2\pi)+\gamma][p-(1-e^{-p})]}{p^2} - \frac{(1-e^{-p})}{2p} \Lambda\left(\frac{p}{2\pi}\right)$$

$$+ 2(1-e^{-p}) \sum_{n=1}^{\infty} \frac{\log n}{4\pi^2 n^2 + p^2}$$

where $\Lambda(z)$ is defined by

(1.2) $$\Lambda(z) := \lim_{n \to \infty} \left( \sum_{j=1}^{n} \frac{j}{j^2 + z^2} - \log n \right)$$

Since $\gamma = \lim_{n \to \infty} \left( \sum_{j=1}^{n} \frac{1}{j} - \log n \right)$, it is easily seen that $\Lambda(0) = \gamma$, where $\gamma$ is Euler's constant.

Using L'Hôpital's rule it is easy to show that (1.1) is valid in the limit as $p \to 0$. In fact, (1.1) is also valid for negative values of $p$. Hence, (1.1) is valid for all $p \in \mathbf{R}$ and we shall consider complex values of $p$ in section 2 of this paper.

We note from (1.2) that

$$\Lambda(z) = \lim_{n \to \infty} \left( \sum_{j=1}^{n} \frac{j}{j^2 + z^2} - \log(n+1) + \log\left(1 + \frac{1}{n}\right) \right)$$

$$= \lim_{n \to \infty} \left( \sum_{j=1}^{n} \frac{j}{j^2 + z^2} - \log(n+1) \right)$$

and employing the telescoping identity

$$\log(n+1) = \sum_{j=1}^{n} \log\left(1 + \frac{1}{j}\right)$$

we obtain the series representation

(1.3) $$\Lambda(z) = \sum_{j=1}^{\infty} \left[ \frac{j}{j^2 + z^2} - \log\left(1 + \frac{1}{j}\right) \right]$$

Differentiation results in

$$\Lambda'(z) = -2z \sum_{j=1}^{\infty} \frac{j}{(j^2 + z^2)^2}$$

and we note that $\Lambda'(0) = 0$. It may also be seen, for example, that $\Lambda''(0) = -2\varsigma(3)$. In fact, we will see below in (5.1) that

(1.3.1) $$\Lambda(z) = \gamma + \sum_{n=1}^{\infty} (-1)^n \varsigma(2n+1) z^{2n}$$

As noted in Section 5 below, the $\Lambda(x)$ function in fact turns out to be rather useful.

We previously noted in [24] that the digamma function may be represented by the infrequently used formula

(1.4) $$\psi(1+x) = \sum_{j=1}^{\infty} \left[ \log\left(1 + \frac{1}{j}\right) - \frac{1}{j+x} \right]$$

and hence we have

(1.4.1) $$\psi(1+iv) + \psi(1-iv) = 2 \sum_{j=1}^{\infty} \left[ \log\left(1 + \frac{1}{j}\right) - \frac{j}{j^2 + v^2} \right]$$

Having regard to (1.3) we then see that

(1.5) $$\psi(1+iv) + \psi(1-iv) = -2\Lambda(v)$$

which was previously obtained in a different manner by Dixit [28] (see Section 5 below). Substituting this in (1.1) gives us the equivalent integral reported by Dixit [28] in 2010

(1.6) $$\int_0^1 e^{-px} \log \Gamma(x) \, dx = \frac{[\log(2\pi) + \gamma][p - (1 - e^{-p})]}{p^2}$$

$$+ \frac{(1 - e^{-p})}{4p} \left[ \psi\left(1 + \frac{ip}{2\pi}\right) + \psi\left(1 - \frac{ip}{2\pi}\right) \right] + 2(1 - e^{-p}) \sum_{n=1}^{\infty} \frac{\log n}{4\pi^2 n^2 + p^2}$$

Simple algebra shows that

(1.7) $$\frac{[p - (1 - e^{-p})]}{p^2} = \left[ \frac{1}{e^p - 1} - \frac{1}{p} + 1 \right] \frac{1 - e^{-p}}{p}$$

$$= \left[ \frac{1}{e^p - 1} - \frac{1}{p} + \frac{1}{2} \right] \frac{1 - e^{-p}}{p} + \frac{1}{2} \frac{1 - e^{-p}}{p}$$

$$= 2(1-e^{-p})\sum_{n=1}^{\infty}\frac{1}{4\pi^2 n^2 + p^2} + \frac{1}{2}\frac{1-e^{-p}}{p}$$

where we have used (2.5) below.

Therefore, we may write (1.6) as

(1.8) $\quad \int_0^1 e^{-px} \log \Gamma(x)\, dx = \frac{1-e^{-p}}{2p}[\gamma + \log(2\pi)] + \frac{1-e^{-p}}{4p}\left[\psi\left(1+\frac{ip}{2\pi}\right)+\psi\left(1-\frac{ip}{2\pi}\right)\right]$

$$+2(1-e^{-p})\sum_{n=1}^{\infty}\frac{\gamma+\log(2\pi n)}{4\pi^2 n^2 + p^2}$$

or equivalently as

(1.8.1) $\quad \dfrac{2p}{1-e^{-p}}\int_0^1 e^{-px}\log\Gamma(x)\,dx = \gamma + \log(2\pi) - \Lambda\left(\dfrac{p}{2\pi}\right) + 4p\sum_{n=1}^{\infty}\dfrac{\gamma+\log(2\pi n)}{4\pi^2 n^2 + p^2}$

or

(1.8.2) $\qquad = \log(2\pi) - \sum_{n=1}^{\infty}\left[\dfrac{4\pi^2 n}{4\pi^2 n^2 + p^2} - \dfrac{1}{n}\right] + 4p\sum_{n=1}^{\infty}\dfrac{\gamma+\log(2\pi n)}{4\pi^2 n^2 + p^2}$

□

Integration by parts results in (provided $p \neq 0$)

$$\int_\varepsilon^1 e^{-px}\log\Gamma(1+x)\,dx = -\frac{1}{p}\log\Gamma(1+x)e^{-px}\Big|_\varepsilon^1 + \frac{1}{p}\int_\varepsilon^1 e^{-px}\psi(1+x)\,dx$$

$$= -\frac{1}{p}\log\Gamma(1+\varepsilon)e^{-p\varepsilon} + \frac{1}{p}\int_\varepsilon^1 e^{-px}\psi(1+x)\,dx$$

$$= -\frac{1}{p}\log\Gamma(1+\varepsilon)e^{-p\varepsilon} + \frac{1}{p}\int_\varepsilon^1 e^{-px}\psi(1+x)\,dx$$

and we have

(1.9) $\quad \int_0^1 e^{-px}\psi(1+x)\,dx = p\int_0^1 e^{-px}\log\Gamma(1+x)\,dx$

$$= p\int_0^1 e^{-px}\log\Gamma(x)\,dx + p\int_0^1 e^{-px}\log x\,dx$$

We have via integration by parts

$$\int e^{-px} \log x \, dx = -\frac{1}{p} e^{-px} \log x + \frac{1}{p} \int \frac{e^{-px}}{x} dx$$

$$= \frac{Ei(-px) - e^{-px} \log x}{p}$$

where $Ei(-x)$ is the exponential integral defined by

(1.10) $$Ei(-x) = \int_x^\infty \frac{e^{-t}}{t} dt$$

We have for $x > 0$

$$Ei(-x) = \gamma + \log x + \sum_{n=1}^\infty \frac{(-1)^n x^n}{n \cdot n!}$$

and hence we see that

$$\lim_{x \to 0}[Ei(-px) - \log x] = \gamma + \log p$$

This results in the definite integral

(1.11) $$\int_0^1 e^{-px} \log x \, dx = -\frac{\gamma + \log p - Ei(-p)}{p}$$

Therefore, we have

$$\int_0^1 e^{-px} \psi(1+x) \, dx = \frac{[\log(2\pi) + \gamma][p - (1-e^{-p})]}{p} - \frac{(1-e^{-p})}{2} \Lambda\left(\frac{p}{2\pi}\right)$$

$$-[\gamma + \log p - Ei(-p)] + 2p(1-e^{-p}) \sum_{n=1}^\infty \frac{\log n}{4\pi^2 n^2 + p^2}$$

or equivalently

(1.12) $$\int_0^1 e^{-px} \psi(1+x) \, dx = \frac{1}{2}(1-e^{-p})[\gamma + \log(2\pi)] + \frac{1-e^{-p}}{4}\left[\psi\left(1+\frac{ip}{2\pi}\right) + \psi\left(1-\frac{ip}{2\pi}\right)\right]$$

$$-[\gamma + \log p - Ei(-p)] + 2p(1-e^{-p}) \sum_{n=1}^\infty \frac{\gamma + \log(2\pi n)}{4\pi^2 n^2 + p^2}$$

See also Section 7 of this paper.

□

The corresponding Laplace transform $\int_0^\infty e^{-px} \log \Gamma(x) \, dx$ for $p > 0$ was considered by

Glasser and Manna [34] who also reported that

(1.13) $$\int_0^\infty e^{-px}\psi(1+x)\,dx = \left(\frac{1}{e^p-1} - \frac{1}{p} + 1\right)\log\frac{2\pi}{p} + 2p\sum_{n=1}^\infty \frac{\log n}{4\pi^2 n^2 + p^2}$$

$$+ \frac{1}{4p}\left[\psi\left(1+\frac{ip}{2\pi}\right) + \psi\left(1-\frac{ip}{2\pi}\right)\right] - \frac{\gamma + \log p}{p}$$

which, using (1.7), may be expressed as

(1.14) $$\int_0^\infty e^{-px}\psi(1+x)\,dx = 2p\sum_{n=1}^\infty \frac{\log(2\pi n/p)}{4\pi^2 n^2 + p^2}$$

$$+ \frac{1}{4p}\left[\psi\left(1+\frac{ip}{2\pi}\right) + \psi\left(1-\frac{ip}{2\pi}\right)\right] - \frac{\gamma + \log p}{p} + \frac{1}{2}\log\frac{2\pi}{p}$$

Reference may also be made to [55].

□

**Proposition 1**

We have for $|p|<1$

(1.15) $$\frac{2\pi^2}{1-e^{-2\pi p}}\int_0^1 e^{-2\pi px}\varsigma'(0,x)\,dx = [\log(2\pi)+\gamma]\sum_{n=0}^\infty \varsigma(2n+2)p^{2n} + \frac{\pi}{2}\sum_{n=0}^\infty \varsigma(2n+3)p^{2n+1}$$

$$-\sum_{n=0}^\infty \varsigma'(2n+2)p^{2n}$$

**Proof**

Espinosa and Moll [31] showed that for $|t|<2\pi$

(1.16)
$$\int_0^1 e^{tx}\varsigma(s,x)\,dx = -2(e^t-1)\Gamma(1-s)(2\pi)^{s-2}\sum_{n=0}^\infty (-1)^n \left(\frac{t}{2\pi}\right)^n \varsigma(n+2-s)\cos\left(\frac{\pi}{2}[s-n]\right)$$

where I have corrected for a minus sign. Differentiating this with respect to $s$ gives us

$$\frac{2\pi^2}{1-e^t}\int_0^1 e^{tx}\varsigma'(s,x)\,dx$$

$$= (2\pi)^s [\Gamma(1-s)\log(2\pi) - \Gamma'(1-s)]\sum_{n=0}^\infty (-1)^n \left(\frac{t}{2\pi}\right)^n \varsigma(n+2-s)\cos\left(\frac{\pi}{2}[s-n]\right)$$

$$-(2\pi)^s \Gamma(1-s) \frac{\pi}{2} \sum_{n=0}^{\infty} (-1)^n \left(\frac{t}{2\pi}\right)^n \varsigma(n+2-s) \sin\left(\frac{\pi}{2}[s-n]\right)$$

$$-(2\pi)^s \Gamma(1-s) \sum_{n=0}^{\infty} (-1)^n \left(\frac{t}{2\pi}\right)^n \varsigma'(n+2-s) \cos\left(\frac{\pi}{2}[s-n]\right)$$

so that with $s=0$ we have

$$\frac{2\pi^2}{1-e^t} \int_0^1 e^{tx} \varsigma'(0,x)\, dx = [\log(2\pi)+\gamma] \sum_{n=0}^{\infty} (-1)^n \left(\frac{t}{2\pi}\right)^n \varsigma(n+2) \cos\left(\frac{n\pi}{2}\right)$$

$$+\frac{\pi}{2} \sum_{n=0}^{\infty} (-1)^n \left(\frac{t}{2\pi}\right)^n \varsigma(n+2) \sin\left(\frac{n\pi}{2}\right)$$

$$-\sum_{n=0}^{\infty} (-1)^n \left(\frac{t}{2\pi}\right)^n \varsigma'(n+2) \cos\left(\frac{n\pi}{2}\right)$$

Letting $t \to -2\pi p$ gives us

$$\frac{2\pi^2}{1-e^{-2\pi p}} \int_0^1 e^{-2\pi px} \varsigma'(0,x)\, dx = [\log(2\pi)+\gamma] \sum_{n=0}^{\infty} \varsigma(n+2) p^n \cos\left(\frac{n\pi}{2}\right)$$

$$+\frac{\pi}{2} \sum_{n=0}^{\infty} \varsigma(n+2) p^n \sin\left(\frac{n\pi}{2}\right) - \sum_{n=0}^{\infty} \varsigma'(n+2) p^n \cos\left(\frac{n\pi}{2}\right)$$

We have for suitably convergent series

$$\sum_{n=0}^{\infty} a_n \cos(n\pi/2) = \sum_{n=0}^{\infty} a_{2n} \cos(n\pi)$$

$$\sum_{n=0}^{\infty} a_n \sin(n\pi/2) = \sum_{n=0}^{\infty} a_{2n+1} \cos(n\pi)$$

and thus

$$\frac{2\pi^2}{1-e^{-2\pi p}} \int_0^1 e^{-2\pi px} \varsigma'(0,x)\, dx = [\log(2\pi)+\gamma] \sum_{n=0}^{\infty} (-1)^n \varsigma(2n+2) p^{2n}$$

$$+\frac{\pi}{2} \sum_{n=0}^{\infty} (-1)^n \varsigma(2n+3) p^{2n+1} - \sum_{n=0}^{\infty} (-1)^n \varsigma'(2n+2) p^{2n}$$

Applying Lerch's identity [53]

$$\varsigma'(0,x) = \log \Gamma(x) - \frac{1}{2}\log(2\pi)$$

we see that

$$\int_0^1 e^{-2\pi px} \varsigma'(0,x)\,dx = \int_0^1 e^{-2\pi px} \log \Gamma(x)\,dx + \frac{1}{2}\log(2\pi)\frac{e^{-2\pi p}-1}{2\pi p}$$

and obtain

(1.17)

$$\frac{2\pi^2}{1-e^{-2\pi p}}\int_0^1 e^{-2\pi px} \log \Gamma(x)\,dx = \frac{\pi}{2p}\log(2\pi) + [\log(2\pi) + \gamma]\sum_{n=0}^{\infty}(-1)^n \varsigma(2n+2)p^{2n}$$

$$+ \frac{\pi}{2}\sum_{n=0}^{\infty}(-1)^n \varsigma(2n+3)p^{2n+1} - \sum_{n=0}^{\infty}(-1)^n \varsigma'(2n+2)p^{2n}$$

To proceed further we require the following lemmas:

**Lemma 1**

(1.18) $$\sum_{n=1}^{\infty}\frac{\log n}{an^s + x^t} = \sum_{m=1}^{\infty}(-1)^m \frac{\varsigma'(ms)}{a^m} x^{t(m-1)}$$

**Proof**

For convenience we define for $\text{Re}(s) > 1$

$$\Omega(x) := \sum_{n=1}^{\infty}\frac{\log n}{an^s + x}$$

and obtain the derivative

$$\Omega^{(k)}(x) = (-1)^k k!\sum_{n=1}^{\infty}\frac{\log n}{(an^s + x)^{k+1}}$$

We then see that

$$\Omega^{(k)}(0) = \frac{(-1)^k k!}{a^{k+1}}\sum_{n=1}^{\infty}\frac{\log n}{n^{(k+1)s}}$$

$$= \frac{(-1)^{k+1} k!}{a^{k+1}}\varsigma'((k+1)s)$$

which results in the Taylor series

$$\Omega(p) = \sum_{k=0}^{\infty} \frac{(-1)^{k+1}}{a^{k+1}} \varsigma'((k+1)s) x^k$$

We have with $x \to x^t$

$$\Omega(x^t) = \sum_{n=1}^{\infty} \frac{\log n}{an^s + x^t}$$

$$= \sum_{k=0}^{\infty} \frac{(-1)^{k+1}}{a^{k+1}} \varsigma'((k+1)s) x^{tk}$$

and accordingly, with re-indexing, we see that

$$\sum_{n=1}^{\infty} \frac{\log n}{an^s + x^t} = \sum_{m=1}^{\infty} (-1)^m \frac{\varsigma'(ms)}{a^m} x^{t(m-1)}$$

It may be noted that trying to construct the Taylor series <u>directly</u> from $\sum_{n=1}^{\infty} \frac{\log n}{an^s + x^t}$ would be considerably more laborious.

With $a = 2\pi$ and $s = t = 2$ we obtain

(1.19) $$\sum_{n=1}^{\infty} \frac{\log n}{4\pi^2 n^2 + x^2} = \sum_{m=1}^{\infty} (-1)^m \frac{\varsigma'(2m)}{(2\pi)^{2m}} x^{2(m-1)}$$

and letting $x = 2\pi u$ results in

(1.20) $$\sum_{n=1}^{\infty} \frac{\log n}{n^2 + u^2} = \sum_{m=1}^{\infty} (-1)^m \varsigma'(2m) u^{2(m-1)}$$

In an analogous manner we may also easily show that

**Lemma 2**

(1.21) $$\sum_{n=1}^{\infty} \frac{1}{an^s + x^t} = \sum_{m=1}^{\infty} (-1)^{m+1} \frac{\varsigma(ms)}{a^m} x^{t(m-1)}$$

and with $a = 2\pi$ and $s = t = 2$ we obtain

(1.22) $$\sum_{n=1}^{\infty} \frac{1}{4\pi^2 n^2 + x^2} = \sum_{m=1}^{\infty} (-1)^{m+1} \frac{\varsigma(2m)}{(2\pi)^{2m}} x^{2(m-1)}$$

and letting $x = 2\pi u$ results in

$$(1.23) \qquad \sum_{n=1}^{\infty} \frac{1}{n^2 + u^2} = \sum_{m=1}^{\infty} (-1)^{m+1} \varsigma(2m) u^{2(m-1)}$$

□

With re-indexing we see from (1.3.1) that

$$(1.24) \qquad \sum_{n=0}^{\infty} (-1)^n \varsigma(2n+3) p^{2n+1} = \frac{\gamma - \Lambda(p)}{p}$$

Substituting (1.20), (1.23) and (1.24) in (1.17) results in

$$\frac{2\pi^2}{1-e^{-2\pi p}} \int_0^1 e^{-2\pi p x} \log \Gamma(x)\, dx = \frac{\pi}{2p} \log(2\pi) + [\log(2\pi) + \gamma] \sum_{n=1}^{\infty} \frac{1}{n^2 + p^2}$$

$$+ \frac{\pi}{2} \frac{\gamma - \Lambda(p)}{p} + \sum_{n=1}^{\infty} \frac{\log n}{n^2 + p^2}$$

and we immediately see that this corresponds with (1.8).

□

We see from (1.16) that

$$\int_0^1 e^{tx} \varsigma(0, x)\, dx = -2(e^t - 1)(2\pi)^{-2} \sum_{n=0}^{\infty} (-1)^n \left( \frac{t}{2\pi} \right)^n \varsigma(n+2) \cos\left( \frac{n\pi}{2} \right)$$

$$= -2(e^t - 1)(2\pi)^{-2} \sum_{n=0}^{\infty} \left( \frac{t}{2\pi} \right)^{2n} (-1)^n \varsigma(2n+2)$$

It is known that [4, p.264]

$$\varsigma(0, x) = \frac{1}{2} - x$$

and therefore

$$\int_0^1 e^{tx} \varsigma(0, x)\, dx = \frac{-t(e^t + 1) + 2(e^t - 1)}{2t^2}$$

Hence we obtain

$$\frac{-t(e^t + 1) + 2(e^t - 1)}{2t^2} = -2(e^t - 1)(2\pi)^{-2} \sum_{n=0}^{\infty} (-1)^n \left( \frac{t}{2\pi} \right)^{2n} \varsigma(2n+2)$$

which easily simplifies to the familiar decomposition formula

(1.25) $$\pi t \coth(\pi t) = 1 - 2\sum_{n=1}^{\infty}(-1)^n t^{2n}\varsigma(2n)$$

## 2. Corresponding trigonometric integrals

Letting $p \to ip\pi$ in (1.6) results in

$$\int_0^1 e^{-ip\pi x}\log\Gamma(x)\,dx = -\frac{[\log(2\pi)+\gamma][ip\pi-(1-e^{-ip\pi})]}{p^2\pi^2}$$

$$+\frac{(1-e^{-ip\pi})}{4ip\pi}\left[\psi\left(\frac{p}{2}\right)+\psi\left(-\frac{p}{2}\right)\right]+\frac{2(1-e^{-ip\pi})}{\pi^2}\sum_{n=1}^{\infty}\frac{\log n}{4n^2-p^2}$$

$$=\frac{[\log(2\pi)+\gamma](1-\cos p\pi)}{p^2\pi^2}+\frac{\sin p\pi}{4p\pi}\left[\psi\left(\frac{p}{2}\right)+\psi\left(-\frac{p}{2}\right)\right]+\frac{2(1-\cos p\pi)}{\pi^2}\sum_{n=1}^{\infty}\frac{\log n}{4n^2-p^2}$$

$$-i\left\langle\frac{[\log(2\pi)+\gamma][p\pi-\sin p\pi]}{p^2\pi^2}+\frac{1-\cos p\pi}{4p\pi}\left[\psi\left(\frac{p}{2}\right)+\psi\left(-\frac{p}{2}\right)\right]-\frac{2\sin p\pi}{\pi^2}\sum_{n=1}^{\infty}\frac{\log n}{4n^2-p^2}\right\rangle$$

Hence we have

(2.1) $$\int_0^1 \log\Gamma(x)\cos p\pi x\,dx$$

$$=\frac{[\log(2\pi)+\gamma](1-\cos p\pi)}{p^2\pi^2}+\frac{\sin p\pi}{4p\pi}\left[\psi\left(\frac{p}{2}\right)+\psi\left(-\frac{p}{2}\right)\right]+\frac{2(1-\cos p\pi)}{\pi^2}\sum_{n=1}^{\infty}\frac{\log n}{4n^2-p^2}$$

and

(2.2) $$\int_0^1 \log\Gamma(x)\sin p\pi x\,dx$$

$$=\frac{[\log(2\pi)+\gamma][p\pi-\sin p\pi]}{p^2\pi^2}+\frac{1-\cos p\pi}{4p\pi}\left[\psi\left(\frac{p}{2}\right)+\psi\left(-\frac{p}{2}\right)\right]-\frac{2\sin p\pi}{\pi^2}\sum_{n=1}^{\infty}\frac{\log n}{4n^2-p^2}$$

Using Hurwitz's formula for the Fourier series of the Hurwitz zeta function $\varsigma(s,x)$

$$\varsigma(s,x) = 2\Gamma(1-s)\left[\sin\left(\frac{\pi s}{2}\right)\sum_{k=1}^{\infty}\frac{\cos 2k\pi x}{(2\pi k)^{1-s}}+\cos\left(\frac{\pi s}{2}\right)\sum_{n=1}^{\infty}\frac{\sin 2k\pi x}{(2\pi k)^{1-s}}\right]$$

we showed in [20] that

$$(2.3) \int_0^1 \log \Gamma(x) \cos p\pi x \, dx$$

$$= \frac{[\gamma + \log(2\pi)]\sin p\pi}{2p\pi} + \frac{\sin p\pi}{4p\pi}\left[\psi\left(\frac{p}{2}\right) + \psi\left(-\frac{p}{2}\right)\right] + \frac{2(1-\cos p\pi)}{\pi^2}\sum_{n=1}^{\infty}\frac{\gamma + \log(2\pi n)}{4n^2 - p^2}$$

and

$$(2.4) \int_0^1 \log \Gamma(x) \sin p\pi x \, dx$$

$$= \frac{(1-\cos p\pi)[\gamma + \log(2\pi)]}{2p\pi} + \frac{1-\cos p\pi}{4p\pi}\left[\psi\left(\frac{p}{2}\right) + \psi\left(-\frac{p}{2}\right)\right] - \frac{2\sin p\pi}{\pi^2}\sum_{n=1}^{\infty}\frac{\gamma + \log(2\pi n)}{4n^2 - p^2}$$

It is easily seen that these correspond with the integrals obtained above by noting that

$$\frac{1-\cos p\pi}{p^2\pi^2} = \frac{\sin p\pi}{2p\pi} + \frac{2(1-\cos p\pi)}{\pi^2}\sum_{n=1}^{\infty}\frac{1}{4n^2 - p^2}$$

which may be deduced by substituting the identity $\coth(ix) = -i\cot x$ in the following well-known series ([12, p.296] and [39, p.378])

$$(2.5) \qquad 2\sum_{n=1}^{\infty}\frac{x}{x^2 + 4\pi^2 n^2} = \frac{1}{e^x - 1} - \frac{1}{x} + \frac{1}{2}$$

$$= \frac{1}{2}\coth\left(\frac{1}{2}x\right) - \frac{1}{x}$$

when we obtain

$$(2.5.1) \qquad \pi \cot\left(\frac{\pi x}{2}\right) = \frac{2}{x} - 4x\sum_{n=1}^{\infty}\frac{1}{4n^2 - x^2}$$

and

$$(2.5.2) \qquad \pi \cot \pi x = \frac{1}{x} - 2x\sum_{n=1}^{\infty}\frac{1}{n^2 - x^2}$$

In fact, *WolframAlpha* provides us with the finite sum

$$(2.5.3) \qquad -2x\sum_{n=1}^{k}\frac{1}{n^2 - x^2} = \psi(1+k+x) - \psi(1+k-x) - \psi(1+x) + \psi(1-x)$$

$$= \psi(1+k+x) - \psi(1+k-x) + \pi \cot \pi x - \frac{1}{x}$$

This may be easily obtained by differentiating the known relationship

$$\log \Gamma(1 \pm x + k) = \log \Gamma(1 \pm x) + \log[(1 \pm x)...(k \pm x)]$$

to get

$$\psi(1 \pm x + k) = \psi(1 \pm x) + \sum_{r=1}^{k} \frac{1}{r \pm x}$$

$\square$

Letting $p = a + ib$ in (1.6) gives us expressions for $I(a,b) = \int_0^1 \log \Gamma(x) e^{-ax} \cos bx \, dx$

and $J(a,b) = \int_0^1 \log \Gamma(x) e^{-ax} \sin bx \, dx$ and, as noted below, by parametric

differentiation we may further increase the family of integrals. For example, we have:

$$\frac{\partial}{\partial a} I(a,b) = -\int_0^1 x \log \Gamma(x) e^{-ax} \cos bx \, dx$$

$$\left. \frac{\partial}{\partial a} I(a,b) \right|_{(a,b)=(0,0)} = -\int_0^1 x \log \Gamma(x) \, dx$$

$$\frac{\partial}{\partial b} I(a,b) = -\int_0^1 x \log \Gamma(x) e^{-ax} \sin bx \, dx$$

$$\left. \frac{\partial}{\partial b} I(a,b) \right|_{(a,b)=(\infty,0)} = -\int_0^1 x \log \Gamma(x) \sin bx \, dx$$

We may differentiate (1.1) parametrically with respect to $p$ to give us evaluations of

$\int_0^1 x^k e^{-px} \log \Gamma(x) \, dx$ and, with $p = 0$, expressions for $\int_0^1 x^k \log \Gamma(x) \, dx$ may be

obtained.

For example, differentiation of (1.1) with respect to $p$ results in

(2.6) $\quad \int_0^1 x e^{-px} \log \Gamma(x) \, dx = -\frac{(1-e^{-p})(p+2) - 2p}{p^3} [\log(2\pi) + \gamma]$

$$+ \frac{pe^{-p} - (1-e^{-p})}{2p^2} \Lambda\left(\frac{p}{2\pi}\right) + \frac{(1-e^{-p})}{4p\pi} \Lambda'\left(\frac{p}{2\pi}\right)$$

$$+4(1-e^{-p})p\sum_{n=1}^{\infty}\frac{\log n}{(4\pi^2 n^2+p^2)^2}-2e^{-p}\sum_{n=1}^{\infty}\frac{\log n}{4\pi^2 n^2+p^2}$$

Letting $p=0$ gives us (with the assistance of several tedious applications of L'Hôpital's rule)

$$\lim_{p\to 0}\frac{(1-e^{-p})(p+2)-2p}{p^3}=-\frac{1}{6}$$

$$\lim_{p\to 0}\frac{pe^{-p}-(1-e^{-p})}{p^2}=-\frac{1}{2}$$

and we easily deduce the known result [31]

(2.7) $$\int_0^1 x\log\Gamma(x)\,dx=\frac{1}{6}[\log(2\pi)+\gamma]-\frac{1}{4}\gamma+\frac{1}{2\pi^2}\varsigma'(2)$$

The integrals $\int_0^1 x^k \log\Gamma(x)\,dx$ are given in [31]; they may also be computed in an efficient manner as follows.

Having regard to (1.1) we note that

$$\frac{p-(1-e^{-p})}{p^2}=\sum_{n=0}^{\infty}(-1)^n\frac{p^n}{(n+2)!}$$

and we therefore obtain

$$\frac{d^k}{dp^k}\frac{p-(1-e^{-p})}{p^2}\bigg|_{p=0}=\frac{(-1)^k}{(k+1)(k+2)}$$

Similarly, we find that

$$\frac{d^k}{dp^k}\frac{1-e^{-p}}{p}\bigg|_{p=0}=\frac{(-1)^k}{k+1}$$

The other relevant terms may be obtained using (1.20), (1.23) and (1.24).

**Proposition 2.1**

(2.4) $$\frac{2p\pi}{1-\cos p\pi}\int_0^1 \log(\sin\pi x)\sin p\pi x\,dx=-\left[2(\gamma+\log 2)+\psi\left(\frac{p}{2}\right)+\psi\left(-\frac{p}{2}\right)\right]$$

**Proof**

We multiply (2.1) by $\sin p\pi$ and multiply (2.2) by $(1-\cos p\pi)$ and add the resulting equations to obtain

$$(2.5) \quad \int_0^1 \log \Gamma(x)[\sin(p\pi(1-x)) + \sin p\pi x]\,dx$$

$$= \left(\frac{\sin^2 p\pi}{2p\pi} + \frac{(1-\cos p\pi)^2}{2p\pi}\right)\left(\gamma + \log(2\pi) + \frac{1}{2}\left[\psi\left(\frac{p}{2}\right) + \psi\left(-\frac{p}{2}\right)\right]\right)$$

$$= \frac{1-\cos p\pi}{p\pi}\left(\gamma + \log(2\pi) + \frac{1}{2}\left[\psi\left(\frac{p}{2}\right) + \psi\left(-\frac{p}{2}\right)\right]\right)$$

An obvious substitution gives us

$$\int_0^1 \log \Gamma(x)[\sin(p\pi(1-x)) + \sin p\pi x]\,dx$$

$$= \int_0^1 \log \Gamma(1-x) \sin p\pi x\,dx + \int_0^1 \log \Gamma(x) \sin p\pi x\,dx$$

and, using Euler's reflection formula for $\Gamma(x)$, we obtain

$$= \log \pi \int_0^1 \sin p\pi x\,dx - \int_0^1 \log(\sin \pi x) \sin p\pi x\,dx$$

$$= \log \pi \frac{1-\cos p\pi}{p\pi} - \int_0^1 \log(\sin \pi x) \sin p\pi x\,dx$$

Hence, we have

$$(2.6) \quad \frac{2p\pi}{1-\cos p\pi}\int_0^1 \log(\sin \pi x) \sin p\pi x\,dx = -\left[2\gamma + 2\log 2 + \psi\left(\frac{p}{2}\right) + \psi\left(-\frac{p}{2}\right)\right]$$

This integral was derived in a different manner in 2000 by Dwilewicz and Mináč [29] and was independently rediscovered by the author by another method some ten years later in [20].

We may express this as

$$(2.7) \quad \frac{4p\pi}{1-\cos 2p\pi}\int_0^1 \log(2\sin \pi x) \sin 2p\pi x\,dx = -[\psi(1+p) + \psi(1-p) + 2\gamma]$$

In passing, in connection with (2.7), we note that [20]

(2.8) $$\psi(1+p)+\psi(1-p)+2\gamma = -2p^2\sum_{n=1}^{\infty}\frac{1}{n}\frac{1}{n^2-p^2}$$

and thus

$$\int_0^1 \log(2\sin\pi x)\sin 2p\pi x\,dx = \frac{(1-\cos 2p\pi)p}{2\pi}\sum_{n=1}^{\infty}\frac{1}{n}\frac{1}{n^2-p^2}$$

Differentiation results in

$$2\pi\int_0^1 x\log(2\sin\pi x)\cos 2p\pi x\,dx$$

$$= \frac{1-\cos 2p\pi + 2p\pi\sin 2p\pi}{2\pi}\sum_{n=1}^{\infty}\frac{1}{n}\frac{1}{n^2-p^2} + \frac{(1-\cos 2p\pi)p^2}{\pi}\sum_{n=1}^{\infty}\frac{1}{n}\frac{1}{(n^2-p^2)^2}$$

from which we see that

$$\int_0^1 x\log(2\sin\pi x)\,dx = 0$$

Dwilewicz and Mináč [29] showed that

$$\sum_{n=1}^{\infty}\varsigma(2n+1)z^{2n} = \frac{2z}{1-\cos 2\pi z}\int_0^{\pi}\log(\sin x)\sin(2zx)\,dx + \log 2$$

$$= \frac{2z}{1-\cos 2\pi z}\int_0^{\pi}\log(2\sin x)\sin(2zx)\,dx$$

Integration gives us

$$\int_{\frac{1}{2}}^{u}\sum_{n=1}^{\infty}\varsigma(2n+1)z^{2n-1}\,dz = 2\int_0^{\pi}\log(2\sin x)\int_{\frac{1}{2}}^{u}\cot(zx)dz\,dx$$

$$= 2\int_0^{\pi}\log(2\sin x)[\log(\sin ux)-\log(\sin\tfrac{1}{2}x)]\,dx$$

**Proposition 2.2**

(2.9) $$\frac{4p\pi}{\sin 2p\pi}\int_0^1\log(2\sin\pi x)\cos 2p\pi x\,dx = -[\psi(1+p)+\psi(1-p)+2\gamma]$$

**Proof**

This may be obtained in a similar manner as Proposition 2.1 by multiplying (2.1) by $(1-\cos p\pi)$ and (2.2) by $\sin p\pi$ and adding the resulting equations.

**Proposition 2.3**

We have for $0 \leq p < 1$

(2.10)  $$\frac{2\pi p}{\sin p\pi} \int_0^{\frac{1}{2}} \log(\sin \pi x) \cos 2p\pi x \, dx = \sum_{n=1}^{\infty} (-1)^n \frac{n}{n^2 - p^2}$$

**Proof**

Integration by parts gives us

$$\int_0^{\frac{1}{2}} \cot(\pi x) \sin 2p\pi x \, dx = \frac{1}{\pi} \log(\sin \pi x) \sin 2p\pi x \Big|_0^{\frac{1}{2}} - 2p \int_0^{\frac{1}{2}} \log(\sin \pi x) \cos 2p\pi x \, dx$$

Writing

$$\log(\sin \pi x) \sin 2p\pi x = x \left[ \log\left(\frac{\sin \pi x}{\pi x}\right) + \log \pi x \right] \frac{\sin 2p\pi x}{x}$$

we easily see that $\lim_{x \to 0}[\log(\sin \pi x) \sin 2p\pi x] = 0$. Hence, we have

$$\int_0^{\frac{1}{2}} \cot(\pi x) \sin 2p\pi x \, dx = -2p \int_0^{\frac{1}{2}} \log(\sin \pi x) \cos 2p\pi x \, dx$$

As mentioned in Proposition 4.4 below we have

(2.11)  $$\int_a^b f(x) \cot(\alpha x / 2) \, dx = 2 \sum_{n=1}^{\infty} \int_a^b f(x) \sin \alpha nx \, dx$$

where $f(x)$ is assumed to be twice continuously differentiable on $[a, b]$. Equation (2.11) is valid provided (i) $\sin(\alpha x / 2) \neq 0 \ \forall \ x \in [a, b]$ or, alternatively, (ii) if $\sin(\alpha t / 2) = 0$, where $t \in [a, b]$, then $f(t) = 0$ also.

For example, with $f(x) = \sin 2p\pi x$ and $\alpha = 2\pi$ we have

$$\int_0^{\frac{1}{2}} \cot(\pi x) \sin 2p\pi x \, dx = 2 \sum_{n=1}^{\infty} \int_0^{\frac{1}{2}} \sin 2p\pi x \sin 2\pi nx \, dx$$

Using the definite integral

$$\int_0^{\frac{1}{2}} \sin 2\pi nx \sin 2p\pi x\, dx = -\frac{(-1)^n n}{2\pi} \frac{\sin p\pi}{n^2 - p^2}$$

we obtain (2.10).

□

We note from (2.6) and (2.9) that

$$\frac{1}{1-\cos 2p\pi}\int_0^1 \log(\sin \pi x)\sin 2p\pi x\, dx = \frac{1}{\sin 2p\pi}\int_0^1 \log(\sin \pi x)\cos 2p\pi x\, dx$$

and therefore

$$\cos p\pi \int_0^1 \log(\sin \pi x)\sin 2p\pi x\, dx = \sin p\pi \int_0^1 \log(\sin \pi x)\cos 2p\pi x\, dx$$

which may be written as

(2.12) $$\int_0^1 \log(\sin \pi x)\sin((2x-1)p\pi)\, dx = 0$$

Differentiation gives us

$$\int_0^1 (2x-1)\log(\sin \pi x)\cos((2x-1)p\pi)\, dx = 0$$

and, more generally, we have for $n \geq 0$

$$\int_0^1 (2x-1)^{2n+1} \log(\sin \pi x)\cos((2x-1)p\pi)\, dx = 0$$

With $p = 0$ we see that

(2.13) $$\int_0^1 (2x-1)^{2n+1} \log(\sin \pi x)\, dx = 0$$

## 3. Some appearances of the sine and cosine integral functions

The sine and cosine integrals, $si(x)$ and $Ci(x)$, are defined [36, p.878] by

$$\text{(3.1)} \quad si(x) = -\int_x^\infty \frac{\sin t}{t} dt$$

and for $x > 0$

$$\text{(3.2)} \quad Ci(x) = -\int_x^\infty \frac{\cos t}{t} dt$$

or equivalently

$$\text{(3.2.1)} \quad Ci(x) = \gamma + \log x + \int_0^x \frac{\cos t - 1}{t} dt$$

Note that $Si(x)$ is a slightly different sine integral which is defined in [36, p.878] and [1, p.231] by

$$\text{(3.3)} \quad Si(x) = \int_0^x \frac{\sin t}{t} dt$$

We have

$$si(x) = -\int_x^\infty \frac{\sin t}{t} dt = \int_0^x \frac{\sin t}{t} dt - \int_0^\infty \frac{\sin t}{t} dt$$

and using the well-known integral from Fourier series analysis

$$\text{(3.4)} \quad \frac{\pi}{2} = \int_0^\infty \frac{\sin t}{t} dt$$

we therefore see that the two sine integrals are intimately related by

$$\text{(3.5)} \quad si(x) = Si(x) - \frac{\pi}{2}$$

Purely on grounds of symmetry, having regard to the definitions in (3.1) and (3.2), it would have been more consistent to have denoted $Ci(x)$ by $ci(x)$, but alas the prevailing literature has not seen fit to do that.

Using Nielsen's representation of $\log \Gamma(x)$

$$\log \Gamma(x) = \frac{1}{2}\log(2\pi) - 1 - \log x + \frac{1}{\pi}\sum_{n=1}^\infty \frac{1}{n}[Ci(2n\pi)\sin 2n\pi x - si(2n\pi)\cos 2n\pi x]$$

we showed in [20] that

(3.6) $$\int_0^1 \log \Gamma(x) \cos p\pi x \, dx = \left[\frac{1}{2}\log(2\pi)-1\right]\frac{\sin p\pi}{p\pi} + \frac{Si(p\pi)}{p\pi}$$

$$+ \frac{2(1-\cos p\pi)}{\pi^2}\sum_{n=1}^{\infty}\frac{Ci(2n\pi)}{4n^2-p^2} + \frac{p\sin p\pi}{\pi^2}\sum_{n=1}^{\infty}\frac{1}{n}\frac{si(2n\pi)}{4n^2-p^2}$$

and

(3.7) $$\int_0^1 \log \Gamma(x) \sin p\pi x \, dx = \left[\frac{1}{2}\log(2\pi)-1\right]\frac{1-\cos p\pi}{p\pi} + \frac{\gamma + \log(p\pi) - Ci(p\pi)}{p\pi}$$

$$- \frac{2\sin p\pi}{\pi^2}\sum_{n=1}^{\infty}\frac{Ci(2n\pi)}{4n^2-p^2} + \frac{p(1-\cos p\pi)}{\pi^2}\sum_{n=1}^{\infty}\frac{1}{n}\frac{si(2n\pi)}{4n^2-p^2}$$

**Proposition 3.1**

(3.8)

$$\frac{2p^2}{\pi}\sum_{n=1}^{\infty}\frac{1}{n}\frac{si(2n\pi)}{4n^2-p^2} = 2 - \frac{Si(p\pi)\sin p\pi}{1-\cos p\pi} + Ci(p\pi) - \log(p\pi) + \frac{1}{2}\left[\psi\left(\frac{p}{2}\right)+\psi\left(-\frac{p}{2}\right)\right]$$

**Proof**

We multiply (3.6) by $\sin p\pi$ and (3.7) by $(1-\cos p\pi)$ to obtain

$$\int_0^1 \log \Gamma(x) \sin p\pi \cos p\pi x \, dx = \left[\frac{1}{2}\log(2\pi)-1\right]\frac{\sin^2 p\pi}{p\pi} + \frac{Si(p\pi)\sin p\pi}{p\pi}$$

$$+ \frac{2(1-\cos p\pi)\sin p\pi}{\pi^2}\sum_{n=1}^{\infty}\frac{Ci(2n\pi)}{4n^2-p^2} + \frac{p\sin^2 p\pi}{\pi^2}\sum_{n=1}^{\infty}\frac{1}{n}\frac{si(2n\pi)}{4n^2-p^2}$$

and

$$\int_0^1 \log \Gamma(x)(1-\cos p\pi)\sin p\pi x \, dx = \left[\frac{1}{2}\log(2\pi)-1\right]\frac{(1-\cos p\pi)^2}{p\pi} + \frac{\gamma + \log(p\pi) - Ci(p\pi)}{p\pi}(1-\cos p\pi)$$

$$- \frac{2(1-\cos p\pi)\sin p\pi}{\pi^2}\sum_{n=1}^{\infty}\frac{Ci(2n\pi)}{4n^2-p^2} + \frac{p(1-\cos p\pi)^2}{\pi^2}\sum_{n=1}^{\infty}\frac{1}{n}\frac{si(2n\pi)}{4n^2-p^2}$$

Therefore, upon adding the resulting equations we have

$$\int_0^1 \log \Gamma(x)\sin p\pi \cos p\pi x \, dx + \int_0^1 \log \Gamma(x)(1-\cos p\pi)\sin p\pi x \, dx$$

$$= \left[\frac{1}{2}\log(2\pi)-1\right]\frac{\sin^2 p\pi}{p\pi} + \left[\frac{1}{2}\log(2\pi)-1\right]\frac{(1-\cos p\pi)^2}{p\pi} + \frac{Si(p\pi)\sin p\pi}{p\pi}$$

$$+ \frac{\gamma+\log(p\pi)-Ci(p\pi)}{p\pi}(1-\cos p\pi) + \left[\frac{p\sin^2 p\pi}{\pi^2} + \frac{p(1-\cos p\pi)^2}{\pi^2}\right]\sum_{n=1}^{\infty}\frac{1}{n}\frac{si(2n\pi)}{4n^2-p^2}$$

We have

$$(3.9)\quad \int_0^1 \log\Gamma(x)\sin p\pi\cos p\pi x\,dx + \int_0^1 \log\Gamma(x)(1-\cos p\pi)\sin p\pi x\,dx$$

$$= \int_0^1 \log\Gamma(x)\sin p\pi x\,dx + \int_0^1 \log\Gamma(x)\sin(p\pi(1-x))\,dx$$

$$= \int_0^1 \log\Gamma(x)\sin p\pi x\,dx + \int_0^1 \log\Gamma(1-x)\sin p\pi x\,dx$$

$$= \log\pi\int_0^1 \sin p\pi x\,dx - \int_0^1 \log(\sin\pi x)\sin p\pi x\,dx$$

$$= \log\pi\,\frac{1-\cos p\pi}{p\pi} - \int_0^1 \log(\sin\pi x)\sin p\pi x\,dx$$

$$= \frac{1-\cos p\pi}{p\pi}\left[\gamma+\log(2\pi)+\frac{1}{2}[\psi(p/2)+\psi(-p/2)]\right]$$

where we have used (2.6).

Therefore, we obtain

$$\frac{1-\cos p\pi}{p\pi}\left[\gamma+\log(2\pi)+\frac{1}{2}[\psi(p/2)+\psi(-p/2)]\right]$$

$$= 2\left[\frac{1}{2}\log(2\pi)-1\right]\frac{(1-\cos p\pi)}{p\pi} + \frac{Si(p\pi)\sin p\pi}{p\pi} + \frac{2p(1-\cos p\pi)}{\pi^2}\sum_{n=1}^{\infty}\frac{1}{n}\frac{si(2n\pi)}{4n^2-p^2}$$

$$+ \frac{[\gamma+\log(p\pi)-Ci(p\pi)](1-\cos p\pi)}{p\pi}$$

A little algebra results in

(3.10)
$$\frac{2p^2}{\pi}\sum_{n=1}^{\infty}\frac{1}{n}\frac{si(2n\pi)}{4n^2-p^2} = 2 - \frac{Si(p\pi)\sin p\pi}{1-\cos p\pi} + Ci(p\pi) - \log(p\pi) + \frac{1}{2}\left[\psi\left(\frac{p}{2}\right)+\psi\left(-\frac{p}{2}\right)\right]$$

We have

$$\lim_{p \to 0} \frac{Si(p\pi)\sin p\pi}{1-\cos p\pi} = 2$$

and using (3.2.1) we see that

$$\lim_{p \to 0} Ci(p\pi) - \log(p\pi) + \frac{1}{2}\left[\psi\left(\frac{p}{2}\right) + \psi\left(-\frac{p}{2}\right)\right] = 0$$

and hence (3.10) is easily seen to be valid at $p = 0$.

For example, with $p = 1$ we have

(3.11) $\quad \dfrac{2}{\pi}\sum_{n=1}^{\infty}\dfrac{1}{n}\dfrac{si(2n\pi)}{4n^2-1} = 3 + Ci(\pi) - \gamma - \log(4\pi)$

which we obtained previously in [20] by a different method.

We recall (3.2.1)

$$Ci(x) = \gamma + \log x + \int_0^x \frac{\cos t - 1}{t} dt$$

and deduce from (3.10) that

(3.12)

$$\int_0^{p\pi} \frac{\cos t - 1}{t} dt = \frac{2p^2}{\pi}\sum_{n=1}^{\infty}\frac{1}{n}\frac{si(2n\pi)}{4n^2-p^2} - 2 - \gamma + \frac{Si(p\pi)\sin p\pi}{1-\cos p\pi} - \frac{1}{2}\left[\psi\left(\frac{p}{2}\right) + \psi\left(-\frac{p}{2}\right)\right]$$

**Proposition 3.2**

(3.13) $\quad \displaystyle\int_0^1 \log \Gamma(x)\sin \pi x\, dx = \dfrac{1}{\pi}\left[\log\dfrac{\pi}{2} + 1\right]$

**Proof**

With $p = 1$ in (3.7) we have

$$\int_0^1 \log \Gamma(x)\sin \pi x\, dx = \left[\frac{1}{2}\log(2\pi) - 1\right]\frac{2}{\pi} + \frac{\gamma + \log \pi - Ci(\pi)}{\pi} + \frac{2}{\pi^2}\sum_{n=1}^{\infty}\frac{1}{n}\frac{si(2n\pi)}{4n^2-1}$$

$$= \left[\frac{1}{2}\log(2\pi) - 1\right]\frac{2}{\pi} + \frac{\gamma + \log \pi - Ci(\pi)}{\pi}$$

$$+\frac{3+Ci(\pi)-\gamma-\log 4\pi}{\pi}$$

where we have used (3.11). Therefore, we obtain

$$\int_0^1 \log \Gamma(x) \sin \pi x \, dx = \frac{1}{\pi}\left[\log\frac{\pi}{2}+1\right]$$

The numerical value of the integral is confirmed by *WolframAlpha* and this corrects the integral in Nielsen's book [48, p.203].

**Proposition**

$$(3.14) \quad \sum_{n=1}^{\infty}\frac{Ci(2n\pi)-\gamma-\log(2\pi n)}{4n^2-p^2} = -\frac{\pi}{4p}\left[Si(p\pi)+[Ci(p\pi)-\gamma-\log(p\pi)]\frac{\sin p\pi}{1-\cos p\pi}\right]$$

**Proof**

Equating (2.3) and (3.6) gives us

$$\frac{Si(p\pi)-\sin p\pi}{p\pi}+\frac{2(1-\cos p\pi)}{\pi^2}\sum_{n=1}^{\infty}\frac{Ci(2n\pi)}{4n^2-p^2}+\frac{p\sin p\pi}{\pi^2}\sum_{n=1}^{\infty}\frac{1}{n}\frac{si(2n\pi)}{4n^2-p^2}$$

$$=\frac{\sin p\pi}{4p\pi}\left[\psi\left(\frac{p}{2}\right)+\psi\left(-\frac{p}{2}\right)+2\gamma\right]+\frac{2(1-\cos p\pi)}{\pi^2}\sum_{n=1}^{\infty}\frac{\gamma+\log(2\pi n)}{4n^2-p^2}$$

Similarly, equating (2.4) and (3.7) results in

$$\frac{\gamma+\log(p\pi)-Ci(p\pi)-(1-\cos p\pi)}{p\pi}-\frac{2\sin p\pi}{\pi^2}\sum_{n=1}^{\infty}\frac{Ci(2n\pi)}{4n^2-p^2}$$

$$+\frac{p(1-\cos p\pi)}{\pi^2}\sum_{n=1}^{\infty}\frac{1}{n}\frac{si(2n\pi)}{4n^2-p^2}$$

$$=\frac{(1-\cos p\pi)}{4p\pi}\left[\psi\left(\frac{p}{2}\right)+\psi\left(-\frac{p}{2}\right)+2\gamma\right]-\frac{2\sin p\pi}{\pi^2}\sum_{n=1}^{\infty}\frac{\gamma+\log(2\pi n)}{4n^2-p^2}$$

and hence we have two simultaneous equations involving $\sum_{n=1}^{\infty}\frac{Ci(2n\pi)}{4n^2-p^2}$ and

$\sum_{n=1}^{\infty}\frac{1}{n}\frac{si(2n\pi)}{4n^2-p^2}$. We obtain (3.8) as one of the solutions. In addition, we have

$$\frac{4(1-\cos p\pi)}{\pi^2}\sum_{n=1}^{\infty}\frac{Ci(2n\pi)}{4n^2-p^2}=\frac{\gamma+\log(p\pi)-Ci(p\pi)-(1-\cos p\pi)}{p\pi}\sin p\pi$$

$$-\frac{Si(p\pi)-\sin p\pi}{p\pi}(1-\cos p\pi)$$

$$+\frac{2(1-\cos p\pi)\sin p\pi}{\pi^2}\sum_{n=1}^{\infty}\frac{\gamma+\log(2\pi n)}{4n^2-p^2}$$

or

$$\sum_{n=1}^{\infty}\frac{Ci(2n\pi)}{4n^2-p^2}=\sum_{n=1}^{\infty}\frac{\gamma+\log(2\pi n)}{4n^2-p^2}-\frac{\pi}{4p}\left[Si(p\pi)+[Ci(p\pi)-\gamma-\log(p\pi)]\frac{\sin p\pi}{1-\cos p\pi}\right] \quad ??$$

With $p=1$ in (3.14) we obtain

(3.15) $$\sum_{n=1}^{\infty}\frac{\gamma+\log(2\pi n)}{4n^2-1}=\frac{\pi}{4}Si(\pi)+\sum_{n=1}^{\infty}\frac{Ci(2n\pi)}{4n^2-1}$$

which concurs with the more general formula noted in [20]

(3.16) $$\sum_{n=1}^{\infty}\frac{\gamma+\log(2\pi nx)}{4n^2-1}=\frac{\pi}{4}Si(\pi x)+\sum_{n=1}^{\infty}\frac{Ci(2n\pi x)}{4n^2-1}$$

Differentiating (3.16) results in

(3.17) $$\sum_{n=1}^{\infty}\frac{1}{4n^2-1}=\frac{\pi}{4}\sin(\pi x)+\sum_{n=1}^{\infty}\frac{\cos(2n\pi x)}{4n^2-1}$$

Using

(3.18) $$\sum_{n=1}^{\infty}\frac{1}{4n^2-1}=\frac{1}{2}$$

we obtain

(3.19) $$\sin x=\frac{2}{\pi}-\frac{4}{\pi}\sum_{n=1}^{\infty}\frac{\cos(2nx)}{4n^2-1}$$

which appears in [5, p.337].

Using (3.18) we may express (3.14) as

$$\frac{1}{2}[\gamma+\log(2\pi)]+\sum_{n=1}^{\infty}\frac{\log n}{4n^2-1}=\frac{\pi}{4}Si(\pi)+\sum_{n=1}^{\infty}\frac{Ci(2n\pi)}{4n^2-1}$$

Referring to (3.2.1)
$$Ci(x) = \gamma + \log x + \int_0^1 \frac{\cos xu - 1}{u} du$$

we make the summation

$$\sum_{n=1}^{\infty} \frac{Ci(2\pi nx) - \log(2\pi nx) - \gamma}{4n^2 - p^2} = \int_0^1 \sum_{n=1}^{\infty} \frac{\cos 2\pi nxu - 1}{4n^2 - p^2} \frac{du}{u}$$

assuming that we may interchange the order of summation and integration.

With $p = 1$ we have

$$\sum_{n=1}^{\infty} \frac{Ci(2\pi nx) - \log(2\pi nx) - \gamma}{4n^2 - 1} = \int_0^1 \sum_{n=1}^{\infty} \frac{\cos 2n\pi xu - 1}{4n^2 - 1} \frac{du}{u}$$

Using (3.18) and (3.19) we have

(3.20) $$\sum_{n=1}^{\infty} \frac{\cos 2n\pi xu - 1}{4n^2 - 1} = -\frac{\pi}{4} \sin \pi xu$$

and

$$\int_0^1 \sum_{n=1}^{\infty} \frac{\cos 2n\pi xu - 1}{4n^2 - 1} \frac{du}{u} = -\frac{\pi}{4} \int_0^1 \frac{\sin \pi xu}{u} du$$

Referring to (3.3) we see that

$$Si(x) = \int_0^x \frac{\sin t}{t} dt = \int_0^1 \frac{\sin xu}{u} du$$

and we therefore see that

$$\sum_{n=1}^{\infty} \frac{Ci(2\pi nx) - \log(2\pi nx) - \gamma}{4n^2 - 1} = -\frac{\pi}{4} Si(\pi x)$$

which concurs with (3.14).

**Remark:**

Letting $u = \frac{1}{2}$ in (3.19) gives us

$$\sum_{n=1}^{\infty} \frac{(-1)^n - 1}{4n^2 - 1} = -\frac{\pi}{4}$$

and we also find that

$$\sum_{n=0}^{\infty} \frac{1}{4(2n+1)^2 - 1} = \frac{\pi}{8}$$

We showed in [20] that

(3.21) $\quad Ci(ux) - \log x = \dfrac{2x \sin \pi x}{\pi} \sum\limits_{n=1}^{\infty} \dfrac{(-1)^{n+1}}{n^2 - x^2}[Ci(nu) - \log n] + \dfrac{(\gamma + \log u)\sin \pi x}{\pi x}$

which corrects a misprint in Nielsen's book [48, p.72, Eq. (5)]. This may be expressed as

$$Ci(ux) - \log x - \cos(2\pi x)\log u = \frac{2x \sin \pi x}{\pi} \sum_{n=1}^{\infty} \frac{(-1)^{n+1}}{n^2 - x^2}[Ci(nu) - \log(nu)] + \frac{\gamma}{\pi x}\sin \pi x$$

Upon differentiating (3.21) with respect to $u$, a modicum of algebra results in

(3.22) $\quad \dfrac{\pi \cos(ux)}{\sin(\pi x)} = \dfrac{1}{x} + 2x \sum\limits_{n=1}^{\infty} \dfrac{(-1)^{n+1} \cos(nu)}{n^2 - x^2}$

which appears in Spiegel's compendium [52, p.315]. This also corrects a misprint in Nielsen's book [48, p.72, Eq. (2)].

With $u = \pi$ in (3.22) we deduce the well-known (ubiquitous!) decomposition formula for the cotangent function

(3.23) $\quad \pi \cot(\pi x) = \dfrac{1}{x} - 2x \sum\limits_{n=1}^{\infty} \dfrac{1}{n^2 - x^2}$

and with $u = 0$ in (3.22) we see that

(3.24) $\quad \dfrac{\pi}{\sin \pi x} = \dfrac{1}{x} + 2x \sum\limits_{n=1}^{\infty} \dfrac{(-1)^{n+1}}{n^2 - x^2}$

which we shall also see in (8.7). We note that it is not permissible to let $u = 2\pi$ in (3.22).

Integrating (3.24) gives us for $0 \leq u < 1$

$$\int_0^u \left[\frac{\pi}{\sin \pi x} - \frac{1}{x}\right] dx = \int_0^u \sum_{n=1}^{\infty} \frac{(-1)^{n+1} 2x}{n^2 - x^2} dx$$

and thus

$$\log\tan\left(\tfrac{\pi u}{2}\right) - \log u - \log\left(\tfrac{\pi}{2}\right) - \lim_{x\to 0}\log\frac{\sin\left(\tfrac{\pi x}{2}\right)}{\tfrac{\pi x}{2}} = \sum_{n=1}^{\infty}(-1)^n \log\left(1-\frac{u^2}{n^2}\right)$$

Hence, we have for $0 \leq u < 1$

$$\log\tan\left(\tfrac{\pi u}{2}\right) - \log\left(\tfrac{\pi u}{2}\right) = \sum_{n=1}^{\infty}(-1)^n \log\left(1-\frac{u^2}{n^2}\right)$$

## 4. Some trigonometric integrals involving $x\log\Gamma(x)$

In this section we consider some trigonometric integrals involving $x\log\Gamma(x)$. The following integral is well known (see for example [31]) and will be derived in a different manner below.

**Proposition 4.1**

(4.1) $$\int_0^1 x\log\Gamma(x)\,dx = \frac{1}{4}\log(2\pi) - \frac{1}{12}[\gamma + \log(2\pi)] + \frac{\varsigma'(2)}{2\pi^2}$$

**Proof**

Letting $p \to 2p$ in (2.4) gives us

(4.2) $$\int_0^1 \log\Gamma(x)\sin 2p\pi x\,dx$$

$$= \frac{1-\cos 2p\pi}{8p\pi}[\psi(1+p) + \psi(1-p) + 2\log(2\pi) + 2\gamma] - \frac{\sin 2p\pi}{2\pi^2}\sum_{n=1}^{\infty}\frac{\gamma + \log(2\pi n)}{n^2 - p^2}$$

and differentiation with respect to $p$ results in

(4.3) $$2\pi\int_0^1 x\log\Gamma(x)\cos 2p\pi x\,dx$$

$$= \frac{1-\cos 2p\pi}{8p\pi}[\psi'(1+p) - \psi'(1-p)]$$

$$+ \frac{2p\pi\sin 2p\pi - (1-\cos 2p\pi)}{8p^2\pi}[\psi(1+p) + \psi(1-p) + 2\log(2\pi) + 2\gamma]$$

$$- \frac{\cos 2p\pi}{\pi}\sum_{m=1}^{\infty}\frac{\gamma + \log(2\pi m)}{m^2 - p^2} - \frac{p\sin 2p\pi}{\pi^2}\sum_{m=1}^{\infty}\frac{\gamma + \log(2\pi m)}{(m^2 - p^2)^2}$$

With $p = 0$ in (4.3) we obtain using L'Hôpital's rule

$$2\pi \int_0^1 x \log \Gamma(x)\,dx = \frac{\pi}{2}\log(2\pi) - \frac{1}{\pi}\sum_{m=1}^{\infty}\frac{\gamma + \log(2\pi m)}{m^2}$$

and hence we easily deduce (4.1) using

$$\sum_{m=1}^{\infty}\frac{\gamma + \log(2\pi m)}{m^2} = [\log(2\pi) + \gamma]\varsigma(2) - \varsigma'(2)$$

Using

$$\varsigma'(-1) = \frac{1}{12}(1 - \gamma - \log 2\pi) + \frac{1}{2\pi^2}\varsigma'(2)$$

we may also express (4.1) in terms of the Glaisher–Kinkelin constant $A$ whose approximate value is A $\approx$ 1.2824271291 … and $\log A = \frac{1}{12} - \varsigma'(-1)$. The value of $\varsigma'(-1) \simeq -0.16542115...$ was computed by Elizalde and Romeo [30] in 1989 (which is close to the value as computed by *WolframAlpha* $\varsigma'(-1) \simeq -0.165421143700...$).

We next give alternative derivations of some integrals which were originally determined by Mező [45] in 2016.

**Proposition 4.2**

(4.4) $$2\int_0^1 x \log \Gamma(x) \cos 2n\pi x\,dx = \frac{1}{4n} - \frac{\gamma + \log(2\pi)}{\pi^2 n^2} - \frac{\log n}{4\pi^2 n^2} - \frac{T_n}{\pi^2}$$

where $n$ is a positive integer and $T_n = \sum_{\substack{m=1 \\ m\neq n}}^{\infty}\frac{\log m}{m^2 - n^2}$.

**Proof**

With $p \to n$ in (4.3), and employing $\psi'(1-p) = \sum_{m=1}^{\infty}\frac{1}{(m-p)^2}$, we have

$$\lim_{p\to n}(1 - \cos 2p\pi)\psi'(1-p) = \lim_{p\to n}\frac{1 - \cos 2p\pi}{(n-p)^2}$$

and two applications of L'Hôpital's rule gives us

$$\lim_{p\to n}(1 - \cos 2p\pi)\psi'(1-p) = -2\pi^2$$

Hence, we obtain

$$\lim_{p\to n}\frac{1 - \cos 2p\pi}{8p\pi}[\psi'(1+p) - \psi'(1-p)] = -\frac{\pi}{4n}$$

Similarly, employing $\psi(1-p) = -\gamma + \sum_{m=1}^{\infty}\left(\frac{1}{m} - \frac{1}{m-p}\right)$, we have

$$\lim_{p \to n}[2p\pi \sin 2p\pi - (1 - \cos 2p\pi)]\psi(1-p) = -\lim_{p \to n}\frac{[2p\pi \sin 2p\pi - (1 - \cos 2p\pi)]}{n-p}$$

and L'Hôpital's rule results in

$$\lim_{p \to n}[2p\pi \sin 2p\pi - (1 - \cos 2p\pi)]\psi(1-p) = 4n\pi^2$$

Hence, we obtain

$$\lim_{p \to n}\frac{2p\pi \sin 2p\pi - (1 - \cos 2p\pi)}{8p^2\pi}\left[\psi(1+p) + \psi(1-p) + 2\log(2\pi) + 2\gamma\right] = \frac{\pi}{2n}$$

We next consider the term

$$-\frac{\cos 2p\pi}{\pi}\sum_{m=1}^{\infty}\frac{\gamma + \log(2\pi m)}{m^2 - p^2} - \frac{p \sin 2p\pi}{\pi^2}\sum_{m=1}^{\infty}\frac{\gamma + \log(2\pi m)}{(m^2 - p^2)^2}$$

and in the limit as $p \to n$ this becomes

$$-\frac{1}{\pi}\sum_{\substack{m=1 \\ m \neq n}}^{\infty}\frac{\gamma + \log(2\pi m)}{m^2 - n^2} - \lim_{p \to n}[\gamma + \log(2\pi n)]\left[\frac{\cos 2p\pi}{\pi}\frac{1}{n^2 - p^2} + \frac{p \sin 2p\pi}{\pi^2}\frac{1}{(n^2 - p^2)^2}\right]$$

We write

$$\frac{\cos 2p\pi}{\pi}\frac{1}{n^2 - p^2} + \frac{p \sin 2p\pi}{\pi^2}\frac{1}{(n^2 - p^2)^2} = \frac{1}{\pi^2}\frac{\pi(n^2 - p^2)\cos 2p\pi + p \sin 2p\pi}{(n^2 - p^2)^2}$$

and L'Hôpital's rule gives us

$$\frac{1}{\pi^2}\lim_{p \to n}\left[\frac{\pi(n^2 - p^2)\cos 2p\pi + p \sin 2p\pi}{(n^2 - p^2)^2}\right] = \frac{1}{4\pi n^2}$$

Assembling the various components, and using the following lemma (see below)

(4.5) $$\sum_{\substack{m=1 \\ m \neq n}}^{\infty}\frac{1}{m^2 - n^2} = \frac{3}{4}\frac{1}{n^2}$$

we see that [45]

$$2\int_0^1 x\log\Gamma(x)\cos 2n\pi x\,dx = \frac{1}{4n} - \frac{\gamma+\log(2\pi)}{\pi^2 n^2} - \frac{\log n}{4\pi^2 n^2} - \frac{1}{\pi^2}\sum_{\substack{m=1\\m\neq n}}^{\infty}\frac{\log m}{m^2-n^2}$$

In particular we have [45]

(4.6) $$2\int_0^1 x\log\Gamma(x)\cos 2\pi x\,dx = \frac{1}{4} - \frac{\gamma+\log(2\pi)}{\pi^2} - \frac{1}{\pi^2}\sum_{m=2}^{\infty}\frac{\log m}{m^2-1}$$

Similarly, with $p = n+\tfrac{1}{2}$, we may evaluate $\int_0^1 x\log\Gamma(x)\cos(2n+1)\pi x\,dx$.

**Lemma**

$$\sum_{\substack{m=1\\m\neq n}}^{\infty}\frac{1}{m^2-n^2} = \frac{3}{4}\frac{1}{n^2}$$

**Proof**

This may be easily deduced by considering (2.5.2) in the form

$$\pi\cot\pi x + \frac{2x}{n^2-x^2} = \frac{1}{x} - 2x\sum_{\substack{m=1\\m\neq n}}^{\infty}\frac{1}{m^2-x^2}$$

and take the limit

$$\lim_{x\to n}\left[\pi\cot\pi x + \frac{2x}{n^2-x^2}\right] = \frac{1}{n} - 2\lim_{x\to n} x\sum_{\substack{m=1\\m\neq n}}^{\infty}\frac{1}{m^2-x^2}$$

With two applications of L'Hôpital's rule we have

$$\lim_{x\to n}\frac{\pi(n^2-x^2)\cos\pi x + 2x\sin\pi x}{(n^2-x^2)\sin\pi x} = -\frac{1}{2n}$$

and the lemma follows directly. A different proof is outlined in [12, p.67].

□

It is known that [17]

(4.7) $$\int_a^b p(x)dx = 2\sum_{n=0}^{\infty}\int_a^b p(x)\cos\alpha nx\,dx$$

Equation (4.7) is valid provided (i) $\sin(\alpha x/2) \neq 0 \ \forall\ x\in[a,b]$ or, alternatively, (ii) if $\sin(\alpha\eta/2) = 0$ for some $\eta\in[a,b]$ then $p(\eta)=0$ also.

Since $x\log\Gamma(x) = x\log\Gamma(1+x) - x\log x$ we see that $\lim_{x\to 0} x\log\Gamma(x) = 0$. Hence $x\log\Gamma(x)$ satisfies the conditions for (4.7) to apply. Therefore, we have

$$-\frac{1}{2}\int_0^1 x\log\Gamma(x)\,dx = 2\sum_{n=1}^{\infty}\int_0^1 x\log\Gamma(x)\cos 2n\pi x\,dx$$

Hence using (4.1) we obtain

$$\frac{1}{4}\log(2\pi) - \frac{1}{12}[\gamma + \log(2\pi)] + \frac{\varsigma'(2)}{2\pi^2} = -4\sum_{n=1}^{\infty}\left[\frac{1}{4n} - \frac{\gamma + \log(2\pi)}{\pi^2 n^2} - \frac{\log n}{4\pi^2 n^2} - \frac{T_n}{\pi^2}\right]$$

$$= \frac{2}{3}[\gamma + \log(2\pi)] - \frac{\varsigma'(2)}{\pi^2} - 4\sum_{n=1}^{\infty}\left[\frac{1}{4n} - \frac{T_n}{\pi^2}\right]$$

and thus

$$\frac{1}{4}\log(2\pi) - \frac{3}{4}[\gamma + \log(2\pi)] + \frac{3\varsigma'(2)}{2\pi^2} = -4\sum_{n=1}^{\infty}\left[\frac{1}{4n} - \frac{T_n}{\pi^2}\right]$$

(albeit the appearance of the term $n^{-1}$ on the right-hand side suggests issues of convergence).

□

Referring to (4.4) we have the summation

$$2\sum_{n=1}^{\infty}\frac{1}{n}\int_0^1 x\log\Gamma(x)\cos 2n\pi x\,dx = \frac{1}{4}\varsigma(2) - \frac{\gamma + \log(2\pi)}{\pi^2}\varsigma(3) + \frac{1}{4\pi^2}\varsigma'(3) - \frac{1}{\pi^2}\sum_{n=1}^{\infty}\frac{T_n}{n}$$

Assuming that

$$\sum_{n=1}^{\infty}\frac{1}{n}\int_0^1 x\log\Gamma(x)\cos 2n\pi x\,dx = \int_0^1 x\log\Gamma(x)\sum_{n=1}^{\infty}\frac{\cos 2n\pi x}{n}\,dx$$

$$= -\int_0^1 x\log\Gamma(x)\log[2\sin \pi x]\,dx$$

we obtain

$$\int_0^1 x\log\Gamma(x)\log[2\sin \pi x]\,dx = -\frac{1}{2}\left[\frac{1}{4}\varsigma(2) - \frac{\gamma + \log(2\pi)}{\pi^2}\varsigma(3) + \frac{1}{4\pi^2}\varsigma'(3) - \frac{1}{\pi^2}\sum_{n=1}^{\infty}\frac{T_n}{n}\right]$$

□

We showed in [17] that

$$\int_0^1 \log\Gamma(1+x)\sin 2n\pi x\,dx = \frac{Ci(2n\pi)}{2n\pi}$$

The following integral was originally determined by Mező [45] by a different method in 2016.

**Proposition 4.3**

$$(4.8) \quad \int_0^1 x \log \Gamma(x) \sin 2n\pi x \, dx = \frac{1}{4\pi} \left[ \frac{\gamma + \log n - H_n}{n} + \frac{1}{n^2} \right]$$

where $n$ is a positive integer and $H_n = \sum_{k=1}^{n} \frac{1}{k}$.

**Proof**

Letting $p \to 2p$ in (2.1) gives us

$$(4.9) \quad \int_0^1 \log \Gamma(x) \cos 2p\pi x \, dx$$

$$= \frac{[\log(2\pi) + \gamma](1 - \cos 2p\pi)}{4p^2\pi^2} + \frac{\sin 2p\pi}{8p\pi} [\psi(1+p) + \psi(1-p)] + \frac{(1 - \cos 2p\pi)}{2\pi^2} \sum_{n=1}^{\infty} \frac{\log n}{n^2 - p^2}$$

and differentiation with respect to $p$ results in

$$(4.10) \quad -2\pi \int_0^1 x \log \Gamma(x) \sin 2p\pi x \, dx$$

$$= \frac{[\gamma + \log(2\pi)][2p\pi \sin 2p\pi - 2(1 - \cos 2p\pi)]}{4p\pi}$$

$$+ \frac{2p\pi \cos 2p\pi - \sin 2p\pi}{8p^2\pi} [\psi(1+p) + \psi(1-p)]$$

$$+ \frac{\sin 2p\pi}{8p\pi} [\psi'(1+p) - \psi'(1-p)]$$

$$+ \frac{\sin 2p\pi}{\pi} \sum_{n=1}^{\infty} \frac{\log n}{n^2 - p^2} + \frac{p(1 - \cos 2p\pi)}{\pi^2} \sum_{n=1}^{\infty} \frac{\log n}{(n^2 - p^2)^2}$$

We see that

$$\lim_{p \to n} \sin 2p\pi \, \psi(1-p) = -\lim_{p \to n} \frac{\sin 2p\pi}{(n-p)} = \frac{2\pi}{n}$$

and

$$\lim_{p \to n} \left[ \frac{2p\pi \cos 2p\pi}{8p^2\pi} \psi(1-p) - \frac{\sin 2p\pi}{8p\pi} \psi'(1-p) \right]$$

$$= \lim_{p \to n} \left[ -\frac{2p\pi \cos 2p\pi}{8p^2\pi(n-p)} - \frac{\sin 2p\pi}{8p\pi(n-p)^2} \right]$$

L'Hôpital's rule easily shows that this limit vanishes.

We also have

$$\lim_{p \to n} \frac{\sin 2p\pi}{\pi} \sum_{m=1}^{\infty} \frac{\log n}{m^2 - p^2} = \frac{\log n}{\pi} \lim_{p \to n} \frac{\sin 2p\pi}{n^2 - p^2}$$

$$= -\frac{\log n}{\pi} \lim_{p \to n} \frac{\pi \cos 2p\pi}{p}$$

$$= -\frac{\log n}{n}$$

Finally, we have

$$\lim_{p \to n} \frac{p(1 - \cos 2p\pi)}{\pi^2} \sum_{m=1}^{\infty} \frac{\log m}{(m^2 - p^2)^2} = \frac{n \log n}{\pi^2} \lim_{p \to n} \frac{1 - \cos 2p\pi}{(n^2 - p^2)^2}$$

$$= \frac{\log n}{2n}$$

Assembling the various components together, using $\psi(1+n) = H_n - \gamma$, we then obtain (4.8).

□

Assuming that the following interchange is valid

$$\sum_{n=1}^{\infty} \frac{1}{n} \int_0^1 x \log \Gamma(x) \sin 2n\pi x \, dx = \int_0^1 x \log \Gamma(x) \sum_{n=1}^{\infty} \frac{\sin 2n\pi x}{n} \, dx$$

and using (4.8) we obtain

$$\int_0^1 x \log \Gamma(x) \sum_{n=1}^{\infty} \frac{\sin 2n\pi x}{n} \, dx = \frac{1}{4\pi} \sum_{n=1}^{\infty} \left[ \frac{\gamma + \log n - H_n}{n^2} + \frac{1}{n^3} \right]$$

Employing the Fourier series

$$\sum_{n=1}^{\infty} \frac{\sin 2n\pi x}{n} = \frac{\pi}{2}(1 - 2x)$$

this results in

$$\int_0^1 x \log \Gamma(x)\,dx - 2\int_0^1 x^2 \log \Gamma(x)\,dx = \frac{1}{2\pi^2} \sum_{n=1}^{\infty}\left[\frac{\gamma + \log n - H_n}{n^2} + \frac{1}{n^3}\right]$$

Using Euler's result

$$\sum_{n=1}^{\infty} \frac{H_n}{n^2} = 2\varsigma(3)$$

we obtain

$$\sum_{n=1}^{\infty}\left[\frac{\gamma + \log n - H_n}{n^2} + \frac{1}{n^3}\right] = \gamma\varsigma(2) - \varsigma'(2) - \varsigma(3)$$

Then, using (4.1), we deduce the well-known result [31]

(4.11) $$\int_0^1 x^2 \log \Gamma(x)\,dx = \frac{1}{12}\log(2\pi) - \frac{1}{12}\gamma + \frac{1}{4\pi^2}\varsigma(3) + \frac{1}{2\pi^2}\varsigma'(2)$$

$\square$

It is well known that the Bernoulli polynomials $B_n(t)$ may be represented by the Fourier series (see for example Milne-Thompson [46, p. 239], Apostol [5] and [20])

(4.12.1) $$B_{2N+1}(t) = (-1)^{N+1} 2(2N+1)! \sum_{n=1}^{\infty} \frac{\sin 2n\pi t}{(2\pi n)^{2N+1}} \quad , N = 0,1,2,...$$

(4.12.2) $$B_{2N}(t) = (-1)^{N+1} 2(2N)! \sum_{n=1}^{\infty} \frac{\cos 2n\pi t}{(2\pi n)^{2N}} \quad , N = 1,2,...$$

Since [4, p.264]

(4.12.3) $$\varsigma(-m,t) = -\frac{B_{m+1}(t)}{m+1}$$

we may deduce the above identities directly from Hurwitz's formula for the Fourier expansion of the Riemann zeta function $\varsigma(s,x)$ as reported in Titchmarsh's treatise [58, p.37].

Adamchik (see [59] and [60]) showed that

(4.12.4) $$\varsigma'(1-2N,1-t) - \varsigma'(1-2N,t) = (-1)^{N+1} 2\pi(2N-1)! \sum_{n=1}^{\infty} \frac{\sin 2n\pi t}{(2\pi n)^{2N}}$$

(4.12.5) $$\varsigma'(-2N,1-t)+\varsigma'(-2N,t)=(-1)^N 2\pi(2N)!\sum_{n=1}^{\infty}\frac{\cos 2n\pi t}{(2\pi n)^{2N+1}}$$

Hence, for example, we obtain

$$\int_0^1 x\log\Gamma(x)\sum_{n=1}^{\infty}\frac{\sin 2n\pi x}{(2\pi n)^{2N+1}}dx=\sum_{n=1}^{\infty}\frac{1}{(2\pi n)^{2N+1}}\frac{1}{4\pi}\left[\frac{\gamma+\log n-H_n}{n}+\frac{1}{n^2}\right]$$

and thus we have

(4.12.6) $$\int_0^1 x B_{2N+1}(x)\log\Gamma(x)\,dx=\frac{(-1)^{N+1}(2N+1)!}{2\pi}\sum_{n=1}^{\infty}\frac{1}{(2\pi n)^{2N+1}}\left[\frac{\gamma+\log n-H_n}{n}+\frac{1}{n^2}\right]$$

This may be compared with the known integral

(4.12.7) $$\int_0^1 B_{2n+1}(x)\cot\pi x\,dx=\frac{(-1)^{n+1}2(2n+1)!\varsigma(2n+1)}{(2\pi)^{2n+1}}$$

which appears in Abramowitz and Stegun [1, p.807]. There are many derivations; for example, see the recent one by Dwilewicz and Mináč [29]. A similar identity was also derived by Espinosa and Moll [31] in the form

(4.12.8) $$\int_0^1 B_{2n}(x)\log\sin\pi x\,dx=(-1)^n\frac{(2n)!\varsigma(2n+1)}{(2\pi)^{2n}}$$

and it is easily shown that equation (4.12.8) above is equivalent to (4.12.7) following a simple integration by parts.

Glasser's formula [33] gives us

(4.12.9) $$\int_0^1 B_{2n+1}(x)\psi(x)\,dx=-\frac{\pi}{2}\int_0^1 B_{2n+1}(x)\cot\pi x\,dx$$

and hence we have

(4.12.10) $$\int_0^1 B_{2n}(x)\log\Gamma(x)\,dx=\frac{\pi}{2(2n+1)}\int_0^1 B_{2n+1}(x)\cot\pi x\,dx$$

**The Fourier Series for** $\log G(x)$

It is known that

(4.13) $$\log G(1+x)-x\log\Gamma(x)=\varsigma'(-1)-\varsigma'(-1,x)$$

where $G(x)$ is the Barnes double gamma function [53, p. 24]. This functional equation was derived by Vardi in 1988 and also by Gosper in 1997 (an elementary derivation is contained in [24]).

We have the well-known Hurwitz's formula for the Fourier expansion of the Riemann zeta function $\varsigma(s,x)$ as reported in Titchmarsh's treatise [58, p.37]

$$(4.14) \qquad \varsigma(s,x) = 2\Gamma(1-s)\left[\sin\left(\frac{\pi s}{2}\right)\sum_{n=1}^{\infty}\frac{\cos 2n\pi x}{(2\pi n)^{1-s}} + \cos\left(\frac{\pi s}{2}\right)\sum_{n=1}^{\infty}\frac{\sin 2n\pi x}{(2\pi n)^{1-s}}\right]$$

where $\mathrm{Re}(s) < 0$ and $0 < x \leq 1$. In 2000, Boudjelkha [8] showed that this formula also applies in the region $\mathrm{Re}(s) < 1$. It may be noted that when $x = 1$ this reduces to Riemann's functional equation for $\varsigma(s)$.

By differentiating (4.14) we obtain with $s = -1$

(4.15)
$$\varsigma'(-1,t) = \frac{1}{4\pi}\sum_{n=1}^{\infty}\frac{\sin 2n\pi t}{n^2} - \frac{1}{2\pi^2}[\log(2\pi)+\gamma-1]\sum_{n=1}^{\infty}\frac{\cos 2n\pi t}{n^2} - \frac{1}{2\pi^2}\sum_{n=1}^{\infty}\frac{\log n}{n^2}\cos 2n\pi t$$

The Fourier coefficients of $x\log\Gamma(x)$ are known from (4.1), (4.4) and (4.8) and we may therefore obtain the Fourier series for $\log G(x)$ for $0 < x < 1$ (as originally determined by Mező [45] by a different method):

$$(4.16) \qquad \log G(x) = a_0 + \sum_{n=1}^{\infty} a_n \cos(2n\pi x) + b_n \sin(2n\pi x)$$

where

$$a_0 = \frac{1}{12} - 2\log A - \frac{1}{4}\log(2\pi)$$

$$a_n = \frac{1}{2\pi^2 n^2}\left[\frac{1}{2}\log n - \gamma - \log(2\pi) - 1\right] - \frac{1}{4n} - \frac{T_n}{n^2}$$

$$b_n = \frac{1}{2\pi n}\left[\frac{1}{2n} - \gamma - \log(4\pi^2 n) - H_n\right]$$

With $t = 1$ in (4.15) we obtain the familiar result

$$\varsigma'(-1) = -\frac{1}{12\pi^2}[\log(2\pi)+\gamma-1] + \frac{\varsigma'(2)}{2\pi^2}$$

**Another derivation of Alexeiewsky's theorem**

Differentiating (4.13) results in

$$\frac{G'(1+x)}{G(1+x)} - x\psi(x) - \log\Gamma(x) = -\frac{d}{dx}\varsigma'(-1,x)$$

$$= \varsigma(0,x) - \varsigma'(0,x)$$

$$= \frac{1}{2} - x - \varsigma'(0,x)$$

since $\dfrac{\partial}{\partial x}\dfrac{\partial}{\partial s}\varsigma(s,x) = -\varsigma(s+1,x) - s\varsigma'(s+1,x)$ and $\varsigma(1-m,x) = -\dfrac{B_m(x)}{m}$ [5, p.264].

Hence, using Lerch's theorem, we obtain

$$\frac{G'(1+x)}{G(1+x)} - x\psi(x) = \frac{1}{2} - x + \frac{1}{2}\log(2\pi)$$

and integrating this we have

$$G(1+u) - \frac{u}{2}\log(2\pi) - \frac{u}{2}(1-u) = \int_0^u x\psi(x)\,dx$$

Integration by parts then gives us Alexeiewsky's theorem [53, p.32]

(4.17.1) $$\int_0^u \log\Gamma(x)\,dx = \frac{u}{2}\log(2\pi) + \frac{u}{2}(1-u) + u\log\Gamma(u) - \log G(1+u)$$

or

(4.17.2) $$\int_0^u \log\Gamma(1+x)\,dx = \frac{1}{2}[\log(2\pi)-1]u - \frac{1}{2}u^2 + u\log\Gamma(1+u) - \log G(1+u)$$

One may derive Lerch's theorem in a similar manner (see [24]).

**Proposition 4.4**

We showed in [17] that

(4.18) $$\int_a^b p(x)\cot(\alpha x/2)\,dx = 2\sum_{n=1}^{\infty}\int_a^b p(x)\sin\alpha nx\,dx$$

where $p(x)$ is assumed to be twice continuously differentiable on $[a,b]$. Equation (4.18) is valid provided (i) $\sin(\alpha x/2) \neq 0 \ \forall \ x \in [a,b]$ or, alternatively, (ii) if $\sin(\alpha\eta/2) = 0$ then $p(\eta) = 0$ also where $\eta \in [a,b]$.

The particular function $p(x) = x \log \Gamma(1+x)$ is twice continuously differentiable on the interval $[0,1]$ and we see that $p(0) = p(1) = 0$. Hence, $p(x) = x \log \Gamma(1+x)$ satisfies the conditions of Proposition 4.4 on the interval $[0,1]$ with $\alpha = 2\pi$. We therefore have

$$(4.19) \qquad \int_0^1 x \log \Gamma(1+x) \cot(\pi x)\, dx = 2 \sum_{n=1}^{\infty} \int_0^1 x \log \Gamma(1+x) \sin(2n\pi x)\, dx$$

As mentioned below in Proposition 4.5, we proved in [19] that the above requirement for $p(x)$ to be a twice continuously differentiable function could be considerably relaxed in certain prescribed circumstances.

**Proposition 4.5**

We assume that $f(x)$ is continuous on $[0,b]$ and that $f(0) = 0$. Let us also assume that $f(x)$ is non-negative and monotonic increasing on $[0,b]$ where $b \leq \pi/\alpha$.

Then we have

$$(4.20) \qquad \int_0^b f(x) \cot(\alpha x/2)\, dx = 2 \sum_{n=1}^{\infty} \int_0^b f(x) \sin n\alpha x\, dx$$

This was proved in [19]. Note that, unlike Proposition 4.4, it is not a requirement for $f(x)$ to be a differentiable function.

We now wish to extend the scope of the above proposition.

**Proposition 4.6 (i)**

Let $q(x)$ be continuous and non-negative on $[0,1]$ with the boundary conditions $q(0) = q(1) = 0$. Let us also assume that $q(x)$ is monotonic increasing in the interval $[0,c]$ and is monotonic decreasing on $[c,1]$ where $c \in (0,1)$. Then we have

$$(4.21) \qquad \int_0^1 q(x) \cot(\pi x)\, dx = 2 \sum_{n=1}^{\infty} \int_0^1 q(x) \sin 2n\pi x\, dx$$

We again note that, unlike Proposition 4.4, it is not a requirement for $q(x)$ to be a differentiable function.

**Proof**

We have

$$\int_0^1 q(x) \cot(\pi x)\, dx = \int_0^c q(x) \cot(\pi x)\, dx + \int_c^1 q(x) \cot(\pi x)\, dx$$

and a change of variables gives us

$$\int_c^1 q(x)\cot(\pi x)\,dx = -\int_{1-c}^0 q(1-u)\cot(\pi(1-u))\,du$$

$$= -\int_0^{1-c} q(1-u)\cot(\pi u)\,du$$

Thus, we have

$$\int_0^1 q(x)\cot(\pi x)\,dx = \int_0^c q(x)\cot(\pi x)\,dx - \int_0^{1-c} q(1-x)\cot(\pi x)\,dx$$

We note that $q(x)$ satisfies the conditions of Proposition 4.5 on the interval $[0,c]$ and we therefore have

$$\int_0^c q(x)\cot(\pi x)\,dx = 2\sum_{n=1}^\infty \int_0^c q(x)\sin(2n\pi x)\,dx$$

Similarly, we see that $q(1-x)$ also satisfies the conditions of Proposition 4.5 on the interval $[0,1-c]$ and hence we have

$$\int_0^{c-1} q(1-x)\cot(\pi x)\,dx = 2\sum_{n=1}^\infty \int_0^{c-1} q(1-x)\sin(2n\pi x)\,dx$$

$$= -2\sum_{n=1}^\infty \int_1^c q(u)\sin(2n\pi(1-u))\,du$$

$$= 2\sum_{n=1}^\infty \int_1^c q(u)\sin(2n\pi u)\,du$$

$$= -2\sum_{n=1}^\infty \int_c^1 q(u)\sin(2n\pi u)\,du$$

Therefore, we obtain

$$\int_0^c q(x)\cot(\pi x)\,dx - \int_0^{1-c} q(1-x)\cot(\pi x)\,dx$$

$$= 2\sum_{n=1}^\infty \int_0^c q(x)\sin(2n\pi x)\,dx + 2\sum_{n=1}^\infty \int_c^1 q(x)\sin(2n\pi x)\,dx$$

and the proof of (4.21) easily follows.

This proposition may be employed for a function such as $\log x \log(1-x)$ whose derivative is infinite at $x=0$ and $x=1$.

**Proposition 4.6 (ii)**

It may be noted that Proposition 4.6(i) also applies in the case where $q(x)$ is monotonic decreasing on $[0,c]$ and is monotonic increasing on $[c,1]$. To see this, we consider the function $Q(x)$ where $Q(x) = -q(x)$. The function $x \log x$ meets these conditions

**Proposition 4.6 (iii)**

Suppose (i) $q(x)$ is continuous on $[c,1]$, (ii) $q(1) = 0$ and (iii) $q(x)$ is monotonic decreasing on $[c,1]$ where $0 < c \leq 1$. Then we have

$$(4.22) \qquad \int_c^1 q(x) \cot(\pi x) \, dx = 2 \sum_{n=1}^{\infty} \int_c^1 q(x) \sin(2n\pi x) \, dx$$

**Proof**

We note from the above that

$$\int_c^1 q(x) \cot(\pi x) \, dx = -\int_0^{1-c} q(1-u) \cot(\pi u) \, du$$

$$= 2 \sum_{n=1}^{\infty} \int_c^1 q(u) \sin(2n\pi u) \, du$$

□

We note that, $q(x) = x \log x$ satisfies the conditions of Proposition 4.6(ii) on the interval $[0,1]$ with $c = e^{-1}$. We therefore have

$$(4.23) \qquad \int_0^1 x \log x \cot(\pi x) \, dx = 2 \sum_{n=1}^{\infty} \int_0^1 x \log x \sin(2n\pi x) \, dx$$

Combining (4.23) with (4.19) we deduce that

$$(4.24) \qquad \int_0^1 x \log \Gamma(x) \cot(\pi x) \, dx = 2 \sum_{n=1}^{\infty} \int_0^1 x \log \Gamma(x) \sin(2n\pi x) \, dx$$

**Remark:**

Letting $q(x) = x \log \Gamma(x)$ we see that (i) $q(0) = q(1) = 0$ and (ii) $q(x) = x \log \Gamma(x)$ is positive and has one turning point at $x_0 \in (0,1)$. The graphical representation of $x \log \Gamma(x)$ is shown below and this indicates that it is monotonic increasing in the interval $[0, x_0)$ and is monotonic decreasing in the interval $(x_0, 1]$.

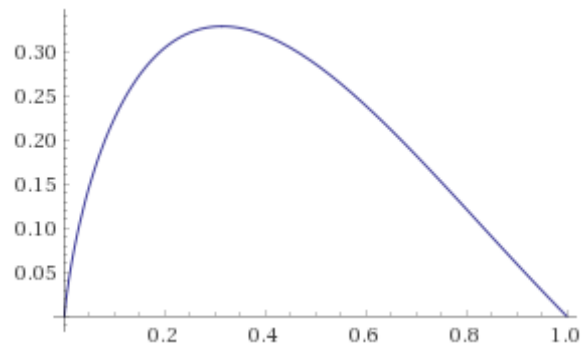

Graph of $q(x) = x \log \Gamma(x)$

Therefore, *subject* to obtaining an analytical proof of the requisite monotonicity conditions, one could have directly used Proposition 4.6(ii) with $q(x) = x \log \Gamma(x)$ to deduce (4.24).

Integration by parts gives us

$$\int x \log x \sin 2n\pi x \, dx = -\frac{2n\pi x \log x \cos 2n\pi x - (1 + \log x) \sin 2n\pi x + Si(2\pi nx)}{4\pi^2 n^2}$$

and the definite integral results

(4.25) $$\int_0^1 x \log x \sin 2n\pi x \, dx = -\frac{Si(2n\pi)}{4\pi^2 n^2}$$

Hence, using Proposition 4.6(ii) we obtain

(4.26) $$\int_0^1 x \log x \cot(\pi x) \, dx = -\frac{1}{2\pi^2} \sum_{n=1}^{\infty} \frac{Si(2n\pi)}{n^2}$$

However, this result is <u>not</u> numerically confirmed by *WolframAlpha* which gives us $\int_0^1 x \log x \cot(\pi x) \, dx \simeq -0.121552...$ for the integral and a slightly different value for the series $\frac{1}{2\pi^2} \sum_{n=1}^{\infty} \frac{Si(2n\pi)}{n^2} \simeq 0.121155...$. It appears that obtaining a closed form for $\sum_{n=1}^{\infty} \frac{Si(2n\pi)}{n^2} \simeq 2.3915...$ is as elusive as $\sum_{n=1}^{\infty} \frac{1}{n^3}$.

We conjectured in [17, Eq. (6.92b)] that

(4.27) $$\sum_{n=1}^{\infty} \frac{\varsigma(2n)}{(2n+1)^2} + \frac{1}{4\pi} \sum_{n=1}^{\infty} \frac{Si(2n\pi)}{n^2} = \frac{1}{2} \quad (?)$$

but *WolframAlpha* gives us a slightly different value for

$$\sum_{n=1}^{\infty} \frac{\varsigma(2n)}{(2n+1)^2} + \frac{1}{4\pi} \sum_{n=1}^{\infty} \frac{Si(2n\pi)}{n^2} \simeq 0.498133...$$

□

An alternative derivation of (4.26) is developed below. Integration by parts results in

$$\int_0^1 x \log x \cot(\pi x) \, dx = \frac{1}{\pi} x \log x \log \sin(\pi x) \Big|_0^1 - \frac{1}{\pi} \int_0^1 (\log x + 1) \log \sin(\pi x) \, dx$$

We have

$$\log x \log \sin(\pi x) = \log x \log[2 \sin(\pi x / 2)] + \log x \log \cos(\pi x / 2)$$

and using

$$\log x \log \cos(\pi x / 2) = \log x \log \frac{\cos(\pi x / 2)}{1-x} + \log x \log(1-x)$$

we see that

$$\lim_{x \to 1} \log \frac{\cos(\pi x / 2)}{1-x} = \log \left[ \lim_{x \to 1} \frac{\cos(\pi x / 2)}{1-x} \right]$$

$$= \log \left[ \lim_{x \to 1} \frac{\sin(\pi x / 2)}{1} \right]$$

We have

$$\log x \log(1-x) = \frac{\log x}{1-x} (1-x) \log(1-x)$$

$$= \frac{\log x}{1-x} \frac{\log(1-x)}{1/(1-x)}$$

Hence, using L'Hôpital's rule, we see that $\lim_{x \to 1} \log x \log(1-x) = 0$ and we determine that $\lim_{x \to 1} \log x \log \sin(\pi x) = 0$.

Using

$$\lim_{x \to 0} x \log x \log \sin(\pi x) = \lim_{x \to 0} x \log x \log \frac{\sin(\pi x)}{\pi x} + \lim_{x \to 0} x \log x \log(\pi x)$$

and
$$\lim_{x\to 0} x\log^2 x = \lim_{x\to 0} \frac{\log^2 x}{1/x} = -2\lim_{x\to 0} x\log x$$

we may easily deduce that

$$\lim_{x\to 0} x\log x \log \sin(\pi x) = 0$$

Hence we have

$$\int_0^1 x\log x \cot(\pi x)\,dx = -\frac{1}{\pi}\int_0^1 (\log x + 1)\log \sin(\pi x)\,dx$$

$$= -\frac{1}{\pi}\int_0^1 \log x \log \sin(\pi x)\,dx + \frac{1}{\pi}\log 2$$

where we have employed Euler's integral [7, p.246]

$$\int_0^1 \log \sin(\pi x)\,dx = -\log 2$$

Using the familiar trigonometric series shown in Carslaw's book [14, p.241]

$$\log[2\sin(\pi x)] = -\sum_{n=1}^\infty \frac{\cos 2n\pi x}{n} \quad (0 < x < 1)$$

and, assuming that term-by-term integration is valid, we have

$$\int_0^1 \log x \log \sin(\pi x)\,dx = -\sum_{n=1}^\infty \frac{1}{n}\int_0^1 \log x \cos 2n\pi x\,dx - \log 2 \int_0^1 \log x\,dx$$

we deduce that

$$\int_0^1 \log x \log \sin(\pi x)\,dx = \log 2 - \sum_{n=1}^\infty \frac{1}{n}\int_0^1 \log x \cos 2n\pi x\,dx$$

Integration by parts gives us

$$\int \log x \cos 2n\pi x\,dx = \frac{\sin 2n\pi x \log x - Si(2n\pi x)}{2n\pi} + c$$

and we accordingly have

$$\int_0^1 \log x \cos 2n\pi x \, dx = -\frac{Si(2n\pi)}{2n\pi}$$

Hence we have

(4.28) $$\int_0^1 \log x \log \sin(\pi x) \, dx = \log 2 + \frac{1}{2\pi} \sum_{n=1}^\infty \frac{Si(2n\pi)}{n^2}$$

and, as before, we obtain

$$\int_0^1 x \log x \cot(\pi x) \, dx = -\frac{1}{2\pi^2} \sum_{n=1}^\infty \frac{Si(2n\pi)}{n^2}$$

$$= -\frac{1}{\pi} \int_0^1 \log x \log \sin(\pi x) \, dx + \frac{1}{\pi} \log 2$$

Curiously, *WolframAlpha* evaluates the integral as $\int_0^1 \log x \log \sin(\pi x) \, dx \simeq 1.07051...$
whereas the value of the right-hand side is given as $\simeq 1.0737...$

$\square$

We have the definition (3.3)

$$Si(x) = \int_0^x \frac{\sin t}{t} \, dt = \int_0^1 \frac{\sin xu}{u} \, du$$

and the summation

$$\sum_{n=1}^\infty \frac{Si(2\pi n)}{n^2} = \sum_{n=1}^\infty \frac{1}{n^2} \int_0^1 \frac{\sin(2\pi n u)}{u} \, du$$

$$= \int_0^1 \sum_{n=1}^\infty \frac{\sin(2\pi n u)}{n^2} \frac{du}{u}$$

This gives us

$$\sum_{n=1}^\infty \frac{Si(2\pi n)}{n^2} = \int_0^1 \frac{Cl_2(2\pi u)}{u} \, du$$

where $Cl_2(x) = \sum_{n=1}^\infty \frac{\sin(nx)}{n^2}$ is the Clausen function [53, p.111].

Using the command *integrate ClausenFunction[2,2Pi x]/x from x=0 to 1*
*WolframAlpha* reports that

$$\int_0^1 \frac{\operatorname{Im} \operatorname{Li}_2(e^{2i\pi x})}{x} dx \simeq 2.39933...$$

but this also differs from $\sum_{n=1}^{\infty} \frac{Si(2n\pi)}{n^2} \simeq 2.3915...$

Witula et al. [57] obtained the following integrals

(4.29) $$\int_0^\pi \log x \log[2\sin(x/2)] dx = \sum_{n=1}^{\infty} \frac{Si(n\pi)}{n^2}$$

(4.30) $$\int_0^\pi \log x \log \cot(x/2) dx = -2 \sum_{n=1}^{\infty} \frac{Si((2n-1)\pi)}{(2n-1)^2}$$

However, *WolframAlpha* different values as shown below:

$$\int_0^\pi \log x \log[2\sin(x/2)] dx \simeq 2.83509...$$

which differs from

$$\sum_{n=1}^{\infty} \frac{Si(n\pi)}{n^2} \simeq 2.82726...$$

or

$$\int_0^{\frac{1}{2}} \log(2\pi x) \log[2\sin(\pi x)] dx = \frac{1}{2\pi} \sum_{n=1}^{\infty} \frac{Si(n\pi)}{n^2}$$

**Proposition 4.7**

(4.31) $$\int_0^1 x \log \Gamma(x) \cot(\pi x) dx = \frac{1}{2\pi} \left[ \gamma_1 + \frac{1}{2}[\varsigma(2) + \gamma^2] \right]$$

where $\gamma_1$ is the first Stieltjes constant defined by the Laurent expansion of the Hurwitz zeta function $\varsigma(s,x)$ about $s=1$

$$\varsigma(s,x) = \frac{1}{s-1} + \sum_{n=0}^{\infty} \frac{(-1)^n}{n!} \gamma_n(x)(s-1)^n$$

and $\gamma_n := \gamma_n(1)$.

**Proof**

Using (4.24) and (4.8) we obtain

$$\int_0^1 x \log \Gamma(x) \cot(\pi x)\, dx = \frac{1}{2\pi} \sum_{n=1}^{\infty} \left[ \frac{\gamma + \log n - H_n}{n} + \frac{1}{n^2} \right]$$

The related infinite series

(4.31.1) $$\sum_{n=1}^{\infty} \frac{\gamma + \log n - H_n}{n} = \gamma_1 - \frac{1}{2}[\varsigma(2) - \gamma^2]$$

was proposed as a problem by Furdui [32] and a proof was independently obtained by the author [18].

Hence, we obtain the definite integral

$$\int_0^1 x \log \Gamma(x) \cot(\pi x)\, dx = \frac{1}{2\pi}\left[ \gamma_1 + \frac{1}{2}[\varsigma(2) + \gamma^2] \right]$$

The numerical value of the right-hand side is 0.145824… in agreement with the value of the integral evaluated by *WolframAlpha*.

In passing, we note from [22] that

(4.32) $$\sum_{n=1}^{\infty} H_n \left[ \log\left(1 + \frac{1}{n}\right) - \frac{1}{n} \right] = -\left[ \gamma_1 + \frac{1}{2}[\varsigma(2) + \gamma^2] \right]$$

We express (4.31.1) as

$$\sum_{n=1}^{\infty} \frac{\gamma + \log n - H_n}{n} = \sum_{n=1}^{\infty} \left\{ \frac{\gamma + \log n}{n} - H_n \log\left(1 + \frac{1}{n}\right) + H_n \left[ \log\left(1 + \frac{1}{n}\right) - \frac{1}{n} \right] \right\}$$

$$= \sum_{n=1}^{\infty} \left\{ \frac{\gamma + \log n}{n} - H_n \log\left(1 + \frac{1}{n}\right) \right\} - \left[ \gamma_1 + \frac{1}{2}[\varsigma(2) + \gamma^2] \right]$$

$$= \gamma_1 - \frac{1}{2}[\varsigma(2) - \gamma^2]$$

Hence we see that

(4.32.1) $$\sum_{n=1}^{\infty} \left\{ \frac{\gamma + \log n}{n} - H_n \log\left(1 + \frac{1}{n}\right) \right\} = 2\gamma_1 + \gamma^2$$

□

We write

$$\int_0^1 x\log\Gamma(x)\cot(\pi x)\,dx = \int_0^1 x\log\Gamma(1+x)\cot(\pi x)\,dx - \int_0^1 x\log x\cot(\pi x)\,dx$$

and using (4.26) we obtain

$$\int_0^1 x\log\Gamma(1+x)\cot(\pi x)\,dx = \frac{1}{2\pi}\left[\gamma_1 + \frac{1}{2}[\varsigma(2)+\gamma^2]\right] + \frac{1}{2\pi^2}\sum_{n=1}^{\infty}\frac{Si(2\pi n)}{n^2}$$

**Proposition 4.8**

(4.33) $$\int_0^1 \log G(1+x)\cot(\pi x)\,dx = \frac{1}{2\pi}\left[\gamma_1 + \frac{1}{2}\gamma^2\right]$$

**Proof**

Using (4.18) we showed in equation (4.4.229u) in [16] (as subsequently corrected) that

(4.34) $$\int_0^1 [\log G(1+x) - x\log\Gamma(x)]\cot(\pi x)\,dx = -\frac{\pi}{24}$$

and using (4.31) we obtain

$$\int_0^1 \log G(1+x)\cot(\pi x)\,dx = \frac{1}{2\pi}\left[\gamma_1 + \frac{1}{2}\gamma^2\right]$$

Using the input logBarnesG[1+x]cot[Pix], *WolframAlpha* computes this integral as 0.0149245… which concurs with the above analytical evaluation.

□

In accordance with (4.18) we have

$$\int_0^1 \log G(1+x)\cot(\pi x)\,dx = 2\sum_{n=1}^{\infty}\int_0^1 G(1+x)\sin 2n\pi x\,dx$$

It may be noted that it is not valid to split the integral as

$$\int_0^1 \log G(x+1)\cot(\pi x)\,dx \neq \int_0^1 \log G(x)\cot(\pi x)\,dx + \int_0^1 \log\Gamma(x)\cot(\pi x)\,dx$$

because the individual integrals on the right-hand side do not converge (as confirmed by *WolframAlpha*). In passing, we may note that $\int_0^1 \log x \cot(\pi x)\, dx$ does not converge.

We may also note that

$$\int_0^1 \log \Gamma(x) \cot(\pi x)\, dx \neq 2\sum_{n=1}^{\infty} \int_0^1 \log \Gamma(x) \sin 2n\pi x\, dx$$

which is rather unfortunate since it is known that [26]

$$\int_0^1 \log \Gamma(x) \sin 2n\pi x\, dx = \frac{\log(2\pi n) + \gamma}{2\pi n}$$

and

$$\int_0^1 \log \Gamma(1+x) \cot(\pi x)\, dx = 2\sum_{n=1}^{\infty} \int_0^1 \log \Gamma(1+x) \sin 2n\pi x\, dx$$

Notwithstanding the above, the following representation is valid

$$\int_0^1 \log G(1+x) \cot(\pi x)\, dx = 2\sum_{n=1}^{\infty} \int_0^1 \log G(1+x) \sin 2n\pi x\, dx$$

$$= 2\sum_{n=1}^{\infty} \int_0^1 [\log G(x) + \log \Gamma(x)] \sin 2n\pi x\, dx$$

and using

$$\int_0^1 \log \Gamma(x) \sin 2n\pi x\, dx = \frac{\log(2\pi n) + \gamma}{2\pi n}$$

and

$$\int_0^1 \log G(x) \sin 2n\pi x\, dx = \frac{1}{4\pi n}\left[\frac{1}{2n} - \gamma - 2\log(2\pi) - \log n - H_n\right]$$

we obtain

$$\int_0^1 [\log G(x) + \log \Gamma(x)] \sin 2n\pi x\, dx = \frac{1}{4\pi n}\left[\frac{1}{2n} + \gamma + \log n - H_n\right]$$

Therefore we have

$$2\sum_{n=1}^{\infty} \int_0^1 [\log G(x) + \log \Gamma(x)] \sin 2n\pi x\, dx = \frac{1}{2\pi} \sum_{n=1}^{\infty}\left[\frac{1}{2n^2} + \frac{\gamma + \log n - H_n}{n}\right]$$

$$= \frac{1}{2\pi}\left[\gamma_1 + \frac{1}{2}\gamma^2\right]$$

where we have used (4.31.1). Hence, we obtain (4.33) again

$$\int_0^1 \log G(1+x)\cot(\pi x)\,dx = \frac{1}{2\pi}\left[\gamma_1 + \frac{1}{2}\gamma^2\right]$$

□

We have

$$\int_0^1 \log\sin 2n\pi x\,dx = \frac{1}{2\pi n}[Ci(2\pi n) - \log(2\pi n) - \gamma]$$

and therefore

$$\int_0^1 \log\Gamma(1+x)\cot(\pi x)\,dx = 2\sum_{n=1}^{\infty}\int_0^1 \log\Gamma(1+x)\sin 2n\pi x\,dx$$

Thus we have

(4.35) $$\int_0^1 \log\Gamma(1+x)\cot(\pi x)\,dx = \sum_{n=1}^{\infty}\frac{Ci(2\pi n)}{\pi n}$$

□

Integration by parts results in

$$\int_0^1 x\log\Gamma(x)\psi(x)\,dx = \frac{1}{2}x\log^2\Gamma(x)\Big|_0^1 - \frac{1}{2}\int_0^1 \log^2\Gamma(x)\,dx$$

and employing L'Hôpital's rule we obtain

$$\lim_{x\to 0} x\log^2\Gamma(x) = \lim_{x\to 0}\frac{\log^2\Gamma(x)}{1/x}$$

$$= \lim_{x\to 0}\frac{2\log\Gamma(x)\psi(x)}{-1/x^2}$$

$$= -2\lim_{x\to 0} x\log\Gamma(x)\lim_{x\to 0} x\psi(x)$$

$$= -2\lim_{x\to 0}\left[x\log\Gamma(1+x) - x\log x\right]\lim_{x\to 0} x\left[\psi(1+x) - \frac{1}{x}\right]$$

Therefore, the integrated part vanishes and we have

(4.36) $$\int_0^1 x \log \Gamma(x) \psi(x)\, dx = -\frac{1}{2}\int_0^1 \log^2 \Gamma(x)\, dx$$

Espinosa and Moll [31] proved that

$$\int_0^1 \log^2 \Gamma(x)\, dx = \frac{1}{4}\log(2\pi) + \frac{1}{8}\varsigma(2) + \frac{1}{12}[\gamma + \log(2\pi)]^2 - [\gamma + \log(2\pi)]\frac{\varsigma'(2)}{\pi^2} + \frac{\varsigma''(2)}{2\pi^2}$$

Since $\psi(x) - \psi(1-x) = -\pi \cot \pi x$ we have

$$\int_0^1 x \log \Gamma(x) \cot(\pi x)\, dx = \frac{1}{\pi}\int_0^1 x \log \Gamma(x)[\psi(1-x) - \psi(x)]\, dx$$

and thus

$$\int_0^1 x \log \Gamma(x) \psi(1-x)\, dx = \frac{1}{2}\left[\gamma_1 + \frac{1}{2}[\varsigma(2) + \gamma^2]\right] + \int_0^1 x \log \Gamma(x) \psi(x)\, dx$$

This results in

(4.37) $$\int_0^1 x \log \Gamma(x) \psi(1-x)\, dx = \frac{1}{2}\left[\gamma_1 + \frac{1}{2}[\varsigma(2) + \gamma^2]\right]$$

$$-\frac{1}{2}\left[\frac{1}{4}\log(2\pi) + \frac{1}{8}\varsigma(2) + \frac{1}{12}[\gamma + \log(2\pi)]^2 - [\gamma + \log(2\pi)]\frac{\varsigma'(2)}{\pi^2} + \frac{\varsigma''(2)}{2\pi^2}\right]$$

**5. Some applications of the $\Lambda(z)$ function**

The $\Lambda(z)$ function was defined in (1.3) as

$$\Lambda(z) = \sum_{j=1}^{\infty}\left[\frac{j}{j^2 + z^2} - \log\left(1 + \frac{1}{j}\right)\right]$$

As noted below, the $\Lambda(z)$ function turns out to be rather useful.

Having regard to the definition of $\Lambda(z)$, for computational convenience we define

$$\lambda(x) = \sum_{j=1}^{\infty}\left[\frac{j}{j^2 + x} - \log\left(1 + \frac{1}{j}\right)\right]$$

and differentiation gives us

$$\lambda^{(n)}(x) = (-1)^n n! \sum_{j=1}^{\infty} \frac{j}{(j^2 + x)^{n+1}}$$

We see that

$$\lambda^{(n)}(0) = (-1)^n n!\varsigma(2n+1)$$

and hence we have the Maclaurin expansion

$$\lambda(x) = \gamma + \sum_{n=1}^{\infty}(-1)^n \varsigma(2n+1)x^n$$

Noting that $\Lambda(x) = \lambda(x^2)$ we obtain

(5.1) $$\Lambda(x) = \gamma + \sum_{n=1}^{\infty}(-1)^n \varsigma(2n+1)x^{2n}$$

and hence we have

$$\Lambda^{(2n+1)}(0) = 0 \text{ for } n \geq 0$$

$$\Lambda^{(2n)}(0) = (-1)^n (2n)!\varsigma(2n+1)$$

Reference to (1.3) shows that

(5.2) $$\sum_{n=1}^{\infty}\left[\frac{n}{n^2+x^2} - \log\left(1+\frac{1}{n}\right)\right] = \gamma + \sum_{n=1}^{\infty}(-1)^n \varsigma(2n+1)x^{2n}$$

and integration gives us

(5.3) $$\sum_{n=1}^{\infty}\left[\tan^{-1}\left(\frac{x}{n}\right) - x\log\left(1+\frac{1}{n}\right)\right] = \gamma x + \sum_{n=1}^{\infty}(-1)^n \frac{\varsigma(2n+1)}{2n+1}x^{2n+1}$$

Letting $x \to ix$ in (5.3) gives us

(5.3.1) $$\sum_{n=1}^{\infty}\left[\tanh^{-1}\left(\frac{x}{n}\right) - x\log\left(1+\frac{1}{n}\right)\right] = \gamma x + \sum_{n=1}^{\infty}\frac{\varsigma(2n+1)}{2n+1}x^{2n+1}$$

Since $\tanh^{-1} x = \frac{1}{2}\log\frac{1+x}{1-x}$ we have

(5.3.2) $$\sum_{n=1}^{\infty}\left[\frac{1}{2}\log\frac{n+x}{n-x} - x\log\left(1+\frac{1}{n}\right)\right] = \gamma x + \sum_{n=1}^{\infty}\frac{\varsigma(2n+1)}{2n+1}x^{2n+1}$$

We note [48, p.130]

$$\tan^{-1}(x/u) = \int_0^{\infty}\frac{e^{-ut}\sin(xt)}{t}dt$$

Employing the Frullani integral [7, p.98]

$$\log(a) = \int_0^\infty \frac{e^{-t} - e^{-at}}{t} dt$$

we see that

$$\log\left(1 + \frac{1}{n}\right) = \int_0^\infty \frac{e^{-nt} - e^{-(n+1)t}}{t} dt$$

We may then write (5.3) as

$$\sum_{n=1}^\infty \left[\tan^{-1}\left(\frac{x}{n}\right) - x \log\left(1 + \frac{1}{n}\right)\right] = \sum_{n=1}^\infty \left[\int_0^\infty \frac{e^{-nt} \sin(xt)}{t} dt - x \int_0^\infty \frac{e^{-nt}[1 - e^{-t}]}{t} dt\right]$$

$$= \int_0^\infty \frac{\sin(xt) - x[1 - e^{-t}]}{t(e^t - 1)} dt$$

$$= \int_0^\infty \left[\frac{\sin(xt)}{t(e^t - 1)} - \frac{x}{te^t}\right] dt$$

Therefore we obtain

(5.4) $$\int_0^\infty \left[\frac{\sin(xt)}{t(e^t - 1)} - \frac{x}{te^t}\right] dt = \gamma x + \sum_{n=1}^\infty (-1)^n \frac{\varsigma(2n+1)}{2n+1} x^{2n+1}$$

$$= -\frac{1}{2i}[\log \Gamma(1 + ix) - \log \Gamma(1 - ix)]$$

where in the last part we have used [53, p.261] (and see (5.40) below).

We note that Candelpergher [13, p.35] showed that

(5.4.1) $$\sum_{n=1}^\infty \left[\tan^{-1}\left(\frac{x}{n}\right) - \frac{x}{n}\right] = -\gamma x - \frac{1}{2i}[\log \Gamma(1 + ix) - \log \Gamma(1 - ix)]$$

which may also be deduced from (5.38) below.

With regard to (5.4) we note a connection with Legendre's relation [56, p.119]

(5.4.2) $$2\int_0^\infty \frac{\sin(xt)}{e^{2\pi t} - 1} dt = \frac{1}{e^x - 1} - \frac{1}{x} + \frac{1}{2} = \frac{1}{2} \coth \frac{x}{2} - \frac{1}{x}$$

A rigorous derivation of this result is shown in Bromwich's book [12, p.501].

Differentiation of (5.4) gives us

(5.5) $$\int_0^\infty \left[\frac{\cos(xt)}{e^t-1} - \frac{1}{te^t}\right] dt = \gamma + \sum_{n=1}^\infty (-1)^n \varsigma(2n+1) x^{2n}$$

$$= \sum_{j=1}^\infty \left[\frac{j}{j^2+x^2} - \log\left(1+\frac{1}{j}\right)\right]$$

See also [53, p.260].

Reference to (1.4.1) then shows that

(5.6) $$\int_0^\infty \left[\frac{\cos(xt)}{e^t-1} - \frac{1}{te^t}\right] dt = -\frac{1}{2}[\psi(1+ix) + \psi(1-ix)]$$

We split the following integral

$$\int_0^\infty \left[\frac{1}{e^t-1} - \frac{1}{t}\right] \cos(xt) dt = \int_0^\infty \left[\frac{\cos(xt)}{e^t-1} - \frac{1}{te^t}\right] dt + \int_0^\infty \left[\frac{e^{-t} - \cos(xt)}{t}\right] dt$$

From Nahin's book [47, p.141] we have

(5.6.1) $$\int_0^\infty \left[\frac{e^{-rt}\cos(pt) - e^{-st}\cos(qt)}{t}\right] dt = \frac{1}{2} \log \frac{q^2+s^2}{p^2+r^2}$$

and with $p = s = 0$ and $r = 1$ we obtain

(5.6.2) $$\int_0^\infty \left[\frac{e^{-t} - \cos(xt)}{t}\right] dt = \log x$$

Hence, for $x > 0$ we end up with

(5.7) $$\int_0^\infty \left[\frac{1}{e^t-1} - \frac{1}{t}\right] \cos(xt) dt = \log x - \frac{1}{2}[\psi(1+ix) + \psi(1-ix)]$$

which appears in [36, p.505, Eq. (3.951.6)].

With $x = 0$ in (5.6) we obtain

$$\int_0^\infty \left[\frac{1}{e^t-1} - \frac{1}{te^t}\right] dt = \gamma$$

and this integral, due to Gauss, appears in [7, p.177].

The first derivative of (5.5) gives us

$$\int_0^\infty \frac{t\sin(xt)}{e^t-1}\,dt = \sum_{n=1}^\infty (-1)^{n+1} 2n\varsigma(2n+1)x^{2n-1}$$

and the second derivative of (5.5) results in

$$\int_0^\infty \frac{t^2\cos(xt)}{e^t-1}\,dt = \sum_{n=1}^\infty (-1)^{n+1} 2n(2n-1)\varsigma(2n+1)x^{2n-2}$$

Letting $x=0$ gives us

$$\int_0^\infty \frac{t^2}{e^t-1}\,dt = 2\varsigma(3)$$

The $2k\,th$ derivative gives us the well-known integral [7, p.223]

$$\int_0^\infty \frac{t^{2k}}{e^t-1}\,dt = (2k)!\varsigma(2k+1)$$

Dividing (5.4) by $x$ and integrating over $[0,1]$ gives us

$$\int_0^\infty \left[\frac{Si(t)}{t(e^t-1)} - \frac{1}{te^t}\right]dt = \gamma + \sum_{n=1}^\infty (-1)^n \frac{\varsigma(2n+1)}{(2n+1)^2}$$

where we have used the sine integral defined in (3.3) as

$$Si(t) = \int_0^t \frac{\sin u}{u}\,du = \int_0^1 \frac{\sin xt}{x}\,dx$$

We note from (5.2) that

(5.8) $$\sum_{n=1}^\infty \left[\frac{n}{n^2+x} - \log\left(1+\frac{1}{n}\right)\right] = \gamma + \sum_{n=1}^\infty (-1)^n \varsigma(2n+1)x^n$$

With $x \to -x$ we have

(5.9) $$\sum_{n=1}^\infty \left[\frac{n}{n^2-x} - \log\left(1+\frac{1}{n}\right)\right] = \gamma + \sum_{n=1}^\infty \varsigma(2n+1)x^n$$

and $x \to x^2$ results in

(5.10) $$\sum_{n=1}^\infty \left[\frac{n}{n^2-x^2} - \log\left(1+\frac{1}{n}\right)\right] = \gamma + \sum_{n=1}^\infty \varsigma(2n+1)x^{2n}$$

We note from [53, p.160] that we have in terms of the Hurwitz zeta function $\varsigma(s,a)$

(5.11) $$\sum_{n=1}^{\infty} \varsigma(2n+1,a)x^{2n} = -\frac{1}{2}[\psi(a+x)+\psi(a-x)]+\psi(a)$$

and thus with $a=1$ we have

(5.12) $$\gamma + \sum_{n=1}^{\infty} \varsigma(2n+1)x^{2n} = -\frac{1}{2}[\psi(1+x)+\psi(1-x)]$$

Hence we see that

(5.13) $$\sum_{n=1}^{\infty}\left[\frac{n}{n^2-x^2} - \log\left(1+\frac{1}{n}\right)\right] = -\frac{1}{2}[\psi(1+x)+\psi(1-x)]$$

which may be deduced directly from (1.4).

We note that *Mathematica* is able to determine a closed form for the summation with $x = \frac{1}{2}$.

$$\sum_{n=1}^{\infty}\left(\frac{n}{n^2 - \frac{1}{4}} - \log\left(\frac{1}{n}+1\right)\right) = -1 + \gamma + \log(4) \approx 0.96351$$

□

We recall [53, p.14]

$$\psi(1+x)-\psi(1-x) = \frac{1}{x} - \pi \cot \pi x$$

and adding this to (5.13) gives us

$$\frac{1}{x} - \pi \cot \pi x - 2\sum_{n=1}^{\infty}\left[\frac{n}{n^2-x^2} - \log\left(1+\frac{1}{n}\right)\right] = 2\psi(1+x)$$

Substituting (1.4)

$$\psi(1+x) = \sum_{n=1}^{\infty}\left[\log\left(1+\frac{1}{n}\right) - \frac{1}{n+x}\right]$$

gives us

$$\pi \cot \pi x = \frac{1}{x} - 2\sum_{n=1}^{\infty}\left[\frac{n}{n^2-x^2} - \log\left(1+\frac{1}{n}\right)\right] - 2\sum_{n=1}^{\infty}\left[\log\left(1+\frac{1}{n}\right) - \frac{1}{n+x}\right]$$

$$= \frac{1}{x} - 2\sum_{n=1}^{\infty}\left[\frac{n}{n^2-x^2} - \frac{1}{n+x}\right]$$

Hence we obtain the well-known decomposition formula for $\cot \pi x$ [53, p.30]

$$\pi \cot \pi x = \frac{1}{x} - 2\sum_{n=1}^{\infty} \frac{x}{n^2 - x^2}$$

□

Integration of (5.9) gives us

$$-\sum_{n=1}^{\infty}\left[n\log\left(1-\frac{x}{n^2}\right) + x\log\left(1+\frac{1}{n}\right)\right] = \gamma x + \sum_{n=1}^{\infty} \frac{\varsigma(2n+1)}{n+1} x^{n+1}$$

and with $x \to x^2$ we have

(5.14) $\quad -\sum_{n=1}^{\infty}\left[n\log\left(1-\frac{x^2}{n^2}\right) + x^2 \log\left(1+\frac{1}{n}\right)\right] = \gamma x^2 + \sum_{n=1}^{\infty} \frac{\varsigma(2n+1)}{n+1} x^{2n+2}$

The Barnes double gamma function $G(x)$ may be defined by [53, p.25]

(5.15) $\quad \log G(1+x) = \frac{1}{2}x\log(2\pi) - \frac{1}{2}x(1+x) - \frac{1}{2}\gamma x^2 + \sum_{n=1}^{\infty}\left[\frac{1}{2n}x^2 - x + n\log\left(1+\frac{x}{n}\right)\right]$

This easily gives us

$$\log[G(1+x)G(1-x)] = -(1+\gamma)x^2 + \sum_{n=1}^{\infty}\left[\frac{1}{n}x^2 + n\log\left(1-\frac{x^2}{n^2}\right)\right]$$

and we write

$$\sum_{n=1}^{\infty}\left[n\log\left(1-\frac{x^2}{n^2}\right) + x^2\log\left(1+\frac{1}{n}\right)\right] = \sum_{n=1}^{\infty}\left[\frac{1}{n}x^2 + n\log\left(1-\frac{x^2}{n^2}\right)\right]$$

$$-x^2 \sum_{n=1}^{\infty}\left[\frac{1}{n} - \log\left(1+\frac{1}{n}\right)\right]$$

and substituting $\gamma = \sum_{n=1}^{\infty}\left[\frac{1}{n} - \log\left(1+\frac{1}{n}\right)\right]$ we obtain

(5.16.1) $\quad \sum_{n=1}^{\infty}\left[n\log\left(1-\frac{x^2}{n^2}\right) + x^2\log\left(1+\frac{1}{n}\right)\right] = \sum_{n=1}^{\infty}\left[\frac{1}{n}x^2 + n\log\left(1-\frac{x^2}{n^2}\right)\right] - \gamma x^2$

(5.16.2) $\hspace{6cm} = x^2 + \log[G(1+x)G(1-x)]$

(5.16.3) $\hspace{6cm} = -\gamma x^2 - \sum_{n=1}^{\infty} \frac{\varsigma(2n+1)}{n+1} x^{2n+2}$

Hence, we obtain for $|x|<1$

$$\text{(5.17)} \qquad \sum_{n=1}^{\infty} \frac{\varsigma(2n+1)}{n+1} x^{2n+2} = -(1+\gamma)x^2 - \log[G(1+x)G(1-x)]$$

which was reported in [53, p.212, Eq. (465)]. This may also be deduced from the Taylor expansion which may be constructed using (5.15)

$$\text{(5.17.1)} \qquad \log G(1+x) = \frac{1}{2}[\log(2\pi)-1]x - \frac{1}{2}(1+\gamma)x^2 + \sum_{n=2}^{\infty}(-1)^n \frac{\varsigma(n)}{n+1} x^{n+1}$$

*WolframAlpha* provides the following partial summation of (5.16.1) for $x = \frac{1}{2}$

**Partial sum formula**

$$\sum_{n=1}^{k} \left( n \log\left(1 - \frac{1}{4n^2}\right) + \frac{1}{4}\log\left(1 + \frac{1}{n}\right) \right) =$$
$$\frac{1}{12}\left(12\zeta^{(1,0)}\left(-1, k+\frac{1}{2}\right) - 24\zeta^{(1,0)}(-1, k+1) + \right.$$
$$12\zeta^{(1,0)}\left(-1, k+\frac{3}{2}\right) - 36\log(A) + 6\log\left(\Gamma\left(k+\frac{1}{2}\right)\right) -$$
$$\left. 3\log(\Gamma(k+1)) - 6\log\left(\Gamma\left(k+\frac{3}{2}\right)\right) + 3\log(\Gamma(k+2)) + 3 + \log(2)\right)$$

and, from the structure of the output, we might assume that the only surviving terms as $k \to \infty$ are

$$\text{(5.18)} \qquad \sum_{n=1}^{\infty}\left[ n\log\left(1-\frac{1}{4n^2}\right) + \frac{1}{4}\log\left(1+\frac{1}{n}\right)\right] = -3\log A + \frac{1}{4} + \frac{1}{12}\log 2 + \lim_{k\to\infty} R_k$$

where $\log A = \frac{1}{12} - \varsigma'(-1)$

and

$$R_k = \varsigma'\left(-1, k+\frac{1}{2}\right) + \varsigma'\left(-1, k+\frac{3}{2}\right) - 2\varsigma'(-1, k+1)$$

$$+ \frac{1}{2}\log\Gamma\left(k+\frac{1}{2}\right) - \frac{1}{2}\log\Gamma\left(k+\frac{3}{2}\right) - \frac{1}{4}\log\Gamma(k+1) + \frac{1}{4}\log\Gamma(k+2)$$

It should be possible to compute $R_k$ by using (4.13)

$$\log G(1+x) - x\log\Gamma(x) = \varsigma'(-1) - \varsigma'(-1, x)$$

*WolframAlpha* also computes the numerical value

$$\sum_{n=1}^{\infty}\left[n\log\left(1-\frac{1}{4n^2}\right)+\frac{1}{4}\log\left(1+\frac{1}{n}\right)\right]\simeq -0.187878...$$

It is known [53, p.26] that

$$\log G(\tfrac{1}{2}) = \frac{1}{24}\log 2 - \frac{1}{4}\log\pi + \frac{3}{2}\varsigma'(-1)$$

and, since $G(1+x) = G(x)\Gamma(x)$, we have

$$\log G(\tfrac{3}{2}) = \frac{1}{24}\log 2 + \frac{1}{4}\log\pi + \frac{3}{2}\varsigma'(-1)$$

and substituting these in (5.16.2) produces

(5.18) $$\sum_{n=1}^{\infty}\left[n\log\left(1-\frac{1}{4n^2}\right)+\frac{1}{4}\log\left(1+\frac{1}{n}\right)\right] = 3\varsigma'(-1)+\frac{1}{4}+\frac{1}{12}\log 2$$

An alternative derivation of (5.16.2) is set out below.

We multiply (5.13) by $x$ and integrate to obtain

$$\int_0^u \sum_{n=1}^{\infty}\left[\frac{nx}{n^2-x^2} - x\log\left(1+\frac{1}{n}\right)\right]dx = -\frac{1}{2}\int_0^u x[\psi(1+x)+\psi(1-x)]dx$$

An elementary derivation of the following integral was given in [24]

(5.19) $$\log G(1+u) = \frac{1}{2}u\log(2\pi) - \frac{1}{2}u(u+1) + \int_0^u x\psi(1+x)\,dx$$

and an obvious substitution shows that

$$\int_0^u x\psi(1-x)\,dx = \int_0^{-u} x\psi(1+x)\,dx$$

Hence we see that

$$\int_0^u x[\psi(1+x)+\psi(1-x)]dx = u^2 + \log[G(1+u)G(1-u)]$$

and we obtain another derivation of (5.16.2).

$$\sum_{n=1}^{\infty}\left[n\log\left(1-\frac{x^2}{n^2}\right)+x^2\log\left(1+\frac{1}{n}\right)\right] = x^2 + \log[G(1+x)G(1-x)]$$

Equating (5.15) and (5.17.1) results in

$$\sum_{n=2}^{\infty}(-1)^n \frac{\varsigma(n)}{n+1} x^{n+1} = \sum_{n=1}^{\infty}\left[\frac{1}{2n}x^2 - x + n\log\left(1+\frac{x}{n}\right)\right]$$

□

Kinkelin (see [24]) showed that for $|x|<1$

(5.20) $\quad \int_0^x \pi t \cot \pi t \, dt = x\log(2\pi) + \log G(1-x) - \log G(1+x)$

and substituting (5.17) we see that

(5.21) $\quad \int_0^x \pi t \cot \pi t \, dt = -\sum_{n=1}^{\infty}\frac{\varsigma(2n+1)}{n+1}x^{2n+2} + x\log(2\pi) - (1+\gamma)x^2 - 2\log G(1+x)$

We showed in [24] that for $|x|<1$

(5.22) $\quad \int_0^x \pi t \cot \pi t \, dt = \frac{1}{2\pi}\sum_{n=1}^{\infty}\frac{\sin 2n\pi x}{n^2} - x\sum_{n=1}^{\infty}\frac{\cos 2n\pi x}{n}$

and we therefore obtain

(5.23) $\quad \frac{1}{2\pi}\sum_{n=1}^{\infty}\frac{\sin 2n\pi x}{n^2} - x\sum_{n=1}^{\infty}\frac{\cos 2n\pi x}{n} = x\log(2\pi) + \log G(1-x) - \log G(1+x)$

**Proposition 5.1**

(5.24) $\quad (x+\tfrac{1}{2})\log(2\pi) + \log G(\tfrac{1}{2}-x) - \log G(\tfrac{3}{2}+x)$

$$= \frac{1}{2}[\log G(1-2x) - \log G(1+2x)] - [\log G(1-x) - \log G(1+x)] + \frac{1}{2}\log(2\cos \pi x)$$

**Proof**

We define

$$C_s(x) := \sum_{n=1}^{\infty}\frac{\cos 2n\pi x}{n^s}$$

$$S_s(x) := \sum_{n=1}^{\infty}\frac{\sin 2n\pi x}{n^s}$$

and we may write (5.23) as

$$\frac{1}{2\pi}S_2(x) - xC_1(x) = x\log(2\pi) + \log G(1-x) - \log G(1+x)$$

It is easily seen that

$$C_s(x+\tfrac{1}{2}) = 2^{1-s}C_s(2x) - C_s(x)$$
$$S_s(x+\tfrac{1}{2}) = 2^{1-s}S_s(2x) - S_s(x)$$

We have

$$\frac{1}{2\pi}S_2(x+\tfrac{1}{2}) - (x+\tfrac{1}{2})C_1(x+\tfrac{1}{2}) = (x+\tfrac{1}{2})\log(2\pi) + \log G(\tfrac{1}{2}-x) - \log G(\tfrac{3}{2}+x)$$

and

$$\frac{1}{2\pi}S_2(x+\tfrac{1}{2}) - C_1(x+\tfrac{1}{2}) = \frac{1}{2\pi}\left[\frac{1}{2}S_2(2x) - S_2(x)\right] - (x+\tfrac{1}{2})\left[C_1(2x) - C_1(x)\right]$$

$$= \frac{1}{2}\left[\frac{1}{2\pi}S_2(2x) - 2xC_1(2x)\right] - \left[\frac{1}{2\pi}S_2(x) - xC_1(x)\right] - \tfrac{1}{2}\left[C_1(2x) - C_1(x)\right]$$

$$= \frac{1}{2}[2x\log(2\pi) + \log G(1-2x) - \log G(1+2x)] - [x\log(2\pi) + \log G(1-x) - \log G(1+x)]$$

$$-\tfrac{1}{2}\left[C_1(2x) - C_1(x)\right]$$

$$= \frac{1}{2}[\log G(1-2x) - \log G(1+2x)] - [\log G(1-x) - \log G(1+x)]$$

$$-\tfrac{1}{2}\left[C_1(2x) - C_1(x)\right]$$

We note [54, p.148]

$$C_1(x) = \sum_{n=1}^{\infty}\frac{\cos 2n\pi x}{n} = -\log[2\sin \pi x]$$

which results in

$$C_1(2x) - C_1(x) = -\log[2\cos \pi x]$$

Hence we obtain

$$(x+\tfrac{1}{2})\log(2\pi) + \log G(\tfrac{1}{2}-x) - \log G(\tfrac{3}{2}+x)$$

$$= \frac{1}{2}[\log G(1-2x) - \log G(1+2x)] - [\log G(1-x) - \log G(1+x)] + \frac{1}{2}\log(2\cos \pi x)$$

which may be compared with [53, p.29].

$\square$

Letting $x \to x+\tfrac{1}{2}$ in (5.23) we have for $-\tfrac{3}{2} < x < \tfrac{1}{2}$

$$\frac{1}{2\pi}\sum_{n=1}^{\infty}(-1)^n \frac{\sin 2n\pi x}{n^2} - (x+\tfrac{1}{2})\sum_{n=1}^{\infty}(-1)^n \frac{\cos 2n\pi x}{n}$$

$$= (x+\tfrac{1}{2})\log(2\pi) + \log G(\tfrac{1}{2}-x) - \log G(\tfrac{3}{2}+x)$$

Letting $x \to x-\tfrac{1}{2}$ we have for $-\tfrac{1}{2} < x < \tfrac{3}{2}$

$$\frac{1}{2\pi}\sum_{n=1}^{\infty}(-1)^n \frac{\sin 2n\pi x}{n^2} - (x-\tfrac{1}{2})\sum_{n=1}^{\infty}(-1)^n \frac{\cos 2n\pi x}{n}$$

$$= (x-\tfrac{1}{2})\log(2\pi) + \log G(\tfrac{3}{2}-x) - \log G(\tfrac{1}{2}+x)$$

Subtracting the above two equations gives us for $-\tfrac{1}{2} < x < \tfrac{1}{2}$

$$\sum_{n=1}^{\infty}(-1)^{n+1}\frac{\cos 2n\pi x}{n} = \log(2\pi) + \log G(\tfrac{1}{2}-x) - \log G(\tfrac{3}{2}+x)$$

$$-\log G(\tfrac{3}{2}-x) + \log G(\tfrac{1}{2}+x)$$

We recall [54, p.148] that for $|x| < \tfrac{1}{2}$

$$\sum_{n=1}^{\infty}(-1)^{n+1}\frac{\cos 2n\pi x}{n} = \log[2\cos \pi x]$$

and we obtain

$$\log \cos \pi x = \log \pi + \log G(\tfrac{1}{2}-x) - \log G(\tfrac{3}{2}+x) - \log G(\tfrac{3}{2}-x) + \log G(\tfrac{1}{2}+x)$$

or

$$\log \frac{\cos \pi x}{\pi} = \log \frac{G(\tfrac{1}{2}+x)G(\tfrac{1}{2}-x)}{G(\tfrac{3}{2}+x)G(\tfrac{3}{2}-x)}$$

Since
$$G(\tfrac{3}{2}+x) = G(1+\tfrac{1}{2}+x)$$

$$= G(\tfrac{1}{2}+x)\Gamma(\tfrac{1}{2}+x)$$

we simply end up with a variant of Euler's reflection formula

$$\Gamma(\tfrac{1}{2}+x)\Gamma(\tfrac{1}{2}-x) = \frac{\pi}{\cos \pi x}$$

□

Letting $x = \tfrac{1}{4}$ in (5.22) gives us

$$\int_0^{\frac{1}{4}} \pi t \cot \pi t \, dt = \frac{1}{2\pi}\sum_{n=1}^{\infty}\frac{\sin(n\pi/2)}{n^2} - \frac{1}{4}\sum_{n=1}^{\infty}\frac{\cos(n\pi/2)}{n}$$

We have for suitably convergent series

$$\sum_{n=1}^{\infty} a_n \cos(n\pi/2) = \sum_{n=1}^{\infty} a_{2n} \cos(n\pi)$$

$$\sum_{n=1}^{\infty} a_n \sin(n\pi/2) = \sum_{n=0}^{\infty} a_{2n+1} \cos(n\pi)$$

and thus

$$\int_0^{\frac{1}{4}} \pi t \cot \pi t \, dt = \frac{1}{2\pi}\sum_{n=0}^{\infty}\frac{(-1)^n}{(2n+1)^2} - \frac{1}{8}\sum_{n=1}^{\infty}\frac{(-1)^n}{n}$$

The Dirichlet beta function $\beta(s)$ is defined for $\operatorname{Re}(s) > 0$ by

$$\beta(s) = \sum_{n=0}^{\infty}\frac{(-1)^n}{(2n+1)^s}$$

and we have $\beta(2) = G$, known as Catalan's constant. This gives us

$$\int_0^{\frac{1}{4}} \pi t \cot \pi t \, dt = \frac{G}{2\pi} + \frac{\log 2}{8}$$

Using (5.20) we have

$$\int_0^{\frac{1}{4}} \pi t \cot \pi t \, dtx = \frac{1}{4}\log(2\pi) + \log G(\tfrac{3}{4}) - \log G(\tfrac{5}{4})$$

and we thereby see that

$$\log G(\tfrac{3}{4}) - \log G(\tfrac{5}{4}) = \frac{G}{2\pi} - \frac{1}{8}\log 2 - \frac{1}{4}\log \pi$$

in agreement with [53, p.30].

## Some further aspects of the digamma function

We note the following expansion for the digamma function [7, p.212]

$$\psi(z) = -\gamma - \sum_{j=0}^{\infty}\left[\frac{1}{j+z} - \frac{1}{j+1}\right]$$

and thus we have

$$\psi(u+iv) = -\gamma - \sum_{j=0}^{\infty}\left[\frac{j+u-iv}{(j+u)^2+v^2} - \frac{1}{j+1}\right]$$

We see that

$$\psi(1+iv) = -\gamma - \sum_{j=1}^{\infty}\left[\frac{j}{j^2+v^2} - \frac{1}{j}\right] + i\sum_{j=1}^{\infty}\frac{v}{j^2+v^2}$$

or equivalently

(5.24.1) $$\psi(1+iv) = -\gamma + \sum_{j=1}^{\infty}\frac{v^2}{j(j^2+v^2)} + i\sum_{j=1}^{\infty}\frac{v}{j^2+v^2}$$

We define $f(x)$ as

$$f(x) = \sum_{n=1}^{\infty}\frac{1}{x+n^2}$$

and we have the derivative

$$f^{(k)}(u) = (-1)^k k!\sum_{n=1}^{\infty}\frac{1}{(x+n^2)^{k+1}}$$

so that

$$f^{(k)}(0) = (-1)^k k!\sum_{n=1}^{\infty}\frac{1}{n^{2k+2}} = (-1)^k k!\varsigma(2k+2)$$

We therefore have the Maclaurin series

$$f(x) = \sum_{n=1}^{\infty}\frac{1}{x+n^2} = \sum_{k=0}^{\infty}(-1)^k \varsigma(2k+2)x^k$$

and we see that

$$\sum_{k=0}^{\infty}(-1)^k \varsigma(2k+2)x^k = \sum_{k=1}^{\infty}(-1)^{k+1}\varsigma(2k)x^{k-1}$$

Hence we obtain

(5.24.2) $$\sum_{n=1}^{\infty}\frac{1}{x+n^2} = \sum_{k=1}^{\infty}(-1)^{k+1}\varsigma(2k)x^{k-1}$$

Letting $x \to x^2$ gives us

(5.24.3) $$\sum_{n=1}^{\infty}\frac{1}{x^2+n^2} = \sum_{n=1}^{\infty}(-1)^{n+1}\varsigma(2n)x^{2n-2}$$

We let $x \to t/2\pi$ to obtain

$$\sum_{n=1}^{\infty}\frac{1}{t^2+4\pi^2 n^2} = \sum_{n=1}^{\infty}(-1)^{n+1}(2\pi)^{-2n}\varsigma(2n)t^{2n-2}$$

It is well known that [4, p.266]

$$\varsigma(2n) = \frac{(-1)^{n+1}(2\pi)^{2n}B_{2n}}{2(2n)!}$$

where $B_n$ are the Bernoulli numbers. Therefore, we obtain

$$2\sum_{n=1}^{\infty}\frac{1}{4\pi^2 t^2+n^2} = \sum_{n=1}^{\infty}\frac{B_{2n}}{(2n)!}t^{2n-2}$$

The Bernoulli numbers $B_n$ are given by the generating function

$$\frac{t}{e^t-1} = \sum_{n=0}^{\infty}B_n \frac{t^n}{n!} \qquad , (|t|<2\pi)$$

Therefore, employing $B_0 = 1$, $B_1 = -\frac{1}{2}$ and $B_{2n+1} = 0$ for all $n \geq 1$, we have

$$\frac{1}{e^x-1} = \frac{1}{x} - \frac{1}{2} + \sum_{n=1}^{\infty}B_{2n}\frac{x^{2n-1}}{(2n)!}$$

and hence we obtain the well-known identity ([12, p.296] and [39, p.378])

(5.24.4) $$2\sum_{n=1}^{\infty}\frac{x}{x^2+4\pi^2 n^2} = \frac{1}{e^x-1} - \frac{1}{x} + \frac{1}{2}$$

$$= \frac{1}{2}\coth\left(\frac{1}{2}x\right) - \frac{1}{x}$$

so that with $x = 2\pi v$ we have

$$\sum_{n=1}^{\infty} \frac{v}{v^2 + n^2} = \frac{\pi}{2}\coth(\pi v) - \frac{1}{2v}$$

Therefore, referring to (5.24.1) we have

(5.25) $$\psi(1+iv) = -\gamma - \sum_{j=1}^{\infty}\left[\frac{j}{j^2+v^2} - \frac{1}{j}\right] + i\left[\frac{\pi}{2}\coth(\pi v) - \frac{1}{2v}\right]$$

or equivalently

(5.26) $$\psi(1+iv) = -\gamma + \sum_{j=1}^{\infty}\frac{v^2}{j(j^2+v^2)} + i\left[\frac{\pi}{2}\coth(\pi v) - \frac{1}{2v}\right]$$

This formula appears in the well-known handbook [1, p.259] and this representation was noted in Bradley's paper [11].

Letting $v \to -iv$ in (5.26) and noting that $\coth(ix) = -i\cot x$ gives us

(5.26.1) $$\psi(1+v) = -\gamma - \sum_{j=1}^{\infty}\frac{v^2}{j(j^2-v^2)} - \frac{\pi}{2}\cot(\pi v) + \frac{1}{2v}$$

or equivalently

(5.26.2) $$\psi(1+v) = -\gamma - \sum_{j=1}^{\infty}\frac{v^2}{j(j^2-v^2)} + \sum_{j=1}^{\infty}\frac{v}{j^2-v^2}$$

$$= -\gamma + \sum_{j=1}^{\infty}\frac{v}{j(j+v)}$$

Since $\psi(1+z) = \psi(z) + \frac{1}{z}$ we obtain from (5.25)

$$\psi(iv) = -\gamma - \sum_{j=1}^{\infty}\left[\frac{j}{j^2+v^2} - \frac{1}{j}\right] + i\left[\frac{\pi}{2}\coth(\pi v) + \frac{1}{2v}\right]$$

As noted by Dixit [28] we see that

$$\sum_{j=1}^{n}\left[\frac{j}{j^2+v^2} - \frac{1}{j}\right] = \sum_{j=1}^{n}\frac{j}{j^2+v^2} - \log n + \log n - \sum_{j=1}^{n}\frac{1}{j}$$

and we have the limit as $n \to \infty$

$$\sum_{j=1}^{\infty}\left[\frac{j}{j^2+v^2}-\frac{1}{j}\right]=\lim_{n\to\infty}\left(\sum_{j=1}^{n}\frac{j}{j^2+v^2}-\log n\right)-\gamma$$
$$=\Lambda(v)-\gamma$$

We see that

$$\psi(-iv)=-\gamma-\sum_{j=1}^{\infty}\left[\frac{j}{j^2+v^2}-\frac{1}{j}\right]-i\left[\frac{\pi}{2}\coth(\pi v)+\frac{1}{2v}\right]$$

and thus we have

$$\psi(iv)+\psi(-iv)=-2\gamma-2\sum_{j=1}^{\infty}\left[\frac{j}{j^2+v^2}-\frac{1}{j}\right]$$

We also have

$$\Lambda(v)=-\frac{1}{2}[\psi(iv)+\psi(-iv)]$$

$$=-\frac{1}{2}[\psi(1+iv)+\psi(1-iv)]$$

which was previously obtained by Dixit [28]. We see that

$$\psi(1+v)+\psi(1-v)=-2\Lambda(iz)$$

and thus

$$\psi(1+v)+\psi(1-v)=-2\sum_{n=1}^{\infty}\left[\frac{n}{n^2-v^2}-\log\left(1+\frac{1}{n}\right)\right]$$

provided $v\neq n$. An alternative derivation of this is shown below.

We recall (2.8)

$$\psi(1+v)+\psi(1-v)+2\gamma=-2\sum_{n=1}^{\infty}\frac{1}{n}\frac{v^2}{n^2-v^2}$$

and substituting

$$\gamma=\sum_{n=1}^{\infty}\left[\frac{1}{n}-\log\left(1+\frac{1}{n}\right)\right]$$

we obtain

$$\psi(1+v)+\psi(1-v) = -2\sum_{n=1}^{\infty}\left[\frac{1}{n}\frac{v^2}{n^2-v^2}+\frac{1}{n}-\log\left(1+\frac{1}{n}\right)\right]$$

$$= -2\sum_{n=1}^{\infty}\left[\frac{n}{n^2-v^2}-\log\left(1+\frac{1}{n}\right)\right]$$

$\square$

We define $g(u)$ as

$$g(u) = \sum_{j=1}^{\infty}\frac{1}{j(j^2+u)}$$

and we have the derivative

$$g^{(k)}(u) = (-1)^k k!\sum_{j=1}^{\infty}\frac{1}{j(j^2+u)^{k+1}}$$

so that

$$g^{(k)}(0) = (-1)^k k!\sum_{j=1}^{\infty}\frac{1}{j^{2k+3}} = (-1)^k k!\varsigma(2k+3)$$

We therefore have the Maclaurin series

$$g(u) = \sum_{j=1}^{\infty}\frac{1}{j(j^2+u)} = \sum_{k=0}^{\infty}(-1)^k \varsigma(2k+3)u^k$$

and we see that

$$\sum_{n=1}^{\infty}(-1)^{n+1}\varsigma(2n+1)u^{n-1} = \sum_{n=0}^{\infty}(-1)^n \varsigma(2n+3)u^n$$

Hence we obtain

(5.27) $$\sum_{j=1}^{\infty}\frac{1}{j(j^2+u)} = \sum_{n=1}^{\infty}(-1)^{n+1}\varsigma(2n+1)u^{n-1}$$

Letting $u = x^2$ gives us

(5.28) $$\sum_{n=1}^{\infty}(-1)^{n+1}\varsigma(2n+1)x^{2n} = \sum_{j=1}^{\infty}\frac{x^2}{j(j^2+x^2)}$$

which, I subsequently noted, appears in [7].

We recall (5.24.3)

$$\sum_{n=1}^{\infty}\frac{x^2}{x^2+n^2}=\sum_{n=1}^{\infty}(-1)^{n+1}\varsigma(2n)x^{2n}$$

If we let $x \to ix$ and combine that with (5.28) we end up with the Maclaurin expansion shown in (5.26.2)

$$\psi(1+x)=-\gamma-\sum_{j=1}^{\infty}\frac{x^2}{j(j^2-x^2)}+\sum_{j=1}^{\infty}\frac{x}{j^2-x^2}$$

With $x \to ix$ in (5.28) we obtain

(5.29) $$\sum_{n=1}^{\infty}\varsigma(2n+1)x^{2n}=\sum_{j=1}^{\infty}\frac{x^2}{j(j^2-x^2)}$$

and we recall (5.26.1)

$$\psi(1+x)=-\gamma-\sum_{j=1}^{\infty}\frac{x^2}{j(j^2-x^2)}-\frac{\pi}{2}\cot(\pi x)+\frac{1}{2x}$$

which then gives us [53, p.260]

(5.30) $$\psi(1+x)=\frac{1}{2x}-\frac{\pi}{2}\cot(\pi x)-\gamma-\sum_{n=1}^{\infty}\varsigma(2n+1)x^{2n}$$

Nielsen ([48] and [36, p.893, 8.363.2]) obtained the equivalent result

(5.31) $$\psi(1+x)=\frac{1}{2x}-\frac{\pi}{2}\cot(\pi x)-\gamma-\frac{x^2}{1-x^2}-\sum_{n=1}^{\infty}[\varsigma(2n+1)-1]x^{2n}$$

which was rediscovered by Zhang and Williams [53, p.260] in 1993.

Integration of (5.30 results in [53, p.261]

(5.32) $$\sum_{n=1}^{\infty}\frac{\varsigma(2n+1)}{2n+1}x^{2n+1}=\frac{1}{2}\log\frac{\pi x}{\sin \pi x}-\gamma x-\log\Gamma(1+x)$$

and with $x \to ix$ we obtain

(5.32.1) $$i\sum_{n=1}^{\infty}(-1)^n\frac{\varsigma(2n+1)}{2n+1}x^{2n+1}=\frac{1}{2}\log\frac{\pi x}{\sinh \pi x}-\gamma ix-\log\Gamma(1+ix)$$

Hence we have

$$\sum_{n=1}^{\infty}(-1)^n \frac{\varsigma(2n+1)}{2n+1}x^{2n+1} = -\gamma x - \operatorname{Im}\log\Gamma(1+ix)$$

or equivalently

(5.33) $$\sum_{n=1}^{\infty}(-1)^n \frac{\varsigma(2n+1)}{2n+1}x^{2n+1} = -\gamma x - \frac{1}{2i}[\log\Gamma(1+ix) - \log\Gamma(1-ix)]$$

Letting $x \to ix$ in (5.4) gives us

(5.34) $$\int_0^{\infty}\left[\frac{\sinh(xt)}{t(e^t-1)} - \frac{x}{te^t}\right]dt = \gamma x + \sum_{n=1}^{\infty}\frac{\varsigma(2n+1)}{2n+1}x^{2n+1}$$

and using (5.32) gives us

(5.35) $$\int_0^{\infty}\left[\frac{\sinh(xt)}{t(e^t-1)} - \frac{x}{te^t}\right]dt = \frac{1}{2}\log\frac{\pi x}{\sin\pi x} - \log\Gamma(1+x)$$

$$= \frac{1}{2}[\log\Gamma(1-x) - \log\Gamma(1+x)]$$

Similarly, letting $x \to ix$ in (5.6) gives us

(5.36) $$\int_0^{\infty}\left[\frac{\cosh(xt)}{e^t-1} - \frac{1}{te^t}\right]dt = -\frac{1}{2}[\psi(1+x) + \psi(1-x)]$$

$\square$

We have the Weierstrass canonical form of the gamma function [53, p.1]

(5.36.1) $$\frac{1}{\Gamma(x)} = xe^{\gamma x}\prod_{n=1}^{\infty}\left\{\left(1+\frac{x}{n}\right)e^{-\frac{x}{n}}\right\}$$

where $\gamma$ is Euler's constant defined by $\gamma = \lim_{n\to\infty}[H_n - \log n]$ and $H_n$ are the harmonic numbers $H_n = \sum_{k=1}^{n}\frac{1}{k}$.

Taking logarithms results in

(5.36.2) $$\log\Gamma(1+x) = -\gamma x - \sum_{n=1}^{\infty}\left[\log\left(1+\frac{x}{n}\right) - \frac{x}{n}\right]$$

Hence we have

(5.37)
$$\log \Gamma(1+ix) = -\gamma ix - \sum_{n=1}^{\infty}\left[\log(n+ix) - \log n - \frac{ix}{n}\right]$$

or equivalently

(5.38)
$$\log \Gamma(1+ix) = -\frac{1}{2}\sum_{n=1}^{\infty}\log\left(1+\frac{x^2}{n^2}\right) - \gamma ix - i\sum_{n=1}^{\infty}\left[\tan^{-1}\left(\frac{x}{n}\right) - \frac{x}{n}\right]$$

We note that

$$\gamma = \sum_{n=1}^{\infty}\left[\frac{1}{n} - \log\left(1+\frac{1}{n}\right)\right]$$

and we therefore obtain

(5.39)
$$\log \Gamma(1+ix) = -\frac{1}{2}\sum_{n=1}^{\infty}\log\left(1+\frac{x^2}{n^2}\right) - i\sum_{n=1}^{\infty}\left[\tan^{-1}\left(\frac{x}{n}\right) - x\log\left(1+\frac{1}{n}\right)\right]$$

We recall (5.3)

$$\sum_{n=1}^{\infty}\left[\tan^{-1}\left(\frac{x}{n}\right) - x\log\left(1+\frac{1}{n}\right)\right] = \gamma x + \sum_{n=1}^{\infty}(-1)^n \frac{\varsigma(2n+1)}{2n+1} x^{2n+1}$$

and deduce

(5.40)
$$\log \Gamma(1+ix) = -\frac{1}{2}\sum_{n=1}^{\infty}\log\left(1+\frac{x^2}{n^2}\right) - \gamma ix - i\sum_{n=1}^{\infty}(-1)^n \frac{\varsigma(2n+1)}{2n+1} x^{2n+1}$$

which corresponds with (5.32.1). See also [48, p.23] and Godefroy [35, p.12]

Employing the Taylor series for $\log(1+t)$ we see that

$$\sum_{n=1}^{\infty}\log\left(1+\frac{x^2}{n^2}\right) = -\sum_{n=1}^{\infty}\sum_{k=1}^{\infty}\frac{1}{k}\left(\frac{x^2}{n^2}\right)^k$$

and thus

(5.40.1)
$$\sum_{n=1}^{\infty}\log\left(1+\frac{x^2}{n^2}\right) = -\sum_{n=1}^{\infty}\frac{\varsigma(2n)}{n} x^{2n}$$

Hence, we see that (5.40) may be expressed as the well-known Taylor series [7, p.201]

(5.41)
$$\log \Gamma(1+ix) = -\gamma ix + \sum_{n=1}^{\infty}(-1)^n \frac{\varsigma(n)}{n}(ix)^n$$

Lerch [53] established the following relationship between the gamma function and the Hurwitz zeta function in 1894

$$\varsigma'(0,x) = \log \Gamma(x) - \frac{1}{2}\log(2\pi)$$

In the definition of the Hurwitz zeta function $\varsigma(s,x)$ the argument $x$ is usually regarded as real but, in principle, the is no reason why it could not be a complex number. In this case we would have

$$\log \Gamma(1+ix) = \varsigma'(0,1+ix) + \frac{1}{2}\log(2\pi)$$

Since

$$\varsigma(s,1+x) = \varsigma(s,x) - \frac{1}{x^s}$$

we have the derivative with respect to $s$

$$\varsigma'(s,1+x) = \varsigma'(s,x) + \frac{\log x}{x^s}$$

and thus

$$\varsigma'(0,1+x) = \varsigma'(0,x) + \log x$$

and

$$\varsigma'(0,1+ix) = \varsigma'(0,ix) + \log ix$$

We showed in [26] that

(5.42) $$\varsigma'(0,x) - \varsigma'(0) = -\log x + \sum_{n=1}^{\infty}\left[x[\log(n+1) - \log n] + \log n - \log(n+x)\right]$$

and therefore

(5.43) $$\varsigma'(0,1+x) - \varsigma'(0) = \sum_{n=1}^{\infty}\left[x[\log(n+1) - \log n] + \log n - \log(n+x)\right]$$

Therefore, we have with $x \to ix$

$$\varsigma'(0,1+ix) - \varsigma'(0) = \sum_{n=1}^{\infty}\left[ix[\log(n+1) - \log n] + \log n - \log(n+ix)\right]$$

$$= -\sum_{n=1}^{\infty}\left[\frac{ix}{n} - ix\log\left(1+\frac{1}{n}\right) + \log(n+ix) - \log n - \frac{ix}{n}\right]$$

$$= -\gamma ix - \sum_{n=1}^{\infty}\left[\log(n+ix) - \log n - \frac{ix}{n}\right]$$

resulting in another derivation of (5.37).

We also obtain from the real part of (5.40)

$$\operatorname{Re} \log \Gamma(1+ix) = \frac{1}{2} \log \frac{\pi x}{\sinh \pi x}$$

which may be easily verified as follows. Euler's reflection formula gives us

$$\log \Gamma(1+z) + \log \Gamma(1-z) = \log \pi z - \log \sin \pi z$$

so that

$$\log \Gamma(1+ix) + \log \Gamma(1-ix) = \log i\pi x - \log \sin i\pi x$$

$$= \log \frac{\pi x}{\sinh \pi x}$$

Therefore, we see that

$$\operatorname{Re} \log \Gamma(1+ix) = \frac{1}{2}[\log \Gamma(1+ix) + \log \Gamma(1-ix)]$$

and hence we have

$$\operatorname{Re} \log \Gamma(1+ix) = \frac{1}{2} \log \frac{\pi x}{\sinh \pi x}$$

Hence, we may express (5.40) as

(5.44) $$\log \Gamma(1+ix) = \frac{1}{2} \log \frac{\pi x}{\sinh \pi x} - \gamma ix - i \sum_{n=1}^{\infty} (-1)^n \frac{\varsigma(2n+1)}{2n+1} x^{2n+1}$$

and we see that

$$\log \Gamma(ix) = \frac{1}{2} \log \frac{\pi x}{\sinh \pi x} - \log x - i\left[\gamma x + \frac{\pi}{2} + \sum_{n=1}^{\infty} (-1)^n \frac{\varsigma(2n+1)}{2n+1} x^{2n+1}\right]$$

We have [48, p.24]

$$|\Gamma(ix)| = \frac{1}{|x|} \sqrt{\frac{2\pi x}{e^{\pi x} - e^{-\pi x}}}$$

It is well known that [53, p.2]

(5.44.1) $$\log \Gamma(1+x) = \sum_{n=1}^{\infty} \left[ x \log\left(1+\frac{1}{n}\right) - \log\left(1+\frac{x}{n}\right) \right]$$

and thus

$$\log \Gamma(1+ix) = \sum_{n=1}^{\infty}\left[ix\log\left(1+\frac{1}{n}\right)-\log\left(1+\frac{ix}{n}\right)\right]$$

Hence we have

(5.44.2) $$\log \Gamma(1+ix) + \log \Gamma(1-ix) = -\sum_{n=1}^{\infty}\log\left(1+\frac{x^2}{n^2}\right)$$

$$= \log\frac{\pi x}{\sinh \pi x}$$

We see that

$$\log \Gamma(1+ix) - \log \Gamma(1-ix) = \sum_{n=1}^{\infty}\left[2ix\log\left(1+\frac{1}{n}\right)-\log\frac{n+ix}{n-ix}\right]$$

and we have

$$\log\frac{n+ix}{n-ix} = \log(n^2 - x^2 + 2inx) - \log(n^2 + x^2)$$

$$= \frac{1}{2}\log[(n^2-x^2)^2 + 4n^2x^2] + i\tan^{-1}\left(\frac{2nx}{n^2-x^2}\right) - \log(n^2+x^2)$$

$$= \frac{1}{2}\log[(n^2+x^2)^2] + i\tan^{-1}\left(\frac{2nx}{n^2-x^2}\right) - \log(n^2+x^2)$$

$$= i\tan^{-1}\left(\frac{2nx}{n^2-x^2}\right)$$

This was a rather convoluted derivation and an easier approach is to consider the expression $\log(n+ix) - \log(n-ix)$ and obtain the equivalent result

$$\log\frac{n+ix}{n-ix} = 2i\tan^{-1}\left(\frac{x}{n}\right)$$

Hence we deduce that

(5.44.3) $$\log \Gamma(1+ix) - \log \Gamma(1-ix) = 2i\sum_{n=1}^{\infty}\left[x\log\left(1+\frac{1}{n}\right)-\tan^{-1}\left(\frac{x}{n}\right)\right]$$

In this regard, see (5.4).

□

Integrating (5.36.2) results in

$$\int_0^u \log \Gamma(x)dx = u - u\log u - \frac{1}{2}\gamma u^2 - \sum_{n=1}^{\infty}\left[(u+n)\log(u+n) - u - (u+n)\log n - \frac{u^2}{2n}\right]$$

and in particular we have

$$\int_0^1 \log \Gamma(x)\,dx = 1 - \frac{1}{2}\gamma - \sum_{n=1}^{\infty}\left[(1+n)\log(1+n) - 1 - (1+n)\log n - \frac{1}{2n}\right]$$

We deduce that

(5.44.4) $$\sum_{n=1}^{\infty}\left[(1+n)\log\left(1+\frac{1}{n}\right) - 1 - \frac{1}{2n}\right] = 1 - \frac{1}{2}[\gamma + \log(2\pi)]$$

We note that

$$\sum_{n=1}^{\infty}\left[(1+n)\log\left(1+\frac{1}{n}\right) - 1 - \frac{1}{2n}\right] = \sum_{n=1}^{\infty}\left[\left(\frac{1}{2}+n\right)\log\left(1+\frac{1}{n}\right) - 1 + \frac{1}{2}\left\{\log\left(1+\frac{1}{n}\right) - \frac{1}{n}\right\}\right]$$

and determine that

(5.44.5) $$\sum_{n=1}^{\infty}\left[\left(\frac{1}{2}+n\right)\log\left(1+\frac{1}{n}\right) - 1\right] = 1 - \frac{1}{2}\log(2\pi)$$

The latter summation features in Hardy's book [38, p.335]. It also appears later in (5.58.1).

**The ubiquitous integral** $\int_0^u \frac{\psi(1+x)+\gamma}{x}\,dx$

Cohen [15, p.142] has reported as an exercise that (see also [22])

(5.45) $$\int_0^1 \frac{\psi(1+x)+\gamma}{x}\,dx = \sum_{n=1}^{\infty}\frac{\log(n+1)}{n(n+1)}$$

Cohen [15, p.142] has also stated that

(5.45.1) $$\sum_{n=1}^{\infty}\frac{\log(n+1)}{n(n+1)} = \int_0^1 \frac{(1-x)\log(1-x)}{x\log x}\,dx$$

(5.45.2) $$= \sum_{n=1}^{\infty}(-1)^{n+1}\frac{\varsigma(n+1)}{n}$$

(5.45.3) $$= -\sum_{n=2}^{\infty}\varsigma'(n)$$

(5.45.4) $$= \sum_{n=1}^{\infty}\frac{1}{n}\log\left(1+\frac{1}{n}\right)$$

We recall (5.31)

(5.46) $$\psi(1+x) = \frac{1}{2x} - \frac{\pi}{2}\cot(\pi x) - \gamma - \frac{x^2}{1-x^2} - \sum_{n=1}^{\infty}[\varsigma(2n+1)-1]x^{2n}$$

and thus we see that

$$\psi(1-x) = -\frac{1}{2x} + \frac{\pi}{2}\cot(\pi x) - \gamma - \frac{x^2}{1-x^2} - \sum_{n=1}^{\infty}[\varsigma(2n+1)-1]x^{2n}$$

We therefore obtain

(5.46.1) $$\psi(1+x)+\gamma+\psi(1-x)+\gamma+\frac{2x^2}{1-x^2} = -2\sum_{n=1}^{\infty}[\varsigma(2n+1)-1]x^{2n}$$

or equivalently

$$\sum_{n=1}^{\infty}\varsigma(2n+1)x^{2n} = -\gamma - \frac{1}{2}[\psi(1+x)+\psi(1-x)]$$

which appears in [53, p.160].

Curiously, *WolframAlpha* indicates that $\sum_{n=1}^{\infty}[\varsigma(2n+1)-1]x^{2n}$ diverges at $x=1$ but, as shown below, this is not correct. We see that

$$\psi(1-x)+\frac{2x^2}{1-x^2} = \frac{(1-x^2)\left[\psi(2-x)-\frac{1}{1-x}\right]+2x^2}{1-x^2}$$

$$= \frac{(1-x^2)\psi(2-x)-(1+x)+2x^2}{1-x^2}$$

and, applying L'Hôpital's rule, we find

$$\lim_{x\to 1}\left[\psi(1-x)+\frac{2x^2}{1-x^2}\right] = \lim_{x\to 1}\frac{-(1-x^2)\psi'(2-x)-2x\psi(2-x)-1+4x}{-2x}$$

Thus we have the finite limit

$$\lim_{x\to 1}\left[\psi(1-x)+\frac{2x^2}{1-x^2}\right] = -\left(\gamma+\tfrac{3}{2}\right)$$

Hence we deduce from (5.46.1) that

(5.46.2) $\sum_{n=1}^{\infty}[\varsigma(2n+1)-1] = \tfrac{1}{4}$

We also note from [53, p.142] that $\sum_{n=2}^{\infty}[\varsigma(n)-1]=1$.

An alternative proof of (5.46.2) is shown below. Referring to

(5.46.3) $$\frac{x^2}{1+x^2}+\sum_{n=1}^{\infty}(-1)^n[1-\varsigma(2n+1)]x^{2n} = \sum_{j=1}^{\infty}\frac{x^2}{j(j^2+x^2)}$$

we easily see that

$$\sum_{n=1}^{\infty}(-1)^n[1-\varsigma(2n+1)]u^n = \sum_{j=1}^{\infty}\frac{u}{j(j^2+u)}-\frac{u}{1+u}$$

and

$$\sum_{n=1}^{\infty}[\varsigma(2n+1)-1]u^n = \sum_{j=1}^{\infty}\frac{u}{j(j^2-u)}-\frac{u}{1-u}$$

or alternatively

(5.46.4) $$\sum_{n=1}^{\infty}[\varsigma(2n+1)-1]u^n = \sum_{j=2}^{\infty}\frac{u}{j(j^2-u)}$$

Thus

$$\sum_{n=1}^{\infty}[\varsigma(2n+1)-1] = \sum_{j=2}^{\infty}\frac{1}{j(j^2-1)}$$

and using

$$\sum_{j=2}^{\infty}\frac{1}{j(j^2-1)}=\frac{1}{4}$$

this completes the proof.

$\square$

Division of (5.46.1) by $x$ followed by integration results in

(5.47) $$\int_0^u \frac{\psi(1+x)+\gamma}{x}dx + \int_0^u\left[\frac{\psi(1-x)+\gamma}{x}+\frac{2x}{1-x^2}\right]dx = -\sum_{n=1}^{\infty}\frac{\varsigma(2n+1)-1}{n}u^{2n}$$

It may be noted that *WolframAlpha* also incorrectly reports $\sum_{n=1}^{\infty}\frac{\varsigma(2n+1)-1}{n}$ as being divergent.

It is known from (5.19) that

$$\frac{G'(1+x)}{G(1+x)} = \frac{1}{2}\log(2\pi) - \frac{1}{2} - x + x\psi(1+x)$$

and with $x \to -x$ we have

$$\frac{G'(1-x)}{G(1-x)} = \frac{1}{2}\log(2\pi) - \frac{1}{2} + x - x\psi(1-x)$$

We then have

(5.47.1) $$\frac{G'(1+x)}{G(1+x)} - \frac{G'(1-x)}{G(1-x)} = -2x + x[\psi(1+x) + \psi(1-x)]$$

We multiply (5.46.1) by $x$ and integrate to obtain for $|x| < 2$

(5.48) $$\sum_{n=1}^{\infty}[\varsigma(2n+1)-1]\frac{x^{2n+2}}{n+1} = -\gamma x^2 - \log\frac{G(1+x)G(1-x)}{1-x^2}$$

where we have used (5.47.1). This appears in [53, p.212, Eq. (472)] in the equivalent form

$$\sum_{n=1}^{\infty}[\varsigma(2n+1)-1]\frac{x^{2n+2}}{n+1} = -\gamma x^2 - \log\frac{G(2+x)G(2-x)}{\Gamma(2+x)\Gamma(2-x)}$$

$\square$

We recall (2.8)

$$\psi(1+x) + \psi(1-x) + 2\gamma = -2x^2 \sum_{n=1}^{\infty}\frac{1}{n}\frac{1}{n^2 - x^2}$$

which we write as

$$\frac{\psi(1+x)+\gamma}{x} + \frac{\psi(1-x)+\gamma}{x} + \frac{2x}{1-x^2} = -2x\sum_{n=2}^{\infty}\frac{1}{n}\frac{1}{n^2 - x^2}$$

Integration gives us

(5.49) $$\int_0^1 \frac{\psi(1+x)+\gamma}{x}dx + \int_0^1\left[\frac{\psi(1-x)+\gamma}{x} + \frac{2x}{1-x^2}\right]dx = \sum_{n=2}^{\infty}\frac{1}{n}\log\left(1 - \frac{1}{n^2}\right)$$

We may also derive this by dividing (5.46.4) by $x$ and then integrating to obtain

(5.49.1) $$\sum_{n=1}^{\infty}\frac{[\varsigma(2n+1)-1]}{n}x^n = -\sum_{n=2}^{\infty}\frac{1}{n}\log\left(1 - \frac{x}{n^2}\right)$$

*WolframAlpha* provides us with the approximation

$$\sum_{n=2}^{\infty} \frac{1}{n} \log\left(1 - \frac{1}{n^2}\right) \simeq -0.223931...$$

Since

$$\sum_{n=2}^{\infty} \frac{1}{n} \log\left(1 - \frac{1}{n^2}\right) - \sum_{n=2}^{\infty} \frac{1}{n} \log\left(1 + \frac{1}{n}\right) = \sum_{n=2}^{\infty} \frac{1}{n} \log\left(1 - \frac{1}{n}\right)$$

referring to (5.49) and using (5.45.4) we see that

(5.50) $$\int_0^1 \left[\frac{\psi(1-x) + \gamma}{x} + \frac{2x}{1-x^2}\right] dx = -\frac{1}{2} \log 2 + \sum_{n=2}^{\infty} \frac{1}{n} \log\left(1 - \frac{1}{n}\right)$$

We have

(5.51)

$$\int_0^1 \frac{\psi(1+x) + \gamma}{x} dx - \int_0^1 \left[\frac{\psi(1-x) + \gamma}{x} + \frac{2x}{1-x^2}\right] dx = \int_0^1 \left[\frac{\psi(x) - \psi(1-x)}{x} + \frac{1}{x^2} - \frac{2x}{1-x^2}\right] dx$$

$$= \int_0^1 \left[\frac{1}{x}\left(\frac{1}{x} - \pi \cot(\pi x)\right) - \frac{2x}{1-x^2}\right] dx$$

We showed in [22, Eq. (3.14)] that

(5.52) $$\int_0^u \frac{\psi(1+x) + \gamma}{x} dx = \sum_{n=1}^{\infty} \frac{1}{n} \log \frac{n+u}{n}$$

Hence, we obtain for $u < 1$

(5.53) $$\int_0^u \frac{1}{x}\left[\frac{1}{x} - \pi \cot(\pi x)\right] dx = 2\sum_{n=1}^{\infty} \frac{\varsigma(2n+1)}{n} u^{2n} + 2\sum_{n=1}^{\infty} \frac{1}{n} \log \frac{n+u}{n}$$

Integration by parts gives us

$$\int_0^u \frac{1}{x}\left[\frac{1}{x} - \pi \cot(\pi x)\right] dx = -\log \frac{\sin(\pi x)}{\pi x}\bigg|_0^u + \int_0^u \frac{1}{x^2} \log \frac{\sin(\pi x)}{\pi x} dx$$

$$= -\log \frac{\sin(\pi u)}{\pi u} + \int_0^u \frac{1}{x^2} \log \frac{\sin(\pi x)}{\pi x} dx$$

Using Euler's reflection formula

$$\Gamma(x)\Gamma(1-x) = \frac{\pi}{\sin(\pi x)}$$

we see that

$$-\log\frac{\sin(\pi x)}{\pi x} = \log\Gamma(1+x) + \log\Gamma(1-x)$$

and thus

$$\int_0^u \frac{1}{x^2}\log\frac{\sin(\pi x)}{\pi x}dx = -\int_0^u \frac{\log\Gamma(1+x)+\log\Gamma(1-x)}{x^2}dx$$

We have

$$\log\Gamma(1+x) = -\gamma x + \sum_{n=2}^{\infty}\frac{(-1)^n \varsigma(n)}{n}x^n$$

and

$$\log\Gamma(1-x) = \gamma x + \sum_{n=2}^{\infty}\frac{\varsigma(n)}{n}x^n$$

so that

$$\log\Gamma(1+x) + \log\Gamma(1-x) = \sum_{n=1}^{\infty}\frac{\varsigma(2n)}{n}x^{2n}$$

Hence we have the integral

$$\int_0^u \frac{\log\Gamma(1+x)+\log\Gamma(1-x)}{x^2}dx = \sum_{n=1}^{\infty}\frac{\varsigma(2n)}{n(2n-1)}u^{2n-1}$$

We then see that for $u < 1$

(5.54) $$\int_0^u \frac{1}{x}\left[\frac{1}{x} - \pi\cot(\pi x)\right]dx = \sum_{n=1}^{\infty}\frac{\varsigma(2n)}{n}u^{2n} + \sum_{n=1}^{\infty}\frac{\varsigma(2n)}{n(2n-1)}u^{2n-1}$$

We recall (5.31)

$$\psi(1+x) = \frac{1}{2x} - \frac{\pi}{2}\cot(\pi x) - \gamma - \frac{x^2}{1-x^2} - \sum_{n=1}^{\infty}[\varsigma(2n+1)-1]x^{2n}$$

whereupon integration gives us

$$2\int_0^u \frac{\psi(1+x)+\gamma}{x}dx = \int_0^u \frac{1}{x}\left[\frac{1}{x} - \pi\cot(\pi x)\right]dx - \sum_{n=1}^{\infty}\frac{\varsigma(2n+1)}{n}u^{2n}$$

It is obvious by comparison with the harmonic series that $\sum_{n=1}^{\infty}\frac{\varsigma(2n+1)}{n}u^{2n}$ is divergent at $u = 1$. *WolframAlpha* confirms that $\int_0^1 \frac{1}{x}\left[\frac{1}{x} - \pi\cot(\pi x)\right]dx$ does not converge. We accordingly consider

$$\int_0^u \frac{\psi(1+x)+\gamma}{x}dx = \int_0^u \left[\frac{1}{2x}\left(\frac{1}{x}-\pi\cot(\pi x)\right)-\frac{x}{1-x^2}\right]dx - \frac{1}{2}\sum_{n=1}^{\infty}\frac{\varsigma(2n+1)-1}{n}u^{2n}$$

and *WolframAlpha* confirms that $\int_0^1 \left[\frac{1}{2x}\left(\frac{1}{x}-\pi\cot(\pi x)\right)-\frac{x}{1-x^2}\right]dx \simeq 1.36971...$ does converge.

We then obtain

(5.55) $$2\int_0^u \frac{\psi(1+x)+\gamma}{x}dx = \sum_{n=1}^{\infty}\frac{\varsigma(2n)}{n}u^{2n} + \sum_{n=1}^{\infty}\frac{\varsigma(2n)}{n(2n-1)}u^{2n-1} - \sum_{n=1}^{\infty}\frac{\varsigma(2n+1)}{n}u^{2n}$$

Substituting (5.52) gives us

$$2\sum_{n=1}^{\infty}\frac{1}{n}\log\frac{n+u}{n} = \sum_{n=1}^{\infty}\frac{\varsigma(2n)}{n}u^{2n} + \sum_{n=1}^{\infty}\frac{\varsigma(2n)}{n(2n-1)}u^{2n-1} - \sum_{n=1}^{\infty}\frac{\varsigma(2n+1)}{n}u^{2n}$$

We have [53, p.161] for $|u|<1$

$$\sum_{n=1}^{\infty}\frac{\varsigma(2n)}{n}u^{2n} = \log\frac{\pi u}{\sin \pi u}$$

and therefore

$$\sum_{n=1}^{\infty}\frac{\varsigma(2n)}{n(2n-1)}u^{2n-1} = \int_0^u \frac{1}{x^2}\log\frac{\pi x}{\sin \pi x}dx$$

With $u=1$ in (5.55) we obtain

$$2\int_0^1 \frac{\psi(1+x)+\gamma}{x}dx = \sum_{n=1}^{\infty}\left[\frac{\varsigma(2n)}{n} + \frac{\varsigma(2n)}{n(2n-1)} - \frac{\varsigma(2n+1)}{n}\right]$$

$$= \sum_{n=1}^{\infty}\left[\frac{2\varsigma(2n)}{2n-1} - \frac{\varsigma(2n+1)}{n}\right]$$

Hence we obtain

$$\int_0^1 \frac{\psi(1+x)+\gamma}{x}dx = \sum_{n=1}^{\infty}\left[\frac{\varsigma(2n)}{2n-1} - \frac{\varsigma(2n+1)}{2n}\right]$$

This may be written as

$$\int_0^1 \frac{\psi(1+x)+\gamma}{x}dx = \sum_{n=1}^{\infty}(-1)^{n+1}\frac{\varsigma(n+1)}{n}$$

as noted by Cohen [15, p.142].

We recall (5.26.2)

$$\psi(1+x)+\gamma = -\sum_{n=1}^{\infty}\frac{x^2}{n(n^2-x^2)}+\sum_{n=1}^{\infty}\frac{x}{n^2-x^2}$$

$$= -\sum_{n=2}^{\infty}\frac{x^2}{n(n^2-x^2)}+\sum_{n=2}^{\infty}\frac{x}{n^2-x^2}-\frac{x^2}{1-x^2}+\frac{x}{1-x^2}$$

$$= -\sum_{n=2}^{\infty}\frac{x^2}{n(n^2-x^2)}+\sum_{n=2}^{\infty}\frac{x}{n^2-x^2}+1-\frac{1}{1+x}$$

$$\int_0^1\frac{\psi(1+x)+\gamma}{x}dx = \sum_{n=2}^{\infty}\int_0^1\left[\frac{1}{n^2-x^2}-\frac{x}{n(n^2-x^2)}\right]dx+1-\log 2$$

$$= \sum_{n=2}^{\infty}\frac{1}{n}\tanh^{-1}\left(\frac{1}{n}\right)+\frac{1}{2}\sum_{n=2}^{\infty}\frac{1}{n}\log\left(1-\frac{1}{n^2}\right)+1-\log 2$$

Noting that $\tanh^{-1}x = \frac{1}{2}\log\frac{1+x}{1-x}$ we end up with

$$= \frac{1}{2}\sum_{n=2}^{\infty}\frac{1}{n}\log\left(\frac{n+1}{n-1}\right)+\frac{1}{2}\sum_{n=2}^{\infty}\frac{1}{n}\log\left(1-\frac{1}{n^2}\right)+1-\log 2$$

$$= \sum_{n=2}^{\infty}\frac{1}{n}\log\left(1+\frac{1}{n}\right)+1-\log 2$$

However, there must be an error here because this does not concur with (5.45.4). Corrections are most welcome.

□

We recall (5.27)
$$\sum_{j=1}^{\infty}\frac{1}{j(j^2+u)} = \sum_{n=1}^{\infty}(-1)^{n+1}\varsigma(2n+1)u^{n-1}$$

Integrating this gives us

$$\sum_{j=1}^{\infty}\frac{1}{j}\log\left(1+\frac{u}{j^2}\right) = \sum_{n=1}^{\infty}(-1)^{n+1}\frac{\varsigma(2n+1)}{n}u^n$$

We see that

$$\sum_{j=1}^{\infty} \frac{1}{j} \log\left(1+\frac{u^2}{j^2}\right) = \sum_{n=1}^{\infty} (-1)^{n+1} \frac{\varsigma(2n+1)}{n} u^{2n}$$

$$= \frac{1}{2} \int_0^u \frac{\psi(1+ix)+\psi(1-ix)+2\gamma}{x} dx$$

where, in the last part, we have used [53, p.160].

**Proposition 5.2**

(5.56) $$\sum_{j=2}^{\infty} \left[ j \log\left(1-\frac{1}{j}\right) + 1 + \frac{1}{2j} \right] = \frac{1}{2}[\gamma + \log(2\pi) - 3]$$

**Proof**

We showed in (5.46.4) that

$$\sum_{n=1}^{\infty} [\varsigma(2n+1)-1] u^n = \sum_{j=2}^{\infty} \frac{u}{j(j^2-u)}$$

It is easily seen that

$$\sum_{n=1}^{\infty} \frac{n^2}{n+1} [\varsigma(2n+1)-1] x^{n+1} = \int_0^x u \frac{d}{du}\left\{ u \frac{df}{du} \right\} du$$

where $f(u) := \sum_{n=1}^{\infty} [\varsigma(2n+1)-1] u^n$. Integration by parts gives us

$$\int_0^x u \frac{d}{du}[uf'(u)] du = u^2 f'(u)\Big|_0^x - \int_0^x u f'(u) du$$

$$= x^2 f'(x) - u f(u)\Big|_0^x + \int_0^x f(u) du$$

$$= x^2 f'(x) - x f(x) + \int_0^x f(u) du$$

We also have

$$f(u) = \sum_{j=2}^{\infty} \frac{u}{j(j^2-u)}$$

and

$$f'(u) = \sum_{j=2}^{\infty} \frac{j}{(j^2-u)^2}$$

Thus

$$\sum_{n=1}^{\infty}\frac{n^2}{n+1}[\varsigma(2n+1)-1]\,x^{n+1} = x^2\sum_{j=2}^{\infty}\frac{j}{(j^2-x)^2} - x^2\sum_{j=2}^{\infty}\frac{1}{j(j^2-x)} + \int_0^x\sum_{j=2}^{\infty}\frac{u}{j(j^2-u)}du$$

We have

$$\int_0^x\sum_{j=2}^{\infty}\frac{u}{j(j^2-u)}du = \int_0^x\sum_{j=2}^{\infty}\frac{1}{j}\left[\frac{j^2}{j^2-u}-1\right]du$$

$$= -\sum_{j=2}^{\infty}\frac{1}{j}\left[j^2\log\left(1-\frac{x}{j^2}\right)+x\right]$$

which results in

(5.57)
$$\sum_{n=1}^{\infty}\frac{n^2}{n+1}[\varsigma(2n+1)-1]\,x^{n+1} = x^2\sum_{j=2}^{\infty}\frac{j}{(j^2-x)^2} - x^2\sum_{j=2}^{\infty}\frac{1}{j(j^2-x)} - \sum_{j=2}^{\infty}\left[j\log\left(1-\frac{x}{j^2}\right)+\frac{x}{j}\right]$$

Letting $x=1$ gives us

$$\sum_{n=1}^{\infty}\frac{n^2}{n+1}[\varsigma(2n+1)-1] = \sum_{j=2}^{\infty}\frac{j}{(j^2-1)^2} - \sum_{j=2}^{\infty}\frac{1}{j(j^2-1)} - \sum_{j=2}^{\infty}\left[j\log\left(1-\frac{1}{j^2}\right)+\frac{1}{j}\right]$$

Substituting

$$\sum_{j=2}^{\infty}\frac{j}{(j^2-1)^2} = \frac{5}{16}$$

$$\sum_{j=2}^{\infty}\frac{1}{j(j^2-1)} = \frac{1}{4}$$

we obtain

$$\sum_{n=1}^{\infty}\frac{n^2}{n+1}[\varsigma(2n+1)-1] = \frac{1}{16} - \sum_{j=2}^{\infty}\left[j\log\left(1-\frac{1}{j^2}\right)+\frac{1}{j}\right]$$

We note [2] and [53, p.155] have reported that

$$\sum_{n=2}^{\infty}\frac{n^2}{n+1}[\varsigma(2n+1)-1] = \frac{9}{16} - \gamma + \log 2 - \frac{1}{2}\varsigma(3)$$

or equivalently

$$\sum_{n=1}^{\infty} \frac{n^2}{n+1}[\varsigma(2n+1)-1] = \frac{1}{16} - \gamma + \log 2$$

It should be noted that *WolframAlpha* incorrectly reports that the series $\sum_{n=1}^{\infty} \frac{n^2}{n+1}[\varsigma(2n+1)-1]$ is divergent but inexplicably states, at the same time, that its value is $\frac{1}{16} - \gamma + \log 2$.

We therefore obtain

$$\sum_{j=2}^{\infty} \left[ j \log\left(1 - \frac{1}{j^2}\right) + \frac{1}{j} \right] = \gamma - \log 2$$

*WolframAlpha* reports that:

Partial sum formula

$$\sum_{n=2}^{k} \left( n \log\left(1 - \frac{1}{n^2}\right) + \frac{1}{n} \right) =$$
$$\frac{1}{2} \left( 2\zeta^{(1,0)}(-1, k) - 4\zeta^{(1,0)}(-1, k+1) + 2\zeta^{(1,0)}(-1, k+2) + \right.$$
$$\left. 2 H_k + 2 \log(\Gamma(k)) - 2 \log(\Gamma(k+2)) - 2 - \log(4) \right)$$

Sum

$$\sum_{n=2}^{10000} \left( n \log\left(1 - \frac{1}{n^2}\right) + \frac{1}{n} \right) \approx -0.115931513158$$

We see that

$$\sum_{j=2}^{\infty} \left[ j \log\left(1 - \frac{1}{j^2}\right) + \frac{1}{j} \right] = \sum_{j=2}^{\infty} \left[ j \log\left(1 - \frac{1}{j^2}\right) + \frac{1}{j} - \log\left(1 + \frac{1}{j}\right) + \log\left(1 + \frac{1}{j}\right) \right]$$

and using the series for Euler's constant

$$\gamma = \sum_{j=1}^{\infty} \left[ \frac{1}{j} - \log\left(1 + \frac{1}{j}\right) \right]$$

we obtain

$$\sum_{j=2}^{\infty} \left[ j \log\left(1 - \frac{1}{j^2}\right) + \frac{1}{j} \right] = \gamma + \log 2 - 1 + \sum_{j=2}^{\infty} \left[ j \log\left(1 - \frac{1}{j^2}\right) + \log\left(1 + \frac{1}{j}\right) \right]$$

Therefore

$$\sum_{j=2}^{\infty}\left[j\log\left(1-\frac{1}{j^2}\right)+\log\left(1+\frac{1}{j}\right)\right]=1-2\log 2$$

It is also apparent that

(5.58)
$$\sum_{j=2}^{\infty}\left[j\log\left(1-\frac{1}{j^2}\right)+\log\left(1+\frac{1}{j}\right)\right]=\sum_{j=2}^{\infty}\left[j\log\left(1-\frac{1}{j}\right)+j\log\left(1+\frac{1}{j}\right)+\log\left(1+\frac{1}{j}\right)\right]$$

Candelpergher [13, p.38] has shown that

(5.58.1) $$\sum_{j=1}^{\infty}\left[j\log\left(1+\frac{1}{j}\right)-1+\frac{1}{2j}\right]=\frac{1}{2}[\gamma-\log(2\pi)]+1$$

or equivalently

$$\sum_{j=2}^{\infty}\left[j\log\left(1+\frac{1}{j}\right)-1+\frac{1}{2j}\right]=\frac{1}{2}[\gamma-\log\pi-3\log 2+3]$$

Candelpergher's formula may be easily deduced from (5.44.5)

$$\sum_{n=1}^{\infty}\left[\left(\frac{1}{2}+n\right)\log\left(1+\frac{1}{n}\right)-1\right]=1-\frac{1}{2}\log(2\pi)$$

Employing (5.58) we then see that

$$\sum_{j=2}^{\infty}\left[j\log\left(1-\frac{1}{j}\right)+1-\frac{1}{2j}+\log\left(1+\frac{1}{j}\right)\right]=1-2\log 2-\frac{1}{2}[\gamma-\log\pi-3\log 2+3]$$

and using the series for Euler's constant this may be expressed as

$$\sum_{j=2}^{\infty}\left[j\log\left(1-\frac{1}{j}\right)+1+\frac{1}{2j}\right]=\frac{1}{2}[\gamma+\log(2\pi)-3]$$

as confirmed by *WolframAlpha*.

We recall (5.16.1)

$$\sum_{n=1}^{\infty}\left[n\log\left(1-\frac{x^2}{n^2}\right)+x^2\log\left(1+\frac{1}{n}\right)\right]=\sum_{n=1}^{\infty}\left[\frac{1}{n}x^2+n\log\left(1-\frac{x^2}{n^2}\right)\right]-\gamma x^2$$

$$=x^2+\log[G(1+x)G(1-x)]$$

and, starting the summation at $n=2$, results in

(5.59)
$$\sum_{n=2}^{\infty}\left[n\log\left(1-\frac{x^2}{n^2}\right)+x^2\log\left(1+\frac{1}{n}\right)\right]=x^2+\log[G(1+x)G(1-x)]-\log(1-x^2)-x^2\log 2$$

We have

$$\frac{G(1-x)}{1-x^2}=\frac{G(2-x)}{(1+x)(1-x)\Gamma(1-x)}$$

$$=\frac{G(2-x)}{(1+x)\Gamma(2-x)}$$

and hence

$$\lim_{x\to 1}\frac{G(1-x)}{1-x^2}=\frac{1}{2}$$

Therefore, we see that

$$\sum_{n=2}^{\infty}\left[n\log\left(1-\frac{1}{n^2}\right)+\log\left(1+\frac{1}{n}\right)\right]=1-2\log 2$$

which we saw above.

## 6. Miscellaneous series involving the logarithmic function

As mentioned by Witula et al. [57], Salaev [50] showed that for $x\in[0,1)$

(6.1) $-\sum_{n=2}^{\infty}\log\left(1-\frac{1}{n^2}\right)\cos 2n\pi x=\log 2-\frac{\pi}{2}\sin 2\pi x+(1-\cos 2\pi x)[\log\pi+\gamma+\psi(1-x)]$

We have [53, p.14]

$$\psi(1-x)=-\gamma+\sum_{m=1}^{\infty}\left(\frac{1}{m}-\frac{1}{m-x}\right)$$

and we therefore see that $\lim_{x\to 1}(1-x)\psi(1-x)=-1$.

Noting that $\lim_{x\to 1}[(1-\cos 2\pi x)\psi(1-x)]=\lim_{x\to 1}\frac{1-\cos 2\pi x}{1-x}\lim_{x\to 1}[(1-x)\psi(1-x)]$ we see from L'Hôpital's rule that $\lim_{x\to 1}(1-\cos 2\pi x)\psi(1-x)=0$. We therefore note that (6.1) also holds in the limit as $x\to 1$.

With $x\to 1-x$ in (6.1) we see that

$$(6.2) \quad -\sum_{n=2}^{\infty} \log\left(1 - \frac{1}{n^2}\right) \cos 2n\pi x = \log 2 + \frac{\pi}{2} \sin 2\pi x + (1 - \cos 2\pi x)[\log \pi + \gamma + \psi(x)]$$

Part of this section was drafted several years ago and I belatedly learned from Boyack's recent paper [9] that (6.2) had in fact been derived much earlier by Lerch [44] in 1897.

Letting $x = 0$ in (6.1) gives us

$$(6.3) \quad \sum_{n=2}^{\infty} \log\left(1 - \frac{1}{n^2}\right) = -\log 2$$

and, as noted below, this is a known result (and was recently reported by Boyack [9] in 2021).

Using the well-known identity originally postulated by Euler

$$\sin \pi x = \pi x \prod_{n=1}^{\infty} \left(1 - \frac{x^2}{n^2}\right)$$

we have

$$\log \frac{\sin \pi x}{\pi x(1-x)} = \log(1+x) + \sum_{n=2}^{\infty} \log\left(1 - \frac{x^2}{n^2}\right)$$

L'Hôpital's rule gives us

$$\lim_{x \to 1} \frac{\sin \pi x}{\pi x(1-x)} = \lim_{x \to 1} \frac{\cos \pi x}{1 - 2x} = 1$$

and from this we easily deduce (6.3). A more direct derivation is shown below:

$$-\sum_{n=2}^{N} \log\left(\frac{n^2 - 1}{n^2}\right) = \sum_{n=2}^{N} \log\left(\frac{n}{n-1}\right) + \sum_{n=2}^{N} \log\left(\frac{n}{n+1}\right)$$

$$= \log N + \log\left(\frac{2}{N+1}\right)$$

$$= \log\left(\frac{2}{1 + 1/N}\right)$$

which approaches $\log 2$ as $N \to \infty$.

□

With $x = 1/2$ in (6.2) we obtain

(6.4) $$\sum_{n=2}^{\infty}(-1)^{n+1}\log\left(1-\frac{1}{n^2}\right)=2\log\pi-3\log 2$$

where we have used [53, p.20] $\psi\left(\frac{1}{2}\right)=-\gamma-2\log 2$.

It is known that [51]

(6.5) $$\log\frac{\pi}{2}=\sum_{n=1}^{\infty}(-1)^{n+1}\log\left(1+\frac{1}{n}\right)$$

being the logarithm of Wallis's product formula which was published in Arithmetica Infinitorum in 1659

$$\frac{\pi}{2}=\prod_{n=1}^{\infty}\frac{(2n)^2}{(2n-1)(2n+1)}$$

Then noting that

$$\log\frac{\pi}{2}=\log 2+\sum_{n=2}^{\infty}(-1)^{n+1}\log\left(1+\frac{1}{n}\right)$$

and

$$\sum_{n=2}^{\infty}(-1)^{n+1}\log\left(1-\frac{1}{n^2}\right)=\sum_{n=2}^{\infty}(-1)^{n+1}\log\left(1-\frac{1}{n}\right)+\sum_{n=2}^{\infty}(-1)^{n+1}\log\left(1+\frac{1}{n}\right)$$

we obtain (what appears to be a new result)

(6.6) $$\sum_{n=2}^{\infty}(-1)^{n+1}\log\left(1-\frac{1}{n}\right)=\log\frac{\pi}{2}$$

which concurs with an evaluation by *WolframAlpha*.

□

Letting $x\to x+\frac{1}{2}$ in (6.1) gives us for $-\frac{1}{2}\le x<\frac{1}{2}$

(6.7)
$$-\sum_{n=2}^{\infty}(-1)^n\log\left(1-\frac{1}{n^2}\right)\cos 2n\pi x=\log 2+\frac{\pi}{2}\sin 2\pi x+(1+\cos 2\pi x)[\log\pi+\gamma+\psi(\tfrac{1}{2}-x)]$$

This also applies in the limit as $x\to\frac{1}{2}$. With $x=0$ we again recover (6.4).

□

We recall

$$\log\left(1+\frac{1}{n}\right) = \int_0^1 \frac{y^n(y-1)}{y \log y} dy$$

and form the summation

$$\sum_{n=1}^{\infty} (-1)^{n+1} \log\left(1+\frac{1}{n}\right) = \sum_{n=1}^{\infty} (-1)^{n+1} \int_0^1 \frac{y^n(y-1)}{y \log y} dy$$

$$= \int_0^1 \frac{y-1}{(y+1)\log y} dy$$

Hence we see that [76]

(6.7) $$\log \frac{\pi}{2} = \int_0^1 \frac{y-1}{(1+y)\log y} dy$$

Three other proofs of this result are given in Sondow's paper [51]. I subsequently discovered that a derivation of this integral was also recorded in 1864 in Bertrand's treatise [6, Book II, p.149], albeit with a trivial sign error being made in the last step of the proof.

□

We easily determine from (6.2) that

(6.7.1) $$\int_0^1 (1-\cos 2\pi x)\psi(1-x) dx = -[\log 2\pi + \gamma]$$

We multiply (6.2) by $x(1-x)$ and integrate over the interval $[0,1]$ to obtain

(6.7.2) $$\frac{1}{2\pi^2} \sum_{n=2}^{\infty} \frac{1}{n^2} \log\left(1-\frac{1}{n^2}\right) = \frac{1}{6}\log 2 + \log \pi + \gamma + \int_0^1 x(1-x)(1-\cos 2\pi x)\psi(x) dx$$

where we have employed some elementary integrals including

$$\int_0^1 x(1-x) \cos 2n\pi x\, dx = -\frac{1}{2\pi^2 n^2}$$

□

Integration by parts results in

$$\int_0^1 (1-\cos 2\pi x)\psi(x) dx = (1-\cos 2\pi x)\log \Gamma(x)\Big|_0^1 - 2\pi \int_0^1 \log \Gamma(x)\sin 2\pi x\, dx$$

$$= -2\pi \int_0^1 \log \Gamma(x) \sin 2\pi x\, dx$$

and we then see that integrating (6.2) results in

(6.8) $$\int_0^1 \log \Gamma(x) \sin 2\pi x \, dx = \frac{\gamma + \log(2\pi)}{2\pi}$$

which is a specific case of the well-known Fourier coefficients [20]

(6.9) $$\int_0^1 \log \Gamma(x) \sin 2k\pi x \, dx = \frac{\gamma + \log(2\pi k)}{2\pi k}$$

$\square$

We use (6.2) to derive the following integral.

**Proposition 6.1**

(6.10) $$\int_0^u \psi(1+x) \cos^2 \pi x \, dx$$

$$= \log \Gamma(1+u) + \frac{1}{4\pi} \sum_{n=2}^{\infty} \frac{1}{n} \log\left(1 - \frac{1}{n^2}\right) \sin 2n\pi u + \frac{1}{2} u \log 2 + \frac{1}{8}[1 - \cos 2\pi u]$$

$$+ \frac{1}{2}[\gamma + \log \pi]\left(u - \frac{1}{2\pi} \sin 2\pi u\right) - \frac{1}{2}[\log(2\pi u) + \gamma - Ci(2\pi u)]$$

**Proof**

Integrating (6.2) gives us

(6.11) $$\int_0^u \psi(x) \sin^2 \pi x \, dx = -\frac{1}{4\pi} \sum_{n=2}^{\infty} \frac{1}{n} \log\left(1 - \frac{1}{n^2}\right) \sin 2n\pi u - \frac{1}{2} u \log 2 - \frac{1}{8}[1 - \cos 2\pi u]$$

$$- \frac{1}{2}\left(u - \frac{1}{2\pi} \sin 2\pi u\right)[\gamma + \log \pi]$$

We have
$$\int_0^u \psi(x) \sin^2 \pi x \, dx = \int_0^u \psi(x)[1 - \cos^2 \pi x] \, dx$$

$$= \int_0^u \psi(1+x)[1 - \cos^2 \pi x] \, dx - \int_0^u \frac{1 - \cos^2 \pi x}{x} \, dx$$

and thus
$$\int_0^u \psi(1+x) \cos^2 \pi x \, dx = \log \Gamma(1+u) - \int_0^u \psi(x) \sin^2 \pi x \, dx - \int_0^u \frac{1 - \cos^2 \pi x}{x} \, dx$$

We obtain using integration by parts

$$\int_\varepsilon^u \frac{1-\cos^2 \pi x}{x} dx = \log x(1-\cos^2 \pi x)\Big|_\varepsilon^u - \pi \int_\varepsilon^u \sin(2\pi x) \log x \, dx$$

and

$$\int_\varepsilon^u \sin(2\pi x) \log x \, dx = -\frac{1}{2\pi} \log x \cos 2\pi x \Big|_\varepsilon^u + \frac{1}{2\pi} \int_\varepsilon^u \frac{\cos(2\pi x)}{x} dx$$

$$= \frac{1}{2\pi} \log x(1-\cos 2\pi x)\Big|_\varepsilon^u + \frac{1}{2\pi} \int_\varepsilon^u \frac{\cos(2\pi x)-1}{x} dx$$

Hence we have

$$\int_0^u \sin(2\pi x) \log x \, dx = \frac{1}{2\pi} \log u(1-\cos 2\pi u) + \frac{1}{2\pi} \int_0^u \frac{\cos(2\pi x)-1}{x} dx$$

Reference to (3.2) then shows that

$$\int_0^u \sin(2\pi x) \log x \, dx = \frac{1}{2\pi}[\log u(1-\cos 2\pi u) + Ci(2\pi u) - \gamma - \log(2\pi u)]$$

and thus

(6.12) $$\int_0^u \sin(2\pi x) \log x \, dx = \frac{1}{2\pi}[Ci(2\pi u) - \log u \cos 2\pi u - \gamma - \log(2\pi)]$$

Hence, we have

$$\int_0^u \frac{1-\cos^2 \pi x}{x} dx = \log u(1-\cos^2 \pi u) - \frac{1}{2}[Ci(2\pi u) - \log u \cos 2\pi u - \gamma - \log(2\pi)]$$

$$= \frac{1}{2}\log u - \frac{1}{2}[Ci(2\pi u) - \gamma - \log(2\pi)]$$

or

(6.13) $$\int_0^u \frac{1-\cos^2 \pi x}{x} dx = \frac{1}{2}[\gamma + \log(2\pi u) - Ci(2\pi u)]$$

Therefore we obtain

$$\int_0^u \psi(1+x) \cos^2 \pi x \, dx$$

$$= \log \Gamma(1+u) + \frac{1}{4\pi} \sum_{n=2}^\infty \frac{1}{n} \log\left(1-\frac{1}{n^2}\right) \sin 2n\pi u + \frac{1}{2} u \log 2 + \frac{1}{8}[1-\cos 2\pi u]$$

$$+ \frac{1}{2}\left(u - \frac{1}{2\pi}\sin 2\pi u\right)[\gamma + \log \pi] - \frac{1}{2}[\log(2\pi u) + \gamma - Ci(2\pi u)]$$

For example, we have

(6.14) $\int_0^{\frac{1}{2}} \psi(1+x)\cos^2 \pi x\, dx = \frac{1}{4}[\log \pi - \gamma + 2Ci(\pi) - 3\log 2 + 1]$

as numerically verified by *WolframAlpha* which reports that

$$\int_0^{\frac{1}{2}} \psi(1+x)\cos^2 \pi x\, dx = -0.09114787415977...$$

and $Ci(\pi) = 0.073667912046425...$

With $u = 1$ in (6.11) we obtain

(6.14.1) $\int_0^1 \psi(x)\sin^2 \pi x\, dx = -\frac{1}{2}[\gamma + \log 2\pi]$

which concurs with [48, p.204].

With $u = \frac{1}{2}$ we have

(6.14.2) $\int_0^{\frac{1}{2}} \psi(x)\sin^2 \pi x\, dx = -\frac{1}{4}[2 + \gamma + \log \pi]$

Other integral identities may be obtained using Legendre's duplication formula for the digamma function [53, p.7].

□

It is obvious from the above that

$$\int_0^u \psi(1+x)\sin^2 \pi x\, dx = \log \Gamma(1+u) - \int_0^u \psi(1+x)\cos^2 \pi x\, dx$$

$$= -\frac{1}{4\pi}\sum_{n=2}^\infty \frac{1}{n}\log\left(1 - \frac{1}{n^2}\right)\sin 2n\pi u - \frac{1}{2}u\log 2 - \frac{1}{8}[1 - \cos 2\pi u]$$

$$-\frac{1}{2}\left(u - \frac{1}{2\pi}\sin 2\pi u\right)[\gamma + \log \pi] + \frac{1}{2}[\log(2\pi u) + \gamma - Ci(2\pi u)]$$

□

We see from (6.10) that

$$\int_0^1 \psi(1+x)\cos^2 \pi x\, dx = \frac{1}{2}Ci(2\pi)$$

Using

$$\int_0^u \psi(1+x)\cos^2 \pi x\, dx = \frac{1}{2}\log \Gamma(1+u) + \frac{1}{2}\int_0^u \psi(1+x)\cos 2\pi x\, dx$$

we have

$$\int_0^1 \psi(1+x)\cos 2\pi x\, dx = Ci(2\pi)$$

in accordance with the known result [20]

$$\int_0^1 \psi(1+x)\cos 2n\pi x\, dx = Ci(2n\pi)$$

$\square$

We note from (3.2) that

$$Ci(x) = \gamma + \log(x) + \int_0^x \frac{\cos t - 1}{t}\, dt$$

$$= \gamma + \log(x) + \int_0^1 \frac{\cos xt - 1}{t}\, dt$$

and hence

$$Ci(2\pi) = \gamma + \log(2\pi) + \int_0^1 \frac{\cos 2\pi t - 1}{t}\, dt$$

Referring to

$$\int_0^1 \psi(1+x)\cos^2 \pi x\, dx = \frac{1}{2}Ci(2\pi)$$

we see that

$$\int_0^1 \psi(1+x)\cos^2 \pi x\, dx - \frac{1}{2}\int_0^1 \frac{\cos 2\pi x - 1}{x}\, dx = \frac{1}{2}[\gamma + \log(2\pi)]$$

or equivalently

$$\int_0^1 \psi(1+x)\cos^2 \pi x\, dx + \int_0^1 \frac{\sin^2 \pi x}{x}\, dx = \frac{1}{2}[\gamma + \log(2\pi)]$$

$$\int_0^1 \left[\psi(x)\cos^2 \pi x + \frac{1}{x}\right] dx = \frac{1}{2}[\gamma + \log(2\pi)]$$

We have

$$\int_0^1 \left[\psi(x)\cos^2 \pi x + \psi(x)\sin^2 \pi x + \frac{1}{x}\right]dx = \int_0^1 \left[\psi(x) + \frac{1}{x}\right]dx$$

$$= \int_0^1 \psi(1+x)\,dx = 0$$

which results in the known integral

$$\int_0^1 \psi(x)\sin^2 \pi x\,dx = -\frac{1}{2}[\gamma + \log(2\pi)]$$

□

We see that

$$\int_0^u \psi(1+x)\cos^2 \pi x\,dx + \int_0^u \psi(x)\sin^2 \pi x\,dx = \int_0^u \psi(1+x)[\cos^2 \pi x + \sin^2 \pi x]\,dx - \int_0^u \frac{\sin^2 \pi x}{x}\,dx$$

$$= \log \Gamma(1+u) - \int_0^u \frac{\sin^2 \pi x}{x}\,dx$$

Using (6.10) and (6.11) we have

$$\int_0^u \psi(1+x)\cos^2 \pi x\,dx + \int_0^u \psi(x)\sin^2 \pi x\,dx = \log \Gamma(1+u) - \frac{1}{2}[\log(2\pi u) + \gamma - Ci(2\pi u)]$$

so that

(6.15) $$\int_0^u \frac{\sin^2 \pi x}{x}\,dx = \frac{1}{2}[\gamma + \log(2\pi u) - Ci(2\pi u)]$$

this shows that $\int_0^\infty \frac{\sin^2 \pi x}{x}\,dx$ does not converge. The above integral may of course be more directly evaluated by using integration by parts, or simply by noting that

$$Ci(x) = \gamma + \log(x) + \int_0^x \frac{\cos t - 1}{t}\,dt$$

so that

$$\gamma + \log(2\pi u) - Ci(2\pi u) = -\int_0^{2\pi u} \frac{\cos t - 1}{t}\,dt$$

and hence

$$\gamma + \log(2\pi u) - Ci(2\pi u) = -\int_0^u \frac{\cos(2\pi x) - 1}{x}\,dx$$

**Remark:**

The integral

$$\int_0^1 \psi(1+x)\cos^2 \pi x\, dx = \frac{1}{2}Ci(2\pi)$$

may be easily derived in the following manner. We see that

$$\int_0^1 \psi(x)\sin^2 \pi x\, dx = \int_0^1 \psi(x)[1-\cos^2 \pi x]\, dx$$

$$= \int_0^1 \psi(1+x)[1-\cos^2 \pi x]\, dx - \int_0^1 \frac{1-\cos^2 \pi x}{x}\, dx$$

Thus we have

$$\int_0^1 \psi(x)\sin^2 \pi x\, dx = -\int_0^1 \psi(1+x)\cos^2 \pi x\, dx - \int_0^1 \frac{1-\cos^2 \pi x}{x}\, dx$$

We obtain using integration by parts

$$\int_\varepsilon^u \frac{1-\cos^2 \pi x}{x}\, dx = \log x(1-\cos^2 \pi x)\Big|_\varepsilon^u - \pi\int_\varepsilon^u \sin(2\pi x)\log x\, dx$$

and

$$\int_\varepsilon^u \sin(2\pi x)\log x\, dx = -\frac{1}{2\pi}\log x\cos 2\pi x\Big|_\varepsilon^u + \frac{1}{2\pi}\int_\varepsilon^u \frac{\cos(2\pi x)}{x}\, dx$$

$$= \frac{1}{2\pi}\log x(1-\cos 2\pi x\Big|_\varepsilon^u + \frac{1}{2\pi}\int_\varepsilon^u \frac{\cos(2\pi x)-1}{x}\, dx$$

Reference to (3.2) then shows that

$$\int_0^u \sin(2\pi x)\log x\, dx = \frac{1}{2\pi}[Ci(2\pi u) - \log u\cos 2\pi u - \gamma - \log(2\pi)]$$

and in particular

$$\int_0^1 \sin(2\pi x)\log x\, dx = \frac{1}{2\pi}[Ci(2\pi) - \gamma - \log(2\pi)]$$

Hence, we have

$$\int_0^u \frac{1-\cos^2 \pi x}{x}\, dx = \log u(1-\cos^2 \pi u) - \frac{1}{2}[Ci(2\pi u) - \log u\cos 2\pi u - \gamma - \log(2\pi)]$$

or

$$\int_0^u \frac{1-\cos^2 \pi x}{x}\,dx = \frac{1}{2}[\log(2\pi u)+\gamma - Ci(2\pi u)]$$

In particular we have

$$\int_0^1 \frac{1-\cos^2 \pi x}{x}\,dx = \frac{1}{2}[\log(2\pi)+\gamma - Ci(2\pi)]$$

We see that

$$\int_0^1 \psi(1+x)\cos^2 \pi x\,dx = -\int_0^1 \psi(x)\sin^2 \pi x\,dx - \int_0^1 \frac{1-\cos^2 \pi x}{x}\,dx$$

$$= \frac{1}{2}[\gamma + \log 2\pi] - \frac{1}{2}[\log(2\pi)+\gamma - Ci(2\pi)]$$

Hence we obtain

$$\int_0^1 \psi(1+x)\cos^2 \pi x\,dx = \frac{1}{2} Ci(2\pi)$$

$\square$

Integration by parts gives us

$$\int_0^u \psi(1+x)\cos^2 \pi x\,dx = \log \Gamma(1+u)\cos^2 \pi u + \pi \int_0^u \log \Gamma(1+x)\sin 2\pi x\,dx$$

and in particular we have

$$\int_0^1 \psi(1+x)\cos^2 \pi x\,dx = \pi \int_0^1 \log \Gamma(1+x)\sin 2\pi x\,dx$$

$$= \pi \int_0^1 \log \Gamma(x)\sin 2k\pi x\,dx + \pi \int_0^1 \sin(2\pi x)\log x\,dx$$

Therefore we obtain as before

$$\int_0^1 \log \Gamma(x)\sin 2\pi x\,dx = \frac{\gamma + \log 2\pi}{2\pi}$$

**Proposition 6.2**

We have

(6.16) $$\int_0^1 x\log \Gamma(x)\sin 2\pi x\,dx = \frac{\gamma}{4\pi}$$

giving us an interesting integral involving both $\gamma$ and $\pi$.

**Proof**

We multiply (6.2) by $x$ and integrate over the interval $[0,1]$ to obtain

$$0 = \frac{1}{2}\log 2 - \frac{1}{4} + \frac{1}{2}[\log \pi + \gamma] + \int_0^1 x(1 - \cos 2\pi x)\psi(x)\,dx$$

We have

$$\int_0^1 x\psi(x)\,dx = x\log\Gamma(x)\Big|_0^1 - \int_0^1 \log\Gamma(x)\,dx$$

$$= -\frac{1}{2}\log 2\pi$$

and

$$\int_0^1 x\cos 2\pi x\,\psi(x)\,dx = x\cos 2\pi x\log\Gamma(x)\Big|_0^1 - \int_0^1 [\cos 2\pi x - 2\pi x\sin 2\pi x]\log\Gamma(x)\,dx$$

$$= -\int_0^1 \log\Gamma(x)\cos 2\pi x\,dx + 2\pi\int_0^1 x\log\Gamma(x)\sin 2\pi x\,dx$$

Employing the known Fourier coefficients [48, p.203]

(6.17) $\quad \int_0^1 \log\Gamma(x)\cos 2k\pi x\,dx = \frac{1}{4k}$

we obtain (6.16) which was previously determined in a different manner by Mező [45].

Mező [45] also showed that

(6.18) $\quad \int_0^1 x\log\Gamma(x)\sin 2\pi x\,dx = \frac{1}{8} - \frac{\gamma + \log 2\pi}{2\pi^2} - \frac{1}{2\pi^2}\sum_{n=1}^{\infty} \frac{\log n}{n^2 - 1}$

**Proposition 6.3**

(6.19) $\quad \sum_{n=2}^{\infty} \psi\left(n + \frac{1}{2}\right)\log\left(1 - \frac{1}{n^2}\right) = [\gamma + 2\log 2 - 2]\log 2 - 4\sum_{n=1}^{\infty}\frac{\log n}{4n^2 - 1}$

**Proof**

We note from (6.3) that $\sum_{n=2}^{\infty}\log\left(1 - \frac{1}{n^2}\right) = -\log 2$ and we may therefore write (6.2) as

$$(6.20) \quad \sum_{n=2}^{\infty} \log\left(1-\frac{1}{n^2}\right)(1-\cos 2n\pi x) = \frac{\pi}{2}\sin 2\pi x + (1-\cos 2\pi x)[\log \pi + \gamma + \psi(x)]$$

or equivalently as

$$2\sum_{n=2}^{\infty} \log\left(1-\frac{1}{n^2}\right)\sin^2 n\pi x = \frac{\pi}{2}\sin 2\pi x + 2\sin^2 \pi x[\log \pi + \gamma + \psi(x)]$$

and

$$\sum_{n=2}^{\infty} \frac{\sin^2 n\pi x}{\sin^2 \pi x} \log\left(1-\frac{1}{n^2}\right) = \frac{\pi}{2}\cot \pi x + \log \pi + \gamma + \psi(x)$$

We multiply (6.20) by $\dfrac{\sin \pi x}{1-\cos 2\pi x}$ to obtain

$$\sum_{n=2}^{\infty} \frac{(1-\cos 2n\pi x)\sin \pi x}{1-\cos 2\pi x} \log\left(1-\frac{1}{n^2}\right) = \frac{\pi}{2}\cos \pi x + [\log \pi + \gamma]\sin \pi x + \psi(x)\sin \pi x$$

since $\dfrac{\sin 2\pi x \sin \pi x}{1-\cos 2\pi x} = \cos \pi x$.

The following indefinite integral may be readily verified by differentiation

$$\int \frac{(1-\cos 2n\pi x)\sin \pi x}{1-\cos 2\pi x} dx = -\sum_{j=0}^{n-1} \frac{\cos(2j+1)\pi x}{(2j+1)\pi}$$

and thus we obtain the definite integral

$$\int_0^1 \frac{(1-\cos 2n\pi x)\sin \pi x}{1-\cos 2\pi x} dx = \frac{2}{\pi}\sum_{j=0}^{n-1}\frac{1}{2j+1}$$

We note [53, p.20]

$$2\sum_{j=0}^{n-1}\frac{1}{2j+1} = \psi\left(n+\frac{1}{2}\right) + \gamma + 2\log 2$$

and obtain

$$\frac{1}{\pi}\sum_{n=2}^{\infty} \psi\left(n+\frac{1}{2}\right)\log\left(1-\frac{1}{n^2}\right) - \frac{1}{\pi}[\gamma + 2\log 2]\log 2 = \frac{2}{\pi}[\log \pi + \gamma] + \int_0^1 \psi(x)\sin \pi x\, dx$$

where we have again employed (6.3). Thus, we see that

$$(6.21) \quad \sum_{n=2}^{\infty} \psi\left(n+\frac{1}{2}\right)\log\left(1-\frac{1}{n^2}\right) = 2\log \pi + 2\gamma + \gamma \log 2 + 2\log^2 2 + \pi\int_0^1 \psi(x)\sin \pi x\, dx$$

Kölbig [40] has shown that

$$(6.22) \quad \int_0^1 \psi(x)\sin \pi x\, dx = -\frac{2}{\pi}\left[\log(2\pi)+\gamma+2\sum_{n=1}^{\infty}\frac{\log n}{4n^2-1}\right]$$

resulting in

$$(6.23) \quad \sum_{n=2}^{\infty}\psi\left(n+\frac{1}{2}\right)\log\left(1-\frac{1}{n^2}\right) = [\gamma+2\log 2-2]\log 2 - 4\sum_{n=1}^{\infty}\frac{\log n}{4n^2-1}$$

*WolframAlpha* computes the finite sum

$$\sum_{n=2}^{10{,}000}\psi\left(n+\frac{1}{2}\right)\log\left(1-\frac{1}{n^2}\right)+4\frac{\log n}{4n^2-1} \simeq -0.02529\ldots$$

which closely approximates the analytic result $[\gamma+2\log 2-2]\log 2$.

We also note that

$$\sum_{j=0}^{n-1}\frac{1}{2j+1} = H_{2n}-\frac{1}{2}H_n$$

$$\sum_{j=1}^{2n}\frac{(-1)^n}{n} = H_{2n}-H_n$$

and thus

$$\sum_{j=0}^{n-1}\frac{1}{2j+1} = \frac{1}{2}H_n+\sum_{j=1}^{2n}\frac{(-1)^n}{n}$$

$$\frac{1-\cos 2n\pi x}{1-\cos 2\pi x} = \frac{\sin^2 n\pi x}{\sin^2 \pi x} = n+2n\sum_{k=1}^{n-1}\left(1-\frac{k}{n}\right)\cos 2\pi x$$

$\square$

**Proposition 6.4**

$$(6.24) \quad \int_0^1 x(1-x)\cos \pi x \cot \pi x\, dx = \frac{1}{\pi^3}[7\varsigma(3)-4]$$

**Proof**

We multiply (6.20) by $\dfrac{x(1-x)\cos \pi x}{1-\cos 2\pi x}$ to obtain

(6.25) $$\sum_{n=2}^{\infty} \frac{x(1-x)(1-\cos 2n\pi x)\cos \pi x}{1-\cos 2\pi x} \log\left(1-\frac{1}{n^2}\right)$$

$$= \frac{\pi}{2} x(1-x)\cos \pi x \cot \pi x + [\log \pi + \gamma] x(1-x)\cos \pi x + x(1-x)\psi(x)\cos \pi x$$

Letting
$$f(x) = \frac{x(1-x)(1-\cos 2n\pi x)\cos \pi x}{1-\cos 2\pi x}$$

we see that $f(x) = -f(1-x)$ and hence $\int_0^1 f(x)dx = 0$.

Therefore, we have

$$\int_0^1 \frac{x(1-x)(1-\cos 2n\pi x)\cos \pi x}{1-\cos 2\pi x} dx = 0$$

and we note that

$$\int_0^1 x(1-x)\cos \pi x \, dx = 0$$

Integrating (6.25) accordingly results in

(6.26) $$\frac{\pi}{2} \int_0^1 x(1-x)\cos \pi x \cot \pi x \, dx = -\int_0^1 x(1-x)\psi(x)\cos \pi x \, dx$$

and we shall see this again in (6.36).

We see from (4.18) that

$$\int_0^1 x(1-x)\cos \pi x \cot \pi x \, dx = 2\sum_{n=1}^{\infty} \int_0^1 x(1-x)\cos \pi x \sin 2n\pi x \, dx$$

and we have

$$2\int_0^1 x(1-x)\cos \pi x \sin 2n\pi x \, dx = \int_0^1 x(1-x)[\sin(2n+1)\pi x + \sin(2n-1)\pi x] dx$$

$$= \frac{16n(4n^2+3)}{\pi^3(4n^2-1)^3}$$

Hence we obtain

$$\int_0^1 x(1-x)\cos \pi x \cot \pi x \, dx = \frac{16}{\pi^3} \sum_{n=1}^{\infty} \frac{n(4n^2+3)}{(4n^2-1)^3}$$

Prima facie, we would therefore like to obtain expressions for the two series

$$\sum_{n=1}^{\infty} \frac{n^3}{(4n^2-1)^3} \text{ and } \sum_{n=1}^{\infty} \frac{n}{(4n^2-1)^3}.$$

Alternatively, we see that

$$\int_0^1 x(1-x)\cos\pi x \cot\pi x\, dx = \frac{16}{\pi^3} \sum_{n=1}^{\infty} \frac{n(4n^2-1+4)}{(4n^2-1)^3}$$

$$= \frac{16}{\pi^3} \sum_{n=1}^{\infty} \frac{n}{(4n^2-1)^2} + \frac{64}{\pi^3} \sum_{n=1}^{\infty} \frac{n}{(4n^2-1)^3}$$

We would therefore like to obtain expressions for the series $\sum_{n=1}^{\infty} \frac{n}{(4n^2-1)^2}$ and

$$\sum_{n=1}^{\infty} \frac{n}{(4n^2-1)^3}.$$

In this connection we recall (2.8)

$$\psi(1+x) + \psi(1-x) + 2\gamma = -2x^2 \sum_{n=1}^{\infty} \frac{1}{n} \frac{1}{n^2-x^2}$$

and differentiation results in

(6.27) $$\psi'(1+x) - \psi'(1-x) = -4x^2 \sum_{n=1}^{\infty} \frac{1}{n} \frac{1}{(n^2-x^2)^2} - 4x \sum_{n=1}^{\infty} \frac{1}{n} \frac{1}{n^2-x^2}$$

which simplifies to [19]

(6.28) $$\psi'(1+x) - \psi'(1-x) = -4x \sum_{n=1}^{\infty} \frac{n}{(n^2-x^2)^2}$$

The second derivative gives us

(6.29) $$\psi''(1+x) + \psi''(1-x) = -16x^2 \sum_{n=1}^{\infty} \frac{n}{(n^2-x^2)^3} - 4 \sum_{n=1}^{\infty} \frac{n}{(n^2-x^2)^2}$$

We may also note that (6.28) and (6.29) may be expressed as

(6.30) $$\psi'(1+x) - \psi'(1-x) - 2\frac{\psi(1+x)+\psi(1-x)+2\gamma}{x} = -4x^2 \sum_{n=1}^{\infty} \frac{1}{n} \frac{1}{(n^2-x^2)^2}$$

(6.31) $$\psi''(1+x)+\psi''(1-x)-\frac{\psi'(1+x)-\psi'(1-x)}{x}=-16x^2\sum_{n=1}^{\infty}\frac{n}{(n^2-x^2)^3}$$

Therefore, reference to (6.30) shows that with $x=\frac{1}{2}$

$$\psi'(\tfrac{3}{2})-\psi'(\tfrac{1}{2})=-32\sum_{n=1}^{\infty}\frac{n}{(4n^2-1)^2}$$

It is well known that [53, p.22]

$$\psi^{(n)}(z)=(-1)^{n+1}n!\varsigma(n+1,z)$$

$$\psi^{(n)}(1+z)=\psi^{(n)}(z)+(-1)^n n!z^{-(n+1)}$$

and thus we have

$$\sum_{n=1}^{\infty}\frac{n}{(4n^2-1)^2}=\frac{1}{8}$$

Reference to (6.31) shows that with $x=\frac{1}{2}$

$$\psi''(\tfrac{3}{2})+\psi''(\tfrac{1}{2})-2[\psi'(\tfrac{3}{2})-\psi'(\tfrac{1}{2})]=-256\sum_{n=1}^{\infty}\frac{n}{(4n^2-1)^3}$$

$$\sum_{n=1}^{\infty}\frac{n}{(4n^2-1)^3}=-\frac{1}{128}[12+\psi''(\tfrac{1}{2})]$$

$$=-\frac{1}{128}[12-2\varsigma(3,\tfrac{1}{2})]$$

Hence we obtain

$$\sum_{n=1}^{\infty}\frac{n}{(4n^2-1)^3}=\frac{1}{64}[7\varsigma(3)-6]$$

and conclude with

$$\int_0^1 x(1-x)\cos\pi x\cot\pi x\,dx=\frac{1}{\pi^3}[7\varsigma(3)-4]$$

From (6.31) we see that

$$\sum_{n=1}^{\infty}\frac{n}{(n^2-x^2)^3}=\frac{1}{16x^2}\left[\frac{\psi'(1+x)-\psi'(1-x)}{x}-\psi''(1+x)+\psi''(1-x)\right]$$

$$:= \varphi(x)$$

Using the relationship

$$\sum_{n=1}^{\infty}(-1)^n a_n = 2\sum_{n=1}^{\infty} a_{2n} - \sum_{n=1}^{\infty} a_n$$

we have

$$\sum_{n=1}^{\infty}\frac{(-1)^n n}{(n^2-x^2)^3} = 2\sum_{n=1}^{\infty}\frac{2n}{(4n^2-x^2)^3} - \sum_{n=1}^{\infty}\frac{n}{(n^2-x^2)^3}$$

We see that

$$2\sum_{n=1}^{\infty}\frac{2n}{(4n^2-x^2)^3} = \frac{1}{16}\sum_{n=1}^{\infty}\frac{n}{(n^2-(x/2)^2)^3}$$

$$= \frac{1}{16}\varphi(x/2)$$

Hence we obtain

(6.32) $$\sum_{n=1}^{\infty}\frac{(-1)^n n}{(n^2-x^2)^3} = \frac{1}{16}\varphi(x/2) - \varphi(x)$$

For example, *WolframAlpha* provides the summation

(6.33) $$\sum_{n=1}^{\infty}\frac{(-1)^n n}{(4n^2-1)^3} = \frac{1}{512}\left[24G + \varsigma\left(2,\frac{5}{4}\right) - 32 - \pi^2\right]$$

where $G$ is Catalan's constant. Catalan's constant features here because

$$\psi'(\tfrac{1}{4}) = \pi^2 + 8G$$

**Proposition 6.5**

(6.34) $$\int_0^1 \psi(x)x(1-x)\cos\pi x\,dx = \frac{1}{\pi^2}\left[2 - \frac{7}{2}\varsigma(3)\right]$$

**Proof**

We showed [25] that for $0 < x < 1$

(6.35) $$\psi(x) = \lim_{s\to 1}\sum_{n=1}^{\infty}\frac{2[\gamma + \log(2\pi n)]\cos(2n\pi x) - \pi\sin(2n\pi x)}{n^{1-s}}$$

Assuming that the interchange of integration and summation is valid, we have

$$\int_0^1 \psi(x)x(1-x)\cos \pi x\, dx = \lim_{s \to 1}\sum_{n=1}^{\infty} \frac{2[\gamma + \log(2\pi n)]I_n - \pi J_n}{n^{1-s}}$$

where

$$I_n = \int_0^1 x(1-x)\cos \pi x \cos(2n\pi x)\, dx$$

and

$$J_n = \int_0^1 x(1-x)\cos \pi x \sin(2n\pi x)\, dx$$

With the substitution $x = 1-t$ in $I_n = \int_0^1 x(1-x)\cos \pi x \cos(2n\pi x)\, dx$ we easily see that $I_n = 0$ for all $n \geq 0$

With regard to $J_n$, *WolframAlpha* conveniently provides us with the indefinite integral

$$\int x(1-x)\cos(\pi x)\sin(2n\pi x)\, dx =$$
$$\frac{(\pi^2(1-2n)^2 x^2 - 2)\cos(\pi(2n-1)x) + 2\pi(1-2n)x\sin(\pi(2n-1)x)}{2\pi^3(2n-1)^3} -$$
$$\frac{\pi(2n-1)x\cos(\pi(2n-1)x) - \sin(\pi(2n-1)x)}{2(\pi - 2\pi n)^2} +$$
$$\frac{\sin(\pi(2n+1)x) - \pi(2n+1)x\cos(\pi(2n+1)x)}{2(2\pi n + \pi)^2} -$$
$$\frac{2\pi(2n+1)x\sin(\pi(2n+1)x) - (\pi^2(2nx+x)^2 - 2)\cos(\pi(2n+1)x)}{2(2\pi n + \pi)^3}$$
$$+ \text{constant}$$

which may be verified by differentiating the given output. We then have

$$J_n = \int_0^1 x(1-x)\cos \pi x \sin(2n\pi x)\, dx = \frac{2}{\pi^3}\left[\frac{1}{(2n-1)^3} + \frac{1}{(2n+1)^3}\right]$$

and thus

$$\int_0^1 \psi(x)x(1-x)\cos \pi x\, dx = -\frac{2}{\pi^2}\sum_{n=1}^{\infty}\left[\frac{1}{(2n-1)^3} + \frac{1}{(2n+1)^3}\right]$$

We see that

$$\sum_{n=1}^{\infty}\frac{1}{(2n-1)^3} = \sum_{n=0}^{\infty}\frac{1}{(2n+1)^3}$$

$$= \frac{1}{8} \sum_{n=0}^{\infty} \frac{1}{(n+\frac{1}{2})^3}$$

$$= \frac{1}{8} \varsigma(3, \tfrac{1}{2})$$

Noting that

$$\sum_{n=1}^{\infty} \frac{1}{(2n+1)^3} = -1 + \sum_{n=0}^{\infty} \frac{1}{(2n+1)^3}$$

we deduce

$$\sum_{n=1}^{\infty} \left[ \frac{1}{(2n-1)^3} + \frac{1}{(2n+1)^3} \right] = -1 + \frac{1}{4} \varsigma(3, \tfrac{1}{2})$$

Since $\varsigma(s, \tfrac{1}{2}) = (2^s - 1)\varsigma(s)$ we obtain

$$\int_0^1 \psi(x) x(1-x) \cos \pi x \, dx = \frac{1}{\pi^2} \left[ 2 - \frac{7}{2} \varsigma(3) \right]$$

as was originally determined by Glasser [33].

$\square$

An alternative derivation of (6.34) is shown below.

We have

$$I = \int_0^1 \psi(x) x(1-x) \cos \pi x \, dx$$

and with $x = 1 - t$ we find that

$$I = -\int_0^1 \psi(1-x) x(1-x) \cos \pi x \, dx$$

Using
$$\psi(x) - \psi(1-x) = -\pi \cot \pi x$$

we obtain

(6.36) $$2I = -\pi \int_0^1 x(1-x) \cos \pi x \cot \pi x \, dx$$

We then employ (4.18) to write

$$\int_0^1 x(1-x) \cos \pi x \cot \pi x \, dx = 2 \sum_{n=1}^{\infty} \int_0^1 x(1-x) \cos \pi x \sin 2n\pi x \, dx$$

Having regard to the definition, we see that

$$\int_0^1 x(1-x)\cos \pi x \cot \pi x\, dx = 2\sum_{n=1}^\infty J_n$$

and the required result easily follows.

□

Glasser [33] expressed the integral as

$$\int_0^1 \cot(\pi x) x(1-x)\cos \pi x\, dx = \int_0^1 x(1-x)\frac{1-\sin^2 \pi x}{\sin \pi x}\, dx$$

$$= \int_0^1 x(1-x)\left[\frac{1}{\sin \pi x} - \sin \pi x\right] dx$$

and noting that

$$\int_0^1 \frac{x(1-x)}{\sin \pi x}\, dx = \int_0^1 \frac{x}{\sin \pi x}\, dx - \int_0^1 \frac{x^2}{\sin \pi x}\, dx$$

he then employed the known integrals

$$\int_0^1 \frac{x}{\sin x}\, dx = 2G \qquad \int_0^1 \frac{x^2}{\sin x}\, dx = 2\pi G - \frac{7}{2}\varsigma(3)$$

where $G$ is Catalan's constant.

□

Letting $x \to 1-x$ in (6.34) gives us

$$\int_0^1 \psi(1-x)x(1-x)\cos \pi x\, dx = -\frac{1}{\pi^2}\left[2 - \frac{7}{2}\varsigma(3)\right]$$

We then have

$$\int_0^1 [\psi(x)-\psi(1-x)]x(1-x)\cos \pi x\, dx = \frac{2}{\pi^2}\left[2 - \frac{7}{2}\varsigma(3)\right]$$

and using $\psi(x)-\psi(1-x) = -\pi \cot \pi x$ we see that

(6.37) $$\int_0^1 x(1-x)\cos \pi x \cot \pi x\, dx = \frac{2}{\pi^3}\left[\frac{7}{2}\varsigma(3) - 2\right]$$

**Proposition 6.6**

(6.38) $$\int_0^1 x(1-x)\cos \pi x\, \psi(x/2)\, dx = -2\log A - \frac{1}{\pi^2}\left[2 + \frac{7}{2}\varsigma(3)\right] + \frac{1}{6}[\gamma + \log \pi]$$

$$-\frac{1}{2\pi^2}\sum_{n=2}^{\infty}\frac{1}{n^2}\log\left(1-\frac{1}{4n^2}\right)$$

**Proof**

Glasser [33] noted that the integral $\int_0^1 x(1-x)\cos\pi x\,\psi(x/2)\,dx$ did not appear to be expressible in a similar form as (6.34).

Assuming that the interchange of integration and summation is valid, using (6.35) we have

$$\int_0^1 \psi(x/2)x(1-x)\cos\pi x\,dx = \lim_{s\to 1}\sum_{n=1}^{\infty}\frac{2[\gamma+\log(2\pi n)]K_n - \pi L_n}{n^{1-s}}$$

where

$$K_n = \int_0^1 x(1-x)\cos\pi x\cos(n\pi x)\,dx$$

and

$$L_n = \int_0^1 x(1-x)\cos\pi x\sin(n\pi x)\,dx$$

Since $K_{2n+1}$ does not vanish, it is clear that the answer in this case will contain the second derivative of the Riemann zeta function.

There is an alternative approach. We set $x \to x/2$ in (6.2), multiply by $x(1-x)$ and then integrate over [0,1]. We thereby obtain

(6.39)
$$\int_0^1 x(1-x)(1-\cos\pi x)\psi(x/2)\,dx = \frac{1}{\pi^2}\sum_{n=2}^{\infty}\frac{1+(-1)^n}{n^2}\log\left(1-\frac{1}{n^2}\right) - \frac{1}{6}[\gamma+\log(2\pi)] + \frac{2}{\pi^2}$$

$$= \frac{1}{2\pi^2}\sum_{n=2}^{\infty}\frac{1}{n^2}\log\left(1-\frac{1}{4n^2}\right) - \frac{1}{6}[\gamma+\log(2\pi)] + \frac{2}{\pi^2}$$

where we have employed the elementary integrals

$$\int_0^1 x(1-x)\cos n\pi x\,dx = -\frac{1+(-1)^n}{\pi^2 n^2}$$

$$\int_0^1 x(1-x)\cos\pi x\,dx = 0$$

$$\int_0^1 x(1-x)\sin \pi x\, dx = \frac{4}{\pi^3}$$

Integrals of the form $\int_0^z x\psi(x+a)\,dx$ and $\int_0^z x^2 \psi(x+a)\,dx$ are given by Srivastava and Choi in [53, p.207 & 209] in terms of the Barnes multiple gamma functions and, using these, we may deduce that

(6.40) $$\int_0^1 x(1-x)\psi(x/2)\,dx = -2\log A - \frac{7\varsigma(3)}{2\pi^2} - \frac{1}{6}\log 2$$

which may also be obtained using *WolframAlpha*. We therefore obtain Glasser's elusive integral

$$\int_0^1 x(1-x)\cos \pi x\, \psi(x/2)\,dx = -2\log A - \frac{1}{\pi^2}\left[2 + \frac{7}{2}\varsigma(3)\right] + \frac{1}{6}[\gamma + \log \pi]$$

$$-\frac{1}{2\pi^2}\sum_{n=2}^{\infty}\frac{1}{n^2}\log\left(1 - \frac{1}{4n^2}\right)$$

## 7. Some connections with the generalized Stieltjes Constants

The digamma function may be regarded as one of the generalized Stieltjes constants because $\psi(x) = -\gamma_0(x)$.

We recall Lerch's trigonometric series expansion for the digamma function for $0 < x < 1$ (see for example Gronwall's paper [37, p.105] and Nielsen's book [48, p.204])

(7.1) $$\psi(x)\sin \pi x + \frac{\pi}{2}\cos \pi x + (\gamma + \log 2\pi)\sin \pi x = -\sum_{n=1}^{\infty}\sin(2n+1)\pi x \cdot \log\left(1 + \frac{1}{n}\right)$$

It may be noted that Boyack [9] has given an extremely elementary derivation of (6.2) simply by multiplying (7.1) by $2\sin \pi x$. We obtain

$$2\psi(x)\sin^2 \pi x = -\frac{\pi}{2}\sin 2\pi x - 2(\gamma + \log 2\pi)\sin^2 \pi x - 2\sin \pi x \sum_{n=1}^{\infty}\sin(2n+1)\pi x \cdot \log\left(1 + \frac{1}{n}\right)$$

Using $2\sin(2n+1)\pi x \sin \pi x = \cos 2n\pi x - \cos(2n+2)\pi x$ we have

$$2\sin \pi x \sum_{n=1}^{\infty}\sin(2n+1)\pi x \cdot \log\left(1 + \frac{1}{n}\right) = \sum_{n=1}^{\infty}\cos 2n\pi x \cdot \log\left(1 + \frac{1}{n}\right)$$

$$-\sum_{n=1}^{\infty}\cos 2(n+1)\pi x\cdot\log\left(1+\frac{1}{n}\right)$$

$$=\log 2.\cos 2\pi x+\sum_{n=2}^{\infty}\cos 2n\pi x\cdot\log\left(1+\frac{1}{n}\right)-\sum_{m=2}^{\infty}\cos 2m\pi x\cdot\log\left(1+\frac{1}{m-1}\right)$$

$$=\log 2.\cos 2\pi x+\sum_{n=2}^{\infty}\cos 2n\pi x\cdot\log\left(\frac{n+1}{n}\frac{n-1}{n}\right)$$

$$=\log 2.\cos 2\pi x+\sum_{n=2}^{\infty}\cos 2n\pi x\cdot\log\left(1-\frac{1}{n^2}\right)$$

This results in

$$2\psi(x)\sin^2\pi x=-\frac{\pi}{2}\sin 2\pi x-2(\gamma+\log 2\pi)\sin^2\pi x-\log 2.\cos 2\pi x-\sum_{n=2}^{\infty}\cos 2n\pi x\cdot\log\left(1-\frac{1}{n^2}\right)$$

and using (6.3) $\sum_{n=2}^{\infty}\log\left(1-\frac{1}{n^2}\right)=-\log 2$ we obtain (6.2)

$$-\sum_{n=2}^{\infty}\log\left(1-\frac{1}{n^2}\right)\cos 2n\pi x=\log 2+\frac{\pi}{2}\sin 2\pi x+(1-\cos 2\pi x)[\log\pi+\gamma+\psi(x)]$$

We also see that

(7.2) $\quad \sum_{n=2}^{\infty}\log\left(1-\frac{1}{n^2}\right)\cos 2n\pi x=-\log 2.\cos 2\pi x+2\sin\pi x\sum_{n=1}^{\infty}\sin(2n+1)\pi x\cdot\log\left(1+\frac{1}{n}\right)$

which we may write as

$$=-\log 2.\cos 4\pi x+\sum_{n=2}^{\infty}[\cos 2n\pi x-\cos(2n+2)\pi x]\cdot\log\left(1+\frac{1}{n}\right)$$

Therefore we see that

$$\sum_{n=2}^{\infty}\log\left(1-\frac{1}{n}\right)\cos 2n\pi x=-\log 2.\cos 4\pi x-\sum_{n=2}^{\infty}\cos 2(n+1)\pi x\cdot\log\left(1+\frac{1}{n}\right)$$

□

Letting $x=\frac{1}{2}$ in (7.2) shows that

(7.3) $\quad \sum_{n=2}^{\infty}(-1)^n\log\left(1-\frac{1}{n^2}\right)=\log 2+2\sum_{n=1}^{\infty}(-1)^n\log\left(1+\frac{1}{n}\right)$

which concurs with (6.5) and (6.6).

□

Letting $x \to x+\tfrac{1}{2}$ in (7.1) results in

(7.4)
$$\psi(x+\tfrac{1}{2})\cos \pi x - \frac{\pi}{2}\sin \pi x + (\gamma + \log 2\pi)\cos \pi x = -\sum_{n=1}^{\infty}(-1)^n \cos(2n+1)\pi x \cdot \log\left(1+\frac{1}{n}\right)$$

**Proposition 7.1**

(7.5)
$$\int_0^u [\gamma_1(1-x)-\gamma_1(x)]\sin \pi x\, dx = [\gamma + \log(2\pi)]\sin \pi u + \sum_{n=1}^{\infty}\frac{1}{2n+1}\log\left(1+\frac{1}{n}\right)\sin(2n+1)\pi u$$

(7.6)
$$\sum_{n=1}^{\infty}\frac{(-1)^n}{2n+1}\log\left(1+\frac{1}{n}\right) = -4\sum_{n=1}^{\infty}\frac{(-1)^n n \log n}{4n^2-1}$$

**Proof**

It was shown in [25] that

(7.7)
$$\sum_{n=1}^{\infty}\log\left(1+\frac{1}{n}\right)\cos(2n+1)\pi x = \frac{1}{\pi}[\gamma_1(1-x)-\gamma_1(x)]\sin \pi x - [\gamma + \log(2\pi)]\cos \pi x$$

and integration gives us (7.5).

We showed in [25] that for $0 < x < 1$

(7.8)
$$\gamma_1(x,s) - \gamma_1(1-x,s) = -2\pi \sum_{n=1}^{\infty}\frac{[\gamma + \log(2\pi n)]}{n^s}\sin(2n\pi x)$$

where
$$\gamma_1(x) = \lim_{s \to 0}\gamma_1(x,s)$$

We formally have using (7.8)

(7.9)
$$\int_0^u [\gamma_1(1-x)-\gamma_1(x)]\sin \pi x\, dx = 2\pi \lim_{s \to 0}\sum_{n=1}^{\infty}\int_0^u \frac{[\gamma + \log(2\pi n)]}{n^s}\sin(2n\pi x)\sin \pi x\, dx$$

and we have the elementary integrals

$$\int \sin(2n\pi x)\sin \pi x\, dx = \frac{\sin((2n-1)\pi x)}{2\pi(2n-1)} - \frac{\sin((2n+1)\pi x)}{2\pi(2n+1)} + c$$

and

$$\int_0^u \sin(2n\pi x)\sin \pi x\, dx = \frac{\sin((2n-1)\pi u)}{2\pi(2n-1)} - \frac{\sin((2n+1)\pi u)}{2\pi(2n+1)}$$

Letting $u = \frac{1}{2}$ we see that

$$\int_0^{\frac{1}{2}} \sin(2n\pi x)\sin \pi x\, dx = -\frac{2(-1)^n n}{(4n^2-1)\pi}$$

and thus

$$\int_0^{\frac{1}{2}} [\gamma_1(1-x) - \gamma_1(x)]\sin \pi x\, dx = -4\sum_{n=1}^{\infty} \frac{(-1)^n n[\gamma + \log(2\pi n)]}{4n^2-1}$$

$$= -4[\gamma + \log(2\pi)]\sum_{n=1}^{\infty} \frac{(-1)^n n}{4n^2-1} - 4\sum_{n=1}^{\infty} \frac{(-1)^n n \log n}{4n^2-1}$$

*WolframAlpha* kindly informs us that

$$\sum_{n=1}^{\infty} \frac{(-1)^n n}{4n^2-1} = -\frac{1}{4}$$

and hence we obtain

(7.10) $$\int_0^{\frac{1}{2}} [\gamma_1(1-x) - \gamma_1(x)]\sin \pi x\, dx = \gamma + \log(2\pi) - 4\sum_{n=1}^{\infty} \frac{(-1)^n n \log n}{4n^2-1}$$

□

We may write (7.9) as follows

$$\int_0^u [\gamma_1(1-x) - \gamma_1(x)]\sin \pi x\, dx = \sum_{n=1}^{\infty} [\gamma + \log(2\pi n)]\left[\frac{\sin((2n-1)\pi u)}{2n-1} - \frac{\sin((2n+1)\pi u)}{2n+1}\right]$$

$$= [\gamma + \log(2\pi)]\sum_{n=1}^{\infty} \left[\frac{\sin((2n-1)\pi u)}{2n-1} - \frac{\sin((2n+1)\pi u)}{2n+1}\right]$$

$$+ \sum_{n=1}^{\infty} \left[\frac{\sin((2n-1)\pi u)}{2n-1} - \frac{\sin((2n+1)\pi u)}{2n+1}\right]\log n$$

We have the Fourier series [54, p.149] for $0 < u < 1$

$$\sum_{n=0}^{\infty} \frac{\sin((2n+1)\pi u)}{2n+1} = \frac{\pi}{4}$$

and therefore

$$\sum_{n=1}^{\infty} \frac{\sin\left((2n+1)\pi u\right)}{2n+1} = \frac{\pi}{4} - \sin \pi u$$

We see that

$$\sum_{n=1}^{\infty} \frac{\sin\left((2n-1)\pi u\right)}{2n-1} = \sum_{m=0}^{\infty} \frac{\sin\left((2m+1)\pi u\right)}{2m+1}$$

and we therefore obtain

(7.10.1) $$\sum_{n=1}^{\infty} \left[ \frac{\sin\left((2n-1)\pi u\right)}{2n-1} - \frac{\sin\left((2n+1)\pi u\right)}{2n+1} \right] = \sin \pi u$$

We also see that

$$\sum_{n=1}^{\infty} \frac{\sin\left((2n-1)\pi u\right)}{2n-1} \log n = \sum_{m=0}^{\infty} \frac{\sin\left((2m+1)\pi u\right)}{2m+1} \log(m+1)$$

$$= \sum_{n=1}^{\infty} \frac{\sin\left((2n+1)\pi u\right)}{2n+1} \log(n+1)$$

Therefore, we obtain

$$\sum_{n=1}^{\infty} \left[ \frac{\sin\left((2n-1)\pi u\right)}{2n-1} - \frac{\sin\left((2n+1)\pi u\right)}{2n+1} \right] \log n = \sum_{n=1}^{\infty} \frac{\sin\left((2n+1)\pi u\right)}{2n+1} \log\left(1 + \frac{1}{n}\right)$$

and we have thereby come full circle to (7.5).

**Remark:**

Letting $u = \frac{1}{2}$ in (7.4) we see that

(7.11) $$\sum_{n=1}^{\infty} \frac{(-1)^n}{2n+1} \log\left(1 + \frac{1}{n}\right) = -4 \sum_{n=1}^{\infty} \frac{(-1)^n n \log n}{4n^2 - 1}$$

We showed in [27] that

(7.12) $$\sum_{n=1}^{\infty} \frac{1}{2n+1} \log\left(1 + \frac{1}{n}\right) = 2 \sum_{n=1}^{\infty} \frac{\log n}{4n^2 - 1}$$

Curiously, *WolframAlpha* evaluates the combined series (7.11) as follows

$$\sum_{n=1}^{\infty}\left(\frac{(-1)^n \log(1+\frac{1}{n})}{2n+1} + \frac{(4(-1)^n)n\log(n)}{4n^2-1}\right) \approx 0.0132252 - 4.55291 \times 10^{-14}\,i$$

and it was initially an upset to find that the real part did not vanish. Undoubtedly an experimentalist would probably regard the imaginary output as evanescent.

On the other hand, *WolframAlpha* does give identical results when the series are evaluated separately:

$$\sum_{n=1}^{\infty}\frac{(-1)^n \log(1+\frac{1}{n})}{2n+1} = -0.176012$$

$$\sum_{n=1}^{\infty}\frac{(4(-1)^n)n\log(n)}{4n^2-1} = 0.176012$$

We see from (1.20) that

$$\sum_{n=1}^{\infty}\frac{\log n}{4n^2-1} = -\sum_{m=1}^{\infty}\frac{\varsigma'(2m)}{2^{2(m+1)}}$$

**Proposition 7.2**

(7.13) $$\int_0^1 \psi(x)\sin^2 \pi x\, dx = -\frac{1}{2}[\gamma + \log(2\pi)]$$

**Proof**

We showed [25] that for $0 < x < 1$

(7.14) $$\psi(x) = \lim_{s \to 1}\sum_{n=1}^{\infty}\frac{2[\gamma + \log(2\pi n)]\cos(2n\pi x) - \pi\sin(2n\pi x)}{n^{1-s}}$$

Assuming that we may interchange the order of integration and summation we tentatively obtain

$$\int_0^1 \psi(x)\sin^2 \pi x\, dx = \lim_{s \to 1}\sum_{n=1}^{\infty}\int_0^1 \frac{2[\gamma + \log(2\pi n)]\cos(2n\pi x) - \pi\sin(2n\pi x)}{n^{1-s}}\sin^2 \pi x\, dx$$

and immediately deduce that (see also (6.14.1))

$$\int_0^1 \psi(x)\sin^2 \pi x\, dx = -\frac{1}{2}[\gamma + \log(2\pi)]$$

where we have employed the elementary integrals

$$\int_0^1 \cos(2\pi x)\sin^2(\pi x)\,dx = -\frac{1}{4}$$

$$\int_0^1 \cos(2n\pi x)\sin^2(\pi x)\,dx = 0 \quad n \neq 1$$

$$\int_0^1 \sin(2n\pi x)\sin^2(\pi x)\,dx = 0 \quad n \geq 0$$

**Proposition 7.3**

We have for $0 \leq u \leq 1$

(7.15) $\displaystyle\int_0^u \psi(x)\sin\pi x\,dx = \frac{2}{\pi}\sum_{n=1}^\infty \left[\frac{2n\sin(2n\pi u)}{4n^2-1}\sin\pi u + \frac{\cos(2n\pi u)}{4n^2-1}[\cos\pi u - 1]\right]\log n$

$$+2[\gamma + \log(2\pi)]\left[\frac{1}{4}\sin\pi u + \frac{1}{2\pi}[\cos\pi u - 1]\right] - \frac{1}{2}\sin\pi u$$

**Proof**

We have using (7.14)

$$\int_0^u \psi(x)\sin\pi x\,dx = \lim_{s\to 1}\sum_{n=1}^\infty \int_0^u \frac{2[\gamma+\log(2\pi n)]\cos(2n\pi x) - \pi\sin(2n\pi x)}{n^{1-s}}\sin\pi x\,dx$$

We have the elementary integrals

$$\int \cos(2n\pi x)\sin\pi x\,dx = \frac{2n\sin(2n\pi x)}{\pi(4n^2-1)}\sin\pi x + \frac{\cos(2n\pi x)}{\pi(4n^2-1)}\cos\pi x + c$$

$$\int \sin(2n\pi x)\sin\pi x\,dx = \frac{\sin(2n-1)\pi x}{2\pi(2n-1)} - \frac{\sin(2n+1)\pi x}{2\pi(2n+1)} + c$$

and thus

$$\int_0^u \cos(2n\pi x)\sin\pi x\,dx = \frac{2n\sin(2n\pi u)}{\pi(4n^2-1)}\sin\pi u + \frac{\cos(2n\pi u)}{\pi(4n^2-1)}[\cos\pi u - 1]$$

$$\int_0^u \sin(2n\pi x)\sin\pi x\,dx = \frac{\sin(2n-1)\pi u}{2\pi(2n-1)} - \frac{\sin(2n+1)\pi u}{2\pi(2n+1)}$$

This results in

$$\int_0^u \psi(x)\sin\pi x\,dx = 2\sum_{n=1}^\infty [\gamma + \log(2\pi n)]\left[\frac{2n\sin(2n\pi u)}{\pi(4n^2-1)}\sin\pi u + \frac{\cos(2n\pi u)}{\pi(4n^2-1)}[\cos\pi u - 1]\right]$$

$$-\frac{1}{2}\sum_{n=1}^{\infty}\left[\frac{\sin\left((2n-1)\pi u\right)}{2n-1}-\frac{\sin\left((2n+1)\pi u\right)}{2n+1}\right]$$

We note the Fourier series [5, p.337] for $0<u<1$

$$\cos \pi u = \frac{8}{\pi}\sum_{n=1}^{\infty}\frac{n\sin(2n\pi u)}{4n^2-1}$$

$$\sin \pi u = \frac{2}{\pi}-\frac{4}{\pi}\sum_{n=1}^{\infty}\frac{\cos(2n\pi u)}{4n^2-1}$$

and, since $\sum_{n=1}^{\infty}\frac{1}{4n^2-1}=\frac{1}{2}$, the latter holds true for $0\leq u \leq \pi$. Therefore, we see that

$$\sum_{n=1}^{\infty}\frac{\cos 2n\pi u - 1}{4n^2-1} = -\frac{\pi}{4}\sin \pi u$$

Accordingly, using (7.10.1) which is valid for $0 \leq u \leq 1$

$$\sum_{n=1}^{\infty}\left[\frac{\sin\left((2n-1)\pi u\right)}{2n-1}-\frac{\sin\left((2n+1)\pi u\right)}{2n+1}\right] = \sin \pi u$$

we have for $0<u<1$

$$\int_0^u \psi(x)\sin \pi x\, dx = \frac{2}{\pi}\sum_{n=1}^{\infty}\left[\frac{2n\sin(2n\pi u)}{4n^2-1}\sin \pi u + \frac{\cos(2n\pi u)}{4n^2-1}[\cos \pi u - 1]\right]\log n$$

$$+2[\gamma+\log(2\pi)]\left[\frac{1}{4}\cos \pi u \sin \pi u + \left(\frac{1}{2\pi}-\frac{1}{4}\sin \pi u\right)[\cos \pi u - 1]\right] - \frac{1}{2}\sin \pi u$$

which simplifies to

$$\int_0^u \psi(x)\sin \pi x\, dx = \frac{2}{\pi}\sum_{n=1}^{\infty}\left[\frac{2n\sin(2n\pi u)}{4n^2-1}\sin \pi u + \frac{\cos(2n\pi u)}{4n^2-1}[\cos \pi u - 1]\right]\log n$$

$$+2[\gamma+\log(2\pi)]\left[\frac{1}{4}\sin \pi u + \frac{1}{2\pi}[\cos \pi u - 1]\right] - \frac{1}{2}\sin \pi u$$

As mentioned above, Kölbig [40] showed that

(7.16) $$\int_0^1 \psi(x)\sin \pi x\, dx = -\frac{2}{\pi}\left[\log(2\pi)+\gamma+2\sum_{n=1}^{\infty}\frac{\log n}{4n^2-1}\right]$$

Kölbig [40] states that the infinite series in (7.16) "does not seem to be expressible in terms of well-known functions". A different proof of (7.16) is given in [con 20].

We therefore see that (7.15) actually also holds for $0 \leq u \leq 1$.

With $u = \tfrac{1}{2}$ we obtain

$$(7.17) \quad \int_0^{\tfrac{1}{2}} \psi(x)\sin \pi x\, dx = -\frac{2}{\pi}\sum_{n=1}^{\infty}\frac{(-1)^n \log n}{4n^2-1} + [\gamma+\log(2\pi)]\left[\frac{1}{2}-\frac{1}{\pi}\right] - \frac{1}{2}$$

**A different attempt at the integral** $\int_0^1 e^{-px}\psi(1+x)\, dx$

We recall (1.12)

$$\int_0^1 e^{-px}\psi(1+x)\, dx = \frac{1}{2}(1-e^{-p})[\gamma+\log(2\pi)] + \frac{1-e^{-p}}{4}\left[\psi\left(1+\frac{ip}{2\pi}\right) + \psi\left(1-\frac{ip}{2\pi}\right)\right]$$

$$-[\gamma+\log p - Ei(-p)] + 2p(1-e^{-p})\sum_{n=1}^{\infty}\frac{\gamma+\log(2\pi n)}{4\pi^2 n^2 + p^2}$$

We showed [25] that for $0 < x < 1$

$$\psi(x) = \lim_{s\to 1}\sum_{n=1}^{\infty}\frac{2[\gamma+\log(2\pi n)]\cos(2n\pi x) - \pi \sin(2n\pi x)}{n^{1-s}}$$

Assuming that the interchange of integration and summation is valid, we have

$$\int_\varepsilon^1 e^{-px}\psi(1+x)\, dx = \int_\varepsilon^1 \frac{e^{-px}}{x}\, dx + \lim_{s\to 1}\sum_{n=1}^{\infty}\int_\varepsilon^1\left[\frac{2[\gamma+\log(2\pi n)]\cos(2n\pi x) - \pi\sin(2n\pi x)}{n^{1-s}}\right]e^{-px}\, dx$$

$$= 2\sum_{n=1}^{\infty}\int_\varepsilon^1 [\gamma+\log(2\pi n)]\cos(2n\pi x)e^{-px}\, dx$$

$$+\int_\varepsilon^1 \frac{e^{-px}}{x}\, dx - \pi\sum_{n=1}^{\infty}\int_\varepsilon^1 \sin(2n\pi x)e^{-px}\, dx$$

We consider the first integral (where we are not concerned with convergence issues and, ab initio, we let $\varepsilon = 0$)

$$2\sum_{n=1}^{\infty}\int_0^1 [\gamma+\log(2\pi n)]\cos(2n\pi x)e^{-px}\, dx = 2(1-e^{-p})p\sum_{n=1}^{\infty}\frac{\gamma+\log(2\pi n)}{4\pi^2 n^2 + p^2}$$

and we note that this term appears in (1.12). The rest of the attempted evaluation is more subtle.

We consider the finite series

$$I(N,\varepsilon) = \int_{\varepsilon}^{1} \frac{e^{-px}}{x} dx - \pi \sum_{n=1}^{N} \int_{\varepsilon}^{1} \sin(2n\pi x) e^{-px} dx$$

$$= \left[\int_{\varepsilon}^{1} \frac{e^{-px}}{x} dx - \frac{1-e^{-p}}{2} \log N\right] + \left[\frac{1-e^{-p}}{2} \log N - \frac{1}{2} \sum_{n=1}^{N} \int_{\varepsilon}^{1} 2\pi \sin(2n\pi x) e^{-px} dx\right]$$

and noting the Frullani integral [48, p.128]

$$\int_{0}^{\infty} \frac{e^{-x} - e^{-ax}}{x} dx = \log a$$

we see that

$$\log N = \int_{0}^{1} \frac{x^{N-1} - 1}{\log x} dx$$

which we express as

$$= \int_{\varepsilon}^{1} \frac{x^{N-1} - 1}{\log x} dx + \eta(\varepsilon) \text{ where } \eta(\varepsilon) \to 0 \text{ as } \varepsilon \to 0$$

We then have

$$I(N,\varepsilon) = \int_{\varepsilon}^{1} \left[\frac{e^{-px}}{x} - \frac{1-e^{-p}}{2} \frac{x^{N-1}-1}{\log x}\right] dx - \eta(\varepsilon) + \frac{1-e^{-p}}{2} \log N - \frac{1}{2} \sum_{n=1}^{N} \int_{\varepsilon}^{1} 2\pi \sin(2n\pi x) e^{-px} dx$$

We have the indefinite integral

$$\int \sin(2n\pi x) e^{-px} dx = -\frac{e^{-px}(p \sin(2n\pi x) + 2n\pi \cos(2n\pi x))}{4\pi^2 n^2 + p^2} + c$$

which results in

$$\int_{0}^{1} e^{-px} \sin(2n\pi x) dx = \frac{2n\pi(1-e^{-p})}{4\pi^2 n^2 + p^2}$$

Letting $(N,\varepsilon) \to (\infty, 0)$ we have

$$\frac{1-e^{-p}}{2} \log N - \frac{1}{2} \sum_{n=1}^{N} \int_{\varepsilon}^{1} 2\pi \sin(2n\pi x) e^{-px} dx \to -\frac{1-e^{-p}}{2} \lim_{N \to \infty} \left[\sum_{n=1}^{N} \frac{4\pi^2 n}{4\pi^2 n^2 + p^2} - \log N\right]$$

$$= -\frac{1-e^{-p}}{2}\Lambda\left(\frac{p}{2\pi}\right)$$

$$= \frac{1-e^{-p}}{4}\left[\psi\left(1+\frac{ip}{2\pi}\right)+\psi\left(1-\frac{ip}{2\pi}\right)\right]$$

where we have used (1.2) and (1.5). We note that this term also appears in (1.12).

At this stage, it appears that we have obtained

$$\int_0^1 e^{-px}\psi(1+x)\,dx = \lim_{(N,\varepsilon)\to(\infty,0)}\int_\varepsilon^1\left[\frac{e^{-px}}{x} - \frac{1-e^{-p}}{2}\frac{x^{N-1}-1}{\log x}\right]dx$$

$$+\frac{1-e^{-p}}{4}\left[\psi\left(1+\frac{ip}{2\pi}\right)+\psi\left(1-\frac{ip}{2\pi}\right)\right]+2p(1-e^{-p})\sum_{n=1}^\infty \frac{\gamma+\log(2\pi n)}{4\pi^2 n^2 + p^2}$$

which we would like to hope is equivalent to (1.12)

$$\int_0^1 e^{-px}\psi(1+x)\,dx = \frac{1}{2}(1-e^{-p})[\gamma+\log(2\pi)]-[\gamma+\log p - Ei(-p)]$$

$$+\frac{1-e^{-p}}{4}\left[\psi\left(1+\frac{ip}{2\pi}\right)+\psi\left(1-\frac{ip}{2\pi}\right)\right]+2p(1-e^{-p})\sum_{n=1}^\infty \frac{\gamma+\log(2\pi n)}{4\pi^2 n^2 + p^2}$$

We have

$$\int\left[\frac{e^{-px}}{x} - \frac{1-e^{-p}}{2}\frac{x^{N-1}-1}{\log x}\right]dx = \frac{1-e^{-p}}{2}[li(x)-Ei(N\log x)]+Ei(-px)$$

where $li(x)$ is the logarithmic integral defined by

$$li(x) = \int_0^x \frac{dt}{\log t}$$

It may be noted that $li(0) = 0$ and $li(x) = Ei(\log x)$. We are "close, but no cigar."

## 8. Miscellaneous results

It was shown in [17] that

(8.1) $$\frac{1}{2}\int_a^b p(x)\,dx = \sum_{n=0}^\infty \int_a^b p(x)\cos\alpha nx\,dx$$

where $p(x)$ is assumed to be twice continuously differentiable on $[a,b]$. Equation (8.1) is valid provided (i) $\sin(\alpha x/2) \neq 0 \ \forall \ x \in [a,b]$ or, alternatively, (ii) if $\sin(\alpha \eta / 2) = 0$ where $\eta \in [a,b]$ then $p(\eta) = 0$ also.

We proved in [21] that the above requirement for $p(x)$ to be a twice continuously differentiable function could be considerably relaxed in certain prescribed circumstances.

For example, we select $p(x) = \sin \mu x$ to give us

$$\frac{1}{2}\int_0^t \sin \mu x \, dx = \sum_{n=0}^{\infty} \int_0^t \sin \mu x \cos \alpha n x \, dx$$

where we require $t < 2\pi/\alpha$. We could have $t = 2\pi/\alpha$ if $\mu = \alpha/2$.

This results in

(8.2) $$\frac{1}{2\mu}(\cos \mu t - 1) = \sum_{n=1}^{\infty} \frac{\alpha n \sin \mu t \sin \alpha n t + \mu \cos \mu t \cos \alpha n t - \mu}{\alpha^2 n^2 - \mu^2}$$

Letting $\alpha = 1$ in (8.2) results in

(8.3) $$\frac{1}{2\mu}(\cos \mu t - 1) = \sum_{n=1}^{\infty} \frac{n \sin \mu t \sin n t + \mu \cos \mu t \cos n t - \mu}{n^2 - \mu^2}$$

and with $t = \pi$ we have

(8.4) $$\frac{1}{2\mu}(\cos \mu \pi - 1) = \mu \sum_{n=1}^{\infty} \frac{(-1)^n \cos \mu \pi - 1}{n^2 - \mu^2}$$

or equivalently

(8.5) $$\frac{\cos \mu \pi}{\mu} = 2\mu \cos \mu \pi \sum_{n=1}^{\infty} \frac{(-1)^n}{n^2 - \mu^2} + \frac{1}{\mu} - \sum_{n=1}^{\infty} \frac{2\mu}{n^2 - \mu^2}$$

Using the cotangent decomposition formula

(8.6) $$\pi \cot(\pi x) = \frac{1}{x} - \sum_{n=1}^{\infty} \frac{2x}{n^2 - x^2}$$

we obtain

$$\frac{\cos \mu \pi}{\mu} = 2\mu \cos \mu \pi \sum_{n=1}^{\infty} \frac{(-1)^n}{n^2 - \mu^2} + \pi \cot(\mu \pi)$$

and dividing this by $\cos\mu\pi$ we see that

(8.7) $$\frac{\pi}{\sin(\mu\pi)} = \frac{1}{\mu} - 2\mu\sum_{n=1}^{\infty}\frac{(-1)^n}{n^2-\mu^2}$$

The following simple trigonometric identities are easily proved

(8.8) $$\cot(x/2) - \cot x = \frac{2}{\sin x}$$

(8.9) $$\cot(x/2) + \tan(x/2) = \frac{2}{\sin x}$$

(8.10) $$\cot(x/2) - \tan(x/2) = 2\cot x$$

Combining (8.6) and (8.8) we see that

$$\pi\cot(\pi x/2) - 2\pi\cot(\pi x) = \frac{2}{x} - 4x\sum_{n=1}^{\infty}\frac{1}{4n^2-x^2} - \frac{1}{x} + 2x\sum_{n=1}^{\infty}\frac{1}{n^2-x^2}$$

$$= \frac{1}{x} - 4x\sum_{n=1}^{\infty}\frac{1}{4n^2-x^2} + 2x\sum_{n=1}^{\infty}\frac{1}{n^2-x^2}$$

and using $\sum_{n=1}^{\infty}(-1)^n a_n = \sum_{n=1}^{\infty}[2a_{2n}-a_n]$ we are able to immediately deduce (8.7).

$\square$

Starting the summation of (8.4) at $n=2$ produces

$$\frac{1}{2\mu}(\cos\mu\pi - 1) = \mu\sum_{n=2}^{\infty}\frac{(-1)^n\cos\mu\pi - 1}{n^2-\mu^2} - \mu\frac{\cos\mu\pi + 1}{1-\mu^2}$$

and employing L'Hôpital's rule we see that

$$\lim_{\mu\to 1}\frac{\cos\mu\pi + 1}{1-\mu^2} = 0$$

and hence we obtain

$$\sum_{n=2}^{\infty}\frac{(-1)^n + 1}{n^2-1} = 1$$

or equivalently

(8.11) $$\sum_{n=1}^{\infty} \frac{1}{4n^2-1} = \frac{1}{2}$$

□

Letting $\alpha t = 2\pi$ and $\mu = \alpha/2$ in (8.2) results in

$$\frac{1}{2\mu}(\cos \mu t - 1) = \sum_{n=1}^{\infty} \frac{\mu \cos \mu t - \mu}{\alpha^2 n^2 - \mu^2}$$

and, with $\mu = \alpha/2$, we obtain (8.11) again.

□

We let $p(x) = \cos \mu x$ in (8.1) with $\alpha = 2\pi$

$$-\frac{1}{2}\int_{\frac{1}{2}}^{t} \cos \mu x \, dx = \sum_{n=1}^{\infty} \int_{\frac{1}{2}}^{t} \cos \mu x \cos 2\pi n x \, dx$$

$$\int \cos \mu x \cos 2\pi n x \, dx = -\frac{\mu \sin \mu x \cos 2\pi n x - 2\pi n \cos \mu x \sin 2\pi n x}{4\pi^2 n^2 - \mu^2}$$

$$\int_{\frac{1}{2}}^{t} \cos \mu x \cos 2\pi n x \, dx = -\frac{\mu[\sin \mu t \cos 2\pi n t - (-1)^n \sin(\mu/2)] - 2\pi n \cos \mu t \sin 2\pi n t}{4\pi^2 n^2 - \mu^2}$$

$$\frac{1}{2\mu}[\sin \mu t - \sin(\mu/2)] = \sum_{n=1}^{\infty} \frac{\mu[\sin \mu t \cos 2\pi n t - (-1)^n \sin(\mu/2)] - 2\pi n \cos \mu t \sin 2\pi n t}{4\pi^2 n^2 - \mu^2}$$

With $\mu = \pi$ we obtain

(8.12) $$\frac{1}{2}[\sin \pi t - 1] = \sum_{n=1}^{\infty} \frac{\sin \pi t \cos 2\pi n t - (-1)^n - 2n \cos \pi t \sin 2\pi n t}{4n^2 - 1}$$

It is an exercise in Apostol's book [5, p.337] to show that for $0 < x < \pi$

(8.13) $$\sin \pi t = \frac{2}{\pi} - \frac{4}{\pi} \sum_{n=1}^{\infty} \frac{\cos(2\pi n t)}{4n^2 - 1}$$

(8.14) $$\frac{\pi}{8} \cos \pi t = \sum_{n=1}^{\infty} \frac{n \sin 2\pi n t}{4n^2 - 1}$$

Letting $t = \frac{1}{2}$ in (8.13) gives us

(8.15) $$1 = \frac{2}{\pi} - \frac{4}{\pi} \sum_{n=1}^{\infty} \frac{(-1)^n}{4n^2 - 1}$$

and simple algebra shows that this is consistent with (8.12).

## 9. Open access to our own work

This paper contains references to various other papers and, rather surprisingly, most of them are currently freely available on the internet. Surely now is the time that all of our work should be freely accessible by all. The mathematics community should lead the way on this by publishing everything on arXiv, or in an equivalent open access repository. We think it, we write it, so why hide it? You know it makes sense.

## 10. Acknowledgement

I thank Iaroslav Blagouchine for providing me with copies of the Russian and English versions of Salaev's paper [50] in October 2019.

## REFERENCES

[1]   M. Abramowitz and I.A. Stegun (Eds.), Handbook of Mathematical Functions with Formulas, Graphs and Mathematical Tables. Dover, New York, 1970.
http://www.math.sfu.ca/~cbm/aands/

[2]   V. Adamchik and H.M. Srivastava, Some series of the zeta and related functions. Analysis, 18(1998), 131--144.
https://www.researchgate.net/publication/228581846_Some_series_of_the_zeta_and_related_functions

[3]   T. Amdeberhan, O. Espinosa and V.H. Moll, The Laplace transform of the digamma function: An integral due to Glasser, Manna and Oloa.
Proc. Amer. Math. Soc., 136 (2008) 3211–3221.
http://129.81.170.14/~vhm/papers_html/official-oloa.pdf
http://arxiv.org/pdf/0707.3663.pdf

[4]   T.M. Apostol, Introduction to Analytic Number Theory.
Springer-Verlag, New York, Heidelberg and Berlin, 1976.

[5]   T.M. Apostol, Mathematical Analysis, Second Ed., Addison-Wesley Publishing Company, Menlo Park (California), London and Don Mills (Ontario), 1974.

[6]   J.L.F. Bertrand, Traité de Calcul Différentiel et de Calcul Intégral.
Gauthier-Villars, Paris, 1864.
http://gallica.bnf.fr/ark:/12148/bpt6k99558p
http://gallica.bnf.fr/ark:/12148/bpt6k995591

[7]   G. Boros and V.H. Moll, Irresistible Integrals: Symbolics, Analysis and Experiments in the Evaluation of Integrals. Cambridge University Press, 2004.

[8]   M.T. Boudjelkha, A proof that extends Hurwitz formula into the critical strip. Applied Mathematics Letters, 14 (2001) 309-403.

[9]  R. Boyack, Summation of certain trigonometric series with logarithmic
     coefficients, 2021.
      https://arxiv.org/pdf/2109.08686.pdf

[10] K. Boyadzhiev and R. Frontczak, Series Involving Euler's Eta (or Dirichlet Eta)
      Function. Journal of Integer Sequences, Vol. 24 (2021).
        https://cs.uwaterloo.ca/journals/JIS/VOL24/Frontczak/front22.pdf

[11] D.M. Bradley, Ramanujan's formula for the logarithmic derivative of the
      gamma function.
       Mathematical Proceedings of the Cambridge Philosophical Society, Vol. 120
      (October 1996), no. 3, pp. 391-401. MR1388195 (97a:11132)
       https://arxiv.org/pdf/math/0505125.pdf

[12] T.J.I'A. Bromwich, Introduction to the theory of infinite series.
       Third edition. AMS Chelsea Publishing, 1991.
        https://archive.org/details/introductiontoth00bromuoft/page/n13/mode/1up

[13] B. Candelpergher, Ramanujan summation of divergent series.
      Lectures notes in mathematics 2185, Springer, 2017.
       https://hal.univ-cotedazur.fr/hal-01150208v2/document

[14] H.S. Carslaw, Introduction to the theory of Fourier Series and Integrals.
       Third Ed. Dover Publications Inc, 1930.

[15] H. Cohen, Number Theory. Volume II: Analytic and modern tools.
        Springer Science, 2007.

[16] D.F. Connon, Some series and integrals involving the Riemann zeta function,
      binomial coefficients and the harmonic numbers. Volume IV, 2007.
       https://arxiv.org/abs/0710.4028

[17] D.F. Connon, Some series and integrals involving the Riemann zeta function,
      binomial coefficients and the harmonic numbers. Volume V, 2007.
       https://arxiv.org/abs/0710.4047

[18] D.F. Connon, Some series and integrals involving the Riemann zeta function,
      binomial coefficients and the harmonic numbers. Volume VI, 2007.
       https://arxiv.org/abs/0710.4032

[19] D.F. Connon, New rapidly converging series representations for values of the
      Riemann zeta function and the Dirichlet beta function. 2010.
       https://arxiv.org/abs/1003.4592

[20] D.F. Connon, Some trigonometric integrals involving the log gamma and the
      digamma function. 2010.
       https://arxiv.org/abs/1005.3469

[21] D.F. Connon, Some applications of the Dirichlet integrals to the summation of series and the evaluation of integrals involving the Riemann zeta function. 2012.
https://arxiv.org/abs/1212.0441

[22] D.F. Connon, On an integral involving the digamma function. 2012.
https://arxiv.org/abs/1212.1432

[23] D.F. Connon, Some integrals and series involving the Stieltjes constants, 2018.
https://arxiv.org/abs/1801.05711

[24] D.F. Connon, An introduction to the Barnes double gamma function with an application to an integral involving the cotangent function, 2018.
https://arxiv.org/abs/1801.08025

[25] D.F. Connon, A Ramanujan enigma involving the first Stieltjes constant, 2019.
https://arxiv.org/abs/1901.03382

[26] D.F. Connon, Some new formulae involving the Stieltjes constants, 2019.
https://arxiv.org/abs/1902.00510

[27] D.F. Connon, A new representation of the Stieltjes constants, 2022.
http://arxiv.org/abs/2201.05084

[28] A. Dixit, The Laplace transform of the psi function.
Proc. Amer. Math. Soc., 138 (2010) 593–603.
http://www.ams.org/journals/proc/2010-138-02/S0002-9939-09-10157-0/S0002-9939-09-10157-0.pdf

[29] R. Dwilewicz and J. Mináč, An introduction to relations between the values of $\varsigma(s)$ in terms of holomorphic functions of two variables.
Proceedings of the Hayama Symposium on Several Complex Variables, Japan, Dec. 2000. Pages 28-38 (2001).

[30] E. Elizalde and A. Romeo, An integral involving the generalized zeta function.
Internat. J. Maths. & Maths. Sci. Vol.13, No.3, (1990) 453-460.
http://www.hindawi.com/GetArticle.aspx?doi=10.1155/S0161171290000679&e=CTA

[31] O. Espinosa and V.H. Moll, On some integrals involving the Hurwitz zeta function: Part I. The Ramanujan Journal, 6,150-188, 2002.
http://arxiv.org/abs/math.CA/0012078

[32] O. Furdui, College Math. Journal, 38, No.1, 61, 2007

[33] M.L. Glasser, Evaluation of some integrals involving the $\psi$ - function.
Math. of Comp., Vol.20, No.94, 332-333, 1966.
http://www.ams.org/journals/mcom/1966-20-094/S0025-5718-66-99934-0/S0025-5718-66-99934-0.pdf

[34] M.L. Glasser and D. Manna, On the Laplace transform of the psi function.
Tapas in Experimental Mathematics (T. Amdeberhan and V. Moll, eds.), Contemporary Mathematics, vol. 457, Amer. Math. Soc., Providence, RI, 205-214 (2008)
https://carma.newcastle.edu.au/resources/jon/Preprints/Papers/By%20Others/Dante-Glasser/ContempMath2.pdf

[35]  M. Godefroy [, p.12] La Fonction Gamma: Théorie, Histoire, Bibliographie.
Gauthier-Villars, 1901.
https://archive.org/details/lafonctiongammat29800gut/page/n19/mode/1up

[36]  I.S. Gradshteyn and I.M. Ryzhik, Tables of Integrals, Series and Products.
Eighth Ed., Academic Press, 2015.
Errata for Eighth Edition http://www.mathtable.com/errata/gr8_errata.pdf

[37]  T.H. Gronwall, The gamma function in integral calculus.
Annals of Math., 20, 35-124, 1918.

[38]  G.H. Hardy, Divergent Series. Chelsea Publishing Company, New York, 1991.

[39]  K. Knopp, Theory and Application of Infinite Series.
Second English Ed. Dover Publications Inc, New York, 1990.

[40]  K.S. Kölbig, On three trigonometric integrals of $\log \Gamma(x)$ or its derivatives.
J. Comput. Appl. Math. 54 (1994) 129-131.
CERN/Computing and Networks Division; CN/94/7, May 1994.
http://doc.cern.ch/tmp/convert_P00023218.pdf

[41]  M. Lerch, Další studie v oboru Malmsténovských řad,
Rozpravy Ceske Akad. 3, no. 28, 1894, 63 pp.
https://dml.cz/handle/10338.dmlcz/501780

[42]  M. Lerch, Sur une relation ayant rapports avec la théorie de la fonction gamma.
Bulletin int. de l'Ac. Prague 2 (1895), 214 - 218.
https://dml.cz/handle/10338.dmlcz/501482

[43]  M. Lerch, Sur la différentiation d'une classe de séries trigonométriques.
*Ann. Scientifiques de l'École Normale Supérieure*, Sér.3, tome 12 (1895),
351-361. http://www.numdam.org/

[44]  Matyáš Lerch, Über eine Formel aus der Théorie der Gammafunction,
Monatsch. Math. Phys. 8 (1897), 187–192
https://dml.cz/bitstream/handle/10338.dmlcz/501506/Lerch_01-0000-135_1.pdf

[45]  I. Mező, The Fourier series of the log-Barnes function, 2016.
http://arxiv.org/pdf/1604.00753.pdf
https://sites.google.com/site/istvanmezo81/publications

[46]  L.M. Milne-Thompson, The calculus of finite differences, Macmillan & Co Ltd,
1933.
https://archive.org/details/calculusoffinite032017mbp/page/n1/mode/2up

[47]  P. J. Nahin, Inside Interesting Integrals, Springer, 2015.

[48]  N. Nielsen, Die Gammafunktion.
Chelsea Publishing Company, Bronx and New York, 1965.

[49] O. Oloa, Some Euler-type integrals and a new rational series for Euler's constant.
Tapas in Experimental Mathematics (T. Amdeberhan and V. Moll, eds.), Contemporary Mathematics, vol. 457, Amer. Math. Soc., Providence, RI, 253-264 (2008)
https://les-mathematiques.net/vanilla/uploads/dump_data/2007/0721/14/6737

[50] B.W. Salaev, A certain class of trigonometric series.
Mat. Zametki, 9 (1971), 533-542 (Russian).

[51] J. Sondow, A faster product for $\pi$ and a new integral for $\log\frac{\pi}{2}$.
Amer. Math. Monthly 112 (2005) 729-734.
https://arxiv.org/ftp/math/papers/0401/0401406.pdf

[52] M.R. Spiegel, Schaum's Outline Series of Theory and Problems of Advanced Calculus. McGraw Hill Book Company, 1963.

[53] H.M. Srivastava and J. Choi, Series Associated with the Zeta and Related Functions. Kluwer Academic Publishers, Dordrecht, the Netherlands, 2001.

[54] G.P. Tolstov, Fourier Series. (Translated from the Russian by R.A. Silverman) Dover Publications Inc, New York, 1976.

[55] N.M. Vildanov, Generalization of the Glasser-Manna-Oloa integral and some new integrals of similar type. 2010.
https://arxiv.org/pdf/1007.3460.pdf

[56] Z.X. Wang and D.R. Guo, Special Functions.
World Scientific Publishing Co Pte Ltd, Singapore, 1989.

[57] R. Witula, D. Jama, D. Slota and P. Zwolénski, Some Fourier series expansions and their applications. International Journal of Pure and Applied Mathematics, Volume 64, No. 2, 2010, 199-223
http://www.ijpam.eu/contents/2010-64-2/4/4.pdf

[58] E.C. Titchmarsh, The Zeta-Function of Riemann.
Oxford University (Clarendon) Press, Oxford, London and New York, 1951; Second Ed. (Revised by D.R. Heath- Brown), 1986.

[59] V.S. Adamchik, A Class of Logarithmic Integrals. Proceedings of the 1997 International Symposium on Symbolic and Algebraic Computation. ACM, Academic Press, 1-8, 2001.
https://www.researchgate.net/publication/221564534_A_Class_of_Logarithmic_Integrals

[60] D.F. Connon, Fourier series representations of the logarithms of the Euler gamma function and the Barnes multiple gamma functions. 2009.
https://arxiv.org/abs/0903.4323